\documentclass[twoside,10pt]{article}
\usepackage[utf8]{inputenc}
\usepackage[french,english]{babel}
\usepackage{csquotes}
\usepackage{times}
\usepackage{authblk}
\usepackage{amsmath,amsfonts,amssymb,amsthm,esvect,envmath,euscript,mathtools,pifont}
\usepackage{bm}
\usepackage{mathrsfs}
\usepackage{soul}
\usepackage[colorlinks=true]{hyperref}
\usepackage{fullpage}
\usepackage{ulem}
\usepackage{bbm} 
\usepackage{verbatim}
\usepackage{graphicx}
\usepackage{float}
\usepackage{subcaption}
\usepackage{mismath}
\usepackage[dvipsnames]{xcolor}
\allowdisplaybreaks     
\usepackage[toc,page]{appendix}

\newcommand{\Ss}{\sum_{i=1}^N}
\newcommand{\Sd}{\sum_{i,j=1}^N}
\newcommand{\inte}[1]{\int_0^{#1}}
\newcommand{\esp}[1]{\mathbb{E}\left(#1\right)}
\newcommand{\espa}[1]{\mathbb{E}\left[#1\right]}

\newcommand{\ens}{[\![1;N]\!]}

\newcommand{\reel}{\mathbb{R}}
\newcommand{\varia}[4][]{{#1}^{{#3},N,{#4}}_{#2}}
\newcommand{\norme}[1]{\big\lVert {#1} \big\rVert}
\newcommand{\Norme}[1]{\big\lvert {#1} \big\rvert}
\newcommand{\teps}{\tau_n^\varepsilon}
\newcommand{\scalp}[2][]{\left\langle {#1}, {#2}  \right\rangle}

\newcommand{\proba}[1]{\mathbb{P}\left(#1\right)}

\newcommand{\supre}[1]{\underset{0\leq t\leq{#1}}{\sup\,}}

\newcommand{\deriv}{\di}

\newcommand{\espacecont}{C^0\left([0;T],\mathbb{R}^{2Nd}\right)}
\newcommand{\espacecontabs}{C_{\mathrm{abs}}^0\left([0;T],\mathbb{R}^{2Nd}\right)}
\newcommand{\probabi}{\mathbb{P}_{\vv{z_0}}}
\newcommand{\esper}{\mathbb{E}_{\vv{z_0}}}
\newcommand{\esperd}{\mathbb{E}_{(\vv{x},\vv{v})}}
\newcommand{\pointeq}{\left(\vv{x}_{\mathrm{eq}}^\ast,0_{\reel^{Nd}}\right)}
\newcommand{\initia}{(\vv{x_0},\vv{v_0})}
\newcommand{\initias}{(\vv{x},\vv{v})}
\newcommand{\probabis}{\mathbb{P}_{\initias}}
\newcommand{\ensinva}[1]{\mathcal{O}_N^{\gamma}\left(#1\right)}
\newcommand{\ensinvab}[1]{\partial\mathcal{O}_N^{\gamma}\left(#1\right)}
\newcommand{\ensinvad}[1]{\overline{\mathcal{O}_N^{\gamma}\left(#1\right)}}
\newcommand{\grandeva}[1]{\big(\vv{X}^{N,\varepsilon}_{#1},\vv{V}^{N,\varepsilon}_{#1}\big)}
\newcommand{\grandevae}[1]{\big(\vv{X}^{N,\varepsilon,h}_{#1},\vv{V}^{N,\varepsilon,h}_{#1}\big)}
\newcommand{\grandevad}[1]{\big(\vv{x}^{N}_{#1},\vv{v}^{N}_{#1}\big)}
\newcommand{\ta}{\tau^{\alpha,\gamma,N}_{\varepsilon}(R)}

\newtheorem{theoreme}{Theorem}[section]
\newtheorem{remarque}[theoreme]{Remark}
\newtheorem{proposition}[theoreme]{Proposition}
\newtheorem{definition}[theoreme]{Definition}
\newtheorem{lemme}[theoreme]{Lemma}
\newtheorem{corollaire}[theoreme]{Corollary}

\numberwithin{equation}{section}
\begin{document}
	\title{Freidlin-Wentzell exit-laws for a stochastic Cucker-Smale model}
\date{\today}
\author{Hetranso AHNI\thanks{Faculté des Sciences et Techniques, Institut Camille Jordan, Université Jean Monnet, France. Email: hetranso.ahni@univ-st-etienne.fr. The research of the author is supported by the doctoral studentship of the university Jean-Monnet.}
\qquad Jean-Fran\c{c}ois JABIR\thanks{Department of Statistics and Data Analysis, HSE University, Moscow, Russia. Email: jjabir@hse.ru. The research of the author is supported by the Russian Science Foundation project (project No. 24-11-00123).}
\qquad  Julian TUGAUT\thanks{Faculté des Sciences et Techniques, Institut Camille Jordan, Université Jean Monnet, France.
Email: julian.tugaut@univ-st-etienne.fr.}}
\maketitle
\begin{abstract}
    This paper investigates the exit-time asymptotics, at a vanishing noise limit, related to a stochastic particles system in mean-field interaction subjected to frictions, confining potential and flocking forces. The dynamics are governed by a kinetic Langevin-type SDE with an idiosyncratic noise acting solely on the velocity coordinates. Using the seminal large deviations principles of Freidlin and Wentzell \cite{Freidlin_Wentzell}, and their extension by Azencott \cite{Azencott} to the framework of non-elliptic diffusions, we establish the Arrhenius' law and Kramers' type law related to the first exit-time of the system from an energy level set $\mathcal O^\gamma_N(R)$ of interest.
    \end{abstract}
   \paragraph{Declaration of AI usage.} The study, the writing and overall analysis conducted in the present paper are the original works of the authors.  No generative AI or AI tools were used. 
\section{Introduction}
 In their seminal work \cite{Cucker-Smale_05}, Cucker and Smale developed a mathematical model describing an interacting system of $N$-agents leading to a ``flocking'' phenomenon -- that is: agents tend to align their change rates/velocities resulting in a agglomeration of their positions -- in the long-time regime. This model consists, for a fixed number of agents $N\geq 2$, of a coupled ODEs system given by:
 \begin{equation}
    \label{CS_o}
     \dot{x}_i(t)=v_i(t),\qquad
    \dot{v}_i(t)=\sum_{j=1}^N \phi(\norme{x_i(t)-x_j(t)})\left(v_j(t)-v_i(t)\right),
\end{equation}
with $\left(x_i(t),v_i(t)\right)\in \reel^d\times \reel^d$ describing the position-velocity of the agent $i$ at time $t$. The function $\phi$ given by 
\begin{equation}
    \label{decr_f}
    \phi(r)=\dfrac{\lambda}{(\kappa+r^2)^\beta},\quad \lambda,\kappa,\beta>0,
\end{equation}
defines the communication rate/kernel and models a set of interaction rules followed by each agent. When agents --which we shall sometimes refer to as particles as well-- are too far apart, interaction between them reduces, while, when closing, their velocity tend to drag one to the other. Defining the quantities 
\begin{equation*}
    \Gamma(t)=\dfrac{1}{2}\Sd\norme{x_i(t)-x_j(t)}^2\quad\mathrm{and}\quad\Lambda(t)=\dfrac{1}{2}\Sd\norme{v_i(t)-v_j(t)}^2\!,
\end{equation*}
that are respectively the agglomeration energy and the alignment energy of \eqref{CS_o}, Cucker and Smale \cite{Cucker-Smale_05} showed that, depending on the value of parameter $\beta$, \eqref{CS_o} exhibits (unconditional or $(x^i_0,v^i_0)$-conditional) flocking of the form:
\begin{align*}
    \lim_{t\to +\infty} \!\Lambda(t)=0 \quad \text{and}\quad \underset{t\geq 0}{\max}\;\Gamma(t)<+\infty.
\end{align*}
By means of suitable Lyapunov functions that highlight the intrinsic structure of \eqref{CS_o}, the authors in \cite{SHa_Liu_1} improved the previous result by lessening the assumptions in the unconditional cases. Since these breakthrough works, many  fruitful extensions have emerged in order to generalize the model to capture or fit more realistic features, such as in biology (swarm of bees, flock of birds, schooling of fish), robotic systems or more generally in engineering contexts \cite{White_Koch, Camazine, Camazine1, Aoki,Naldi_Par_Tos}. More specifically, we can mention:
\begin{itemize}
    \item The presence of heterogeneous ``dominating'' agents and leadership effects \cite{Ha-Li, Li-Ha-Xue,RenLi}, where some agents exert a stronger influence on the dynamics.
    \item Nonlinear velocity couplings cases where, instead of the linear drag force $v_j-v_i$ in \eqref{CS_o}, more general interaction rules beyond linear alignment, \cite{SHa_Kim, Naldi_Par_Tos,LiuRenZhang}, or matrix-valued communication kernels are taken into account, \cite{Shu_Tadmor, Ad_Cotil}. This enriches the interaction structure of the systems \eqref{CS_o} and allows for anisotropic effects.
    \item  Control problems where external agents (predator for instance) act on the interacting system \cite{Capo_Fornasier,Labi_Pareschi,ChoiOh,LiuRenZhang}.
    \item Mean-field limits, characterizing the statistical description of the multi-agents/particles system at the very large population limit ($N\rightarrow\infty$), \cite{SHa_Liu_1,nguyenShvyTan}.
    \item Flocking dynamics in random environments and other stochastic extensions,  \cite{Ref_Cattiaux,Ha_Je_Noh,Cari_For_Ro_Tos,Erban_Haskovec}: an additional noise force, modeling  possible random disturbances or small-scale fluctuations, is incorporated in the velocity dynamics -- the long-time robustness of the flocking takes then different forms, depending on whether the noise is in an additive, multiplicative, idiosyncratic, common, ... etc framework.
\end{itemize}
 In this work, we study the behavior of a particular noisy Cucker-Smale model with an external quadratic confining potential $U_i(x):=\dfrac{\gamma}{2}\norme{x-x_i^\ast}^2$, $\gamma>0$, and friction force. The model is described by the following system of second-order SDE systems:
 \begin{equation}
    \label{Cucker-Smale-bruite}
    \left\lbrace
    \begin{aligned}
        &\deriv X_t^{i,N,\varepsilon}=V_t^{i,N,\varepsilon}\deriv t,\\
        &\deriv V_t^{i,N,\varepsilon}=b^{i,N}\left[\left(X_t^{j,N,\varepsilon},V_t^{j,N,\varepsilon}\right)_{j}\right]\deriv t -\alpha\varia[V]{t}{i}{\varepsilon}\deriv t-\gamma\left(\varia[X]{t}{i}{\varepsilon}-x_i^\ast\right)\di t+\varepsilon \deriv W_t^i,\\
        &\left(X_0^{i,N,\varepsilon},V_0^{i,N,\varepsilon}\right)=\left(x_0^i,v_0^i\right),\qquad 1\le i\le N,
    \end{aligned}
    \right.
\end{equation}
where $\left(W^i\right)_{i\leq N}$ are independent $d$-dimensional Brownian motions, representing idiosyncratic noises, 
and the initial states $\left(x_0^i,v_0^i\right)$ and the concentrating points  $\{x_i^\ast\}_i\in\mathbb R^{dN}$ are all deterministic. For each $i$, the vector $x_i^*$ can be viewed as an individual target state to the agent $i$. The positive scalar $\varepsilon$ defines the shared noise intensity acting on each agent, the friction coefficient $\alpha$ and the confinement coefficient $\gamma$
are positive and fixed throughout. Up to a $1/N$-rescaling, the interaction kernels $b^{i,N}$ are defined (as originally considered in  \cite{Cucker-Smale_05}) by:
\begin{equation}
    \label{fonction_communis}
    b^{i,N}\left(\left(x_j,v_j\right)_j\right)=\dfrac{1}{N}\sum_{j=1}^N\dfrac{\lambda}{\left(\kappa+\lVert x_i-x_j\rVert^2\right)^\beta}\left(v_j-v_i\right),\qquad \kappa,\beta,\lambda>0.
\end{equation}

\subsection{The importance of friction and confinement: from mathematics to ecology}\leavevmode\par\vspace{0.5em}

In this paper, we move beyond the classical setting in \cite{Cucker-Smale_05,SHa_Liu_1,Ad_Cotil}. In fact, a key aspect of our work is the choice of adding a quadratic confining landscape and friction. Far from being a mere technical point, this choice comes from both physico-biological reasons and technical imperatives:
\begin{enumerate}
    \item[($\mathrm{I}$)] Physical realism: we add to the original model a friction term to embed a dissipation of the agent group's energy in order to fit more real-world context (like friction induced by the environment). Considering quadratic potential is very natural for a localized flock tied to a geographical anchor or a vital ecological niche. For instance, a school of fish may align to swim in order to reduce threats of predators but should remain near a coral reef for food. We choose the quadratic confinement as a baseline rather than a generic potential to reduce the technicality of the work and highlight the crux of the problem.
    \item [($\mathrm{II}$)] Metastability and exit-times: without confinement and friction, as shown in \cite{Cucker-Smale_05,SHa_Liu_1}, the center of mass of the particles system follows a uniform rectilinear motion making the associated exit-problem for domains of interest beyond the scope of Freidlin-Wentzell framework. 
   By adding confining and friction forces, we create a basin of attraction with one well in the phase space. Over long-time horizons, even an infinitesimal noise can induce rare events for the model regarding the behavior of its deterministic counterpart (when noise vanishes). For example, noise could drive the model out of its equilibrium, exiting from an invariant domain for the deterministic counterpart. This phenomenon is known as metastability.
  Then, the central question becomes resilience: how long interacting agents trapped in a well resist external small perturbation before a rare fluctuation takes away the entire group out of the well ?
\end{enumerate}
\subsection{Technical challenges: Degeneracy and locally Lipschitz continuous drifts}\leavevmode\par\vspace{0.5em}

From a technical point of view, the study of the dynamics \eqref{Cucker-Smale-bruite} in a small-noise regime can be done with the theory of large deviations principles (LDP). This approach was firstly developed in this context in Freidlin and Wentzell \cite{Freidlin_Wentzell} and was initially aimed to estimate rare events -- in the small-noise $\varepsilon\rightarrow 0$-regime -- related to the processes
\begin{equation}
    \label{overdamped}
     \di X_t^\varepsilon=b\left(X_t^\varepsilon\right)\di t+\varepsilon\di W_t,\qquad
    X_0^\varepsilon=x_0,
\end{equation}
For the sake of simplicity, we assume in this section $b$ is globally Lipschitz continuous. Rare events are estimated by means of a rate function defined by 
$$
\mathcal{I}_{0,T}(\varphi)=\frac{1}{2}\inte{T}\lVert \dot{\varphi}(t)- b(\varphi(t))\rVert^2\di t,\qquad \varphi\in \mathcal C^1([0,T])
$$
which should be understood as a similarity measure to the deterministic counterpart of \eqref{overdamped} given by the dynamic $\{X^0(t)=x_0+\int_0^tb(X^0(s))\,\di s\}_{t\ge 0}$ (corresponding to the case $\varepsilon=0$). When $X^0$ has an invariant domain $D$ (see~Definition~\ref{invar_domain}) containing a unique equilibrium point $x_{\text{eq}}$ that attracts the deterministic system starting for any point $D\cup\partial D$, one can guess, provided that $X^\varepsilon_0=x_0\in D$, an exit out of $D$ occurs at very large time because of the proximity between $X^0$ and $X^\varepsilon$ in small-noise regimes. Building upon the rate function, this guess can be measured by the \textit{quasi-potential}: 
\begin{equation*}
  V\left(x,y\right)=\underset{\substack{T>0,\\\varphi\in K_{x,y}^T}}{\inf}\,\mathcal{I}_{0,T}(\varphi),\qquad K_{x,y}^T:=\{\varphi\in H^1 :\,\varphi(0)=x,\;\varphi(T)=y\},  
\end{equation*}
 which characterizes the infrequency of the excursions of $X^\varepsilon$ out of $D$. As established in \cite{Freidlin_Wentzell}, the quasi-potential $V$ specifically characterizes the $\epsilon\rightarrow 0$ behavior of the exit--time 
$$
\tau^\varepsilon=\inf\{t>0: X^\varepsilon_t\notin D\},
$$ 
for starting points $x\in D$ with the following Arrhenius' law and Kramers' type law
\begin{align}
    \label{Arh_over}
    &\lim_{\varepsilon\to 0}\varepsilon^2\log\mathbb{E}_x(\tau^\varepsilon)=\underset{y\in \partial D}{\inf}V(x_{\text{eq}},y)=:V_0\quad\text{and}\quad\lim_{\varepsilon\to 0}\varepsilon^2\log(\tau^\varepsilon)=V_0\text{ in probability}. 
\end{align}
Comparatively, in our context, the dynamics \eqref{Cucker-Smale-bruite} take the form of a second order SDE where the communication part is only locally Lipschitz continuous. This induces an issue of  partially-degeneracy of the noise as Brownian perturbation acts solely on the velocity variables whereas positions evolve through the kinematic relation $\dot{X}=V$. This framework is nevertheless embedded in the work of Azencott \cite{Azencott} which extend Freidlin and Wentzell's LDP theory to non-homogeneous, eventually degenerate, diffusion matrices of class $C^1$. Note that since then, many works have been done to lessen assumptions ensuring LDPs for \eqref{overdamped} or integrate more complex situations \cite{Boue_Dupuis,Baldi_Maurel,Dembo_Zeitouni,Priouret,David_Lipschutz,NZi_Dakaou,Julian_1,Nua_Ro,Jac_Pa,Zhan_Xi}. However, a LDP does not systematically guarantee Arrhenius' law and Kramers' type law such as in \eqref{Arh_over}. Extending Freidlin-Wentzell's exit-laws \eqref{Arh_over} to the particles system \eqref{Cucker-Smale-bruite} is the crux of the paper. Notably, proving the controllability of the rate function (see Proposition \ref{ppdeva}) of our system under the structural constraint $\dot{X}=V$ is one of the fundamental bottleneck to overcome, in order to make the rate function sufficiently close to the exit-cost related to \eqref{Cucker-Smale-bruite}.

\subsection{Recent advances in the hypo-elliptic context}\leavevmode\par\vspace{0.5em}
The hypo-elliptic kinetic Langevin dynamic is generically described by the following SDE:
\begin{equation}
    \label{Langevin_cinetique}
    	\di X_t^{\varepsilon,\gamma}=V_t^{\varepsilon,\gamma}\di t,\quad
    	\di V_t^{\varepsilon,\gamma}=-\nabla F(X^{\varepsilon,\gamma}_t)\di t-\alpha V_t^{\varepsilon,\gamma} \di t+ \varepsilon\di W_t,
\end{equation}
with the potential $F:\reel^d\to\reel$ being $C^2$, non-convex, and admitting two (or a finite number of) minima. A major motivation behind the study of metastability in these systems or their long-time behavior is their broad range of applications, notably in simulated annealing (global optimization) \cite{Guillin_Mon,Monmarche}. In these settings, small random perturbations accelerate convergence towards the Gibbs measure $C_{\varepsilon,\alpha}\exp\{-\dfrac{2\alpha}{\varepsilon^2}\big(F(x)+\dfrac{|v|^2}{2}\big)\}$ and enable transitions between local minima, thereby increasing the probability of reaching the global minimum of a multi-wells potential and the speed of convergence too, compared to their elliptic versions.
 In this context, the rigorous analysis of metastable transitions between two equilibria $\left(m_1,0\right)$ and $\left(m_2,0\right)$,  with $m_1$, $m_2$  being the minima of $F$, has seen significant developments in recent years. While the potential theoretic approach  developed by the authors in \cite{bovier,bovier2,BovierHo} provides sharper asymptotic framework for the overdamped (elliptic) Langevin equation  \eqref{overdamped} (taking $b=\nabla F$ therein) --  namely a sharp prefactor in the Arrhenius'law \eqref{Arh_over} recovering the so-called Eyring-Kramers' law --, its extension to the underdamped (kinetic) regime has long been precluded by the structural degeneracy of the diffusion matrix in \eqref{Langevin_cinetique}. Recently, breakthrough works have targeted the global transition landscape of kinetic Langevin processes. Specifically, Lee, Ramil, and Seo \cite{Lee_Ramil_Seo} successfully established the full Eyring-Kramers' law for underdamped Langevin dynamics featuring a double-wells potential. Their approach relies on a tricky extension of the potential theory approach to hypo-elliptic diffusions enabling to compute the prefactor of Arrhenius' law. Furthermore, recent extensions have pushed these techniques toward even higher-order hypo-elliptic, such as third-order Langevin diffusions, \cite{wangzhu}. 

\subsection{Main contributions of the paper}\leavevmode\par\vspace{0.5em}
In our setting, we investigate the asymptotic behavior of the first exit-time of \eqref{Cucker-Smale-bruite} from an invariant domain of the phase space containing a stable equilibrium. Compared to \eqref{Langevin_cinetique}, the drifts $b^{i,N}$ in  \eqref{fonction_communis} depend both on the velocity and position; the methodology in \cite{Lee_Ramil_Seo} can not be directly applied. Moreover, our focus is on a probabilistic approach of the asymptotics of the exit-time which we provide in a three-step analysis of the metastable system \eqref{Cucker-Smale-bruite}. First, in Section~\ref{sec1}, we study the properties of the deterministic counterpart  -- when $\varepsilon$ is equal to $0$ -- by introducing suitable Lyapunov functions representing the agent group's energy, leading to the definition of an explicit invariant domain of interest and convergence to the equilibrium point $\left((x_i^\ast)_{1\le i\leq N},0_{\reel^{Nd}}\right)$ in the phase space. Second, in Section~\ref{sec2}, we establish estimates of the noisy particles system and its quantitative $L^1$-convergence toward its deterministic counterpart for the uniform topology on any finite horizon $[0;T]$ (the rate of convergence being $\mathcal{O}(\varepsilon)$). Finally, in Section \ref{sec3}, we establish Arrhenius' law and Kramers' type law to our framework considering the first exit-time from the stable domain; particularly, we show that, in small-noise regimes, the first exit-time scales exponentially as $\exp{\left({V}/\varepsilon^2\right)}$, in mean and in probability. The scalar $V$ represents the (positive) cost of exiting the domain in some finite time. In Section \ref{sec5}, we further illustrate these asymptotics with simulations. 

\subsection{Notations and conventions}\leavevmode\par\vspace{0.5em}

Throughout this paper, we adopt the following notations:
\begin{itemize}
    \item For $N\ge 1$, we denote by $\ens$ the set $\{1,2,3,\cdots,N\}$.
    \item We denote by $\lfloor r\rfloor:=\max\{k\in \mathbb{N}: k\leq r\}$ and $\lceil r\rceil:=\min\{k\in \mathbb{N}: k\geq r\}$, respectively the floor and ceiling functions on $r\geq 0$. For $x,y\in \reel$, we also define $x\land y:=\min(x,y)$ and $x\vee y:=\max(x,y)$.
    \item On $\reel^d$ or $\reel^{2d}$, we consider the Euclidean norm  $\norme{\cdot}$, the corresponding inner product $\langle\cdot,\cdot\rangle$ and denote related open and closed balls as $B(x,r):=\{y\in\reel^k:\norme{y-x}< r \}$ and $\overline{B(x,r)}:=\{y\in\reel^k:\norme{y-x}\leq r\} $.
    \item  On $\reel^{Nd}$ or $\reel^{2Nd}$, given $x_i,v_i\in \reel^d$ for $i\in \ens$, we put $\vv{x}:=(x_i)_{1\le i \leq N}\in \reel^{Nd}$, $\vv{v}:=(v_i)_{1\le i \leq N}\in \reel^{Nd}$ and $\vv{z}=[(x_i,v_i)_{1\le i\le n}]$. These spaces will be endowed with the corresponding $1$-norms: $\norme{\vv{x}}_1=\sum_{i=1}^N \norme{x_i}$ and $\norme{\vv{z}}_1=\norme{\vv{x}}_1+\norme{\vv{v}}_1$ or with the Euclidean norm $\norme{\cdot}$.
  \item We equip the set $\mathcal{B}\left([0;T],\reel^{Nd}\times\reel^{Nd}\right)$ of bounded function on $[0;T]$ taking values in $\reel^{Nd}\times \reel^{Nd}$ with the distance
 \begin{align*}
     &\mathcal{D}^T\big((f_1,f_2),(g_1,g_2)\big):=\underset{t\in [0;T]}{\max}\norme{(f_1(t),f_2(t))-(g_1(t),g_2(t))},
 \end{align*}
 for every function $(f_1,f_2),(g_1,g_2)\in \mathcal{B}\big([0;T],\reel^{Nd}\times \reel^{Nd}\big) $. We shall identify the set $\mathcal{B}\big([0;T],\reel^{Nd}\times\reel^{Nd}\big)$ to $\mathcal{B}\left([0;T],\reel^{2Nd}\right)$. Let $A$ and $B$ subsets of $\mathcal{B}\big([0;T],\reel^{2Nd}\big)$. We define
 \begin{equation*}
 	\mathcal{D}^T(\varphi,A):=\underset{\psi\in A}{\inf}\mathcal{D}^T(\varphi,\psi)
    \quad\mathrm{and}\quad\mathcal{D}^T(A,B):=\underset{\psi\in A,\varphi\in B}{\inf}\mathcal{D}^T(\varphi,\psi).
 \end{equation*}
 \item $C_{\vv{z}}^0\left([0;T],\mathbb{R}^{2Nd}\right)$ is the space of continuous functions from $[0;T]$ to $\mathbb{R}^{2Nd}$ starting in $\vv{z}=(\vv{x},\vv{v})$ at $t=0$, and $\espacecont$ is the set of continuous functions taking values in $\reel^{2Nd}$. We equip both spaces with $\mathcal{D}^T$ and with the time-shift operator $\theta_t$ and entering time functional $\tau_A$ (for $A\subset\reel^{2dN}$)
 \begin{equation}
 	\label{fonctionnelle}
 \theta_t((x_s,y_s)_{s\geq 0})=(x_{s+t},y_{s+t})_{s\geq 0}\quad\mathrm{and}\quad
 	\tau_A\left(\left(x_t,y_t\right)_{t\geq 0}\right)=\inf\{t\geq 0\,\,:\,\,(x_t,y_t)\in A\}.
 \end{equation}
 \item $C_{\text{abs}}^0\left([0;T],\mathbb{R}^{Nd}\right)$ denotes the space of absolutely continuous functions on $[0;T]$ taking values in $\reel^{Nd}$. The set  $H_2\left([0;T]\right)$ denotes the Sobolev space
 $$
 H_2\left([0;T]\right)~:=\{\varphi \in C^1\left([0;T],\mathbb{R}^{Nd}\right):\,\dot{\varphi}\in C_{\mathrm{abs}}^0\left([0;T],\mathbb{R}^{Nd}\right),\inte{T}\norme{\ddot{\varphi}(t)}^2\di t<+\infty\}.
 $$
   \item Finally, when appropriate, we shall write  $\probabi{}$, or $\mathbb P_{\initia}$, for the probability measure given $(X^{N,\varepsilon}_0,V^{N,\varepsilon}_0)=\initia=\vv{z_0}$ and write $\mathbb E_{(\vv{x_0},\vv{v_0})}$ or $ \esper{}$ for the corresponding expectation.
\end{itemize}
\section{Analysis of the deterministic system and related invariant domains}\label{sec1}
In this section, we study the deterministic counterpart of \eqref{CS_o}, given, for $N\geq 2$, by the following ODE systems 
\begin{equation}
    \label{C-S-conf_1}
    \left\lbrace
        \begin{aligned}
            \dot{x}_i(t)&=v_i(t),\\
            \dot{v}_i(t)&=\dfrac{1}{N}\sum_{j=1}^N\phi(\lVert x_i(t)-x_j(t)\rVert)(v_j(t)-v_i(t))-\alpha v_i(t)-\gamma(x_i(t)-x^*_i),\quad (x_i(0),v_i(0))=(x_0^i,v_0^i),
        \end{aligned}
    \right.
\end{equation}
the communication kernel $\phi$ being as in \eqref{decr_f}. Since $\phi$ and other friction and confining forces are locally Lipschitz continuous with at most linear growth, this system admits a unique non-exploding solution at any time $t\in[0,\infty)$. Note that the ``communicating'' part $\displaystyle \phi_{i,j}(t):=\phi(\lVert x_i(t)-x_j(t)\rVert)$ is symmetric with $\phi_{i,j}=\phi_{j,i}$ and so we have $\displaystyle\Sd\phi(\lVert x_i(t)-x_j(t)\rVert)(v_j(t)-v_i(t))=0$. Hence, the center of mass $\displaystyle\overline{x}(t)=\frac1{N}\sum_{i=1}^Nx_i(t)$ satisfies the second order linear ODE: $\displaystyle
\ddot{\overline{x}}+\alpha\dot{\overline{x}}+\gamma\overline{x}=\gamma\overline{x}^\ast$ for $\displaystyle\overline{x}^*:=\frac1{N}\sum_{i=1}^Nx^*_i$. Then, it is easy to show the macroscopic limit:
\begin{equation*}
    \Big(\overline{x}(t),\overline{v}(t):=\frac1{N}\sum_{i=1}^Nv_i(t)\Big)\underset{t\to+\infty}{\longrightarrow}\left(\overline{x}^\ast,0_{\mathbb{R}^{d}}\right).
\end{equation*}
These macroscopic components provide a reference understanding of the system's behavior, which we show to hold true also for individual agents by introducing suitable Lyapunov functions. From here on we denote by $\pointeq$, $\vv{x}_{\mathrm{eq}}^\ast=\left(x_i^*\right)_{i\in\ens}$, the equilibrium point of \eqref{C-S-conf_1}.

\subsection{Energy dissipation and construction of stable domains}\leavevmode\par\vspace{0.5em}
Consider the following Lyapunov 
functionals: for $\vv{z}=(\vv{x},\vv{v})=([(x_i,v_i)]_{1\le i\le N})\in\mathbb R^{2Nd}$, 
\begin{equation}
    \label{ener_u} 
    \mathcal{U}^{\gamma}_{\mathrm{tot}}(\vv{z}):=\sum_{i=1}^N\lVert v_i\rVert^2+\gamma\sum_{i=1}^N\left\lVert x_i -x_i^\ast\right\rVert^2\quad\mathrm{and}\quad\mathcal{U}^{\gamma,N}_{\mathrm{tot}}(t):=\mathcal{U}^{\gamma}_{\mathrm{tot}}(\left[x_i(t),v_i(t)\right]_{1\le i\leq N}), 
\end{equation}
which, from a physical standpoint, corresponds to the total energy, namely the alignment plus the agglomeration energy, of the system~\eqref{C-S-conf_1}. Here, the parameter $\gamma$ defines a homogenization coefficient. The level set of the total energy allows us to construct a stable domain in the following sense.
\begin{definition}
\label{invar_domain}
 We say that an open set $D_{\mathrm{inv}}\subset\reel^{2Nd}$ is a stable domain (or invariant domain) for the system~\eqref{C-S-conf_1} if any trajectory starting in $D_{\mathrm{inv}}$ remains in $D_{\mathrm{inv}}$ at any time, namely:
 \begin{equation*}
     \initia \in D_{\mathrm{inv}} \implies \left[x_i(t),v_i(t)\right]_{1\le i\leq N} \in D_{\mathrm{inv}},\quad  \text{for all}\quad t\geq 0.
 \end{equation*}
\end{definition}
The following lemma states a first straightforward estimate of the dissipation of $\mathcal{U}^{\gamma,N}_{\mathrm{tot}}$.
\begin{lemme}
    \label{Lyapou_cible_naive}
    The total energy $\mathcal{U}^{\gamma,N}_{\mathrm{tot}}$ of the particles system \eqref{C-S-conf_1} is non-increasing with, for all $t\geq 0$, 
    \begin{equation}
        \label{derivee_th3}
        \dfrac{\di}{\di t}\mathcal{U}^{\gamma,N}_{\mathrm{tot}}(t)=-\dfrac{1}{N}\sum_{i,j=1}^N\phi(\lVert x_i(t)-x_j(t)\rVert)\lVert v_j(t)-v_i(t)\rVert^2-2\alpha\sum_{i=1}^N\lVert v_i(t)\rVert^2.
    \end{equation}
\end{lemme}
\begin{proof}
    The function $t\longmapsto\mathcal{U}^{\gamma,N}_{\mathrm{tot}}(t)$ is differentiable and, for all $t\geq 0$, we have 
    \begin{align*}
        \dfrac{\di}{\di t}\mathcal{U}^{\gamma,N}_{\mathrm{tot}}(t)&=2\sum_{i=1}^N\langle \dot{v}_i(t),v_i(t)\rangle+2\gamma\sum_{i=1}^N \left\langle v_i(t),x_i(t)-x_i^\ast\right\rangle\\
        &=\dfrac{2}{N}\sum_{i,j=1}^N\phi(\lVert x_i(t)-x_j(t)\rVert) \langle v_j(t)-v_i(t),v_i(t)\rangle-2\alpha\sum_{i=1}^N\langle v_i(t),v_i(t)\rangle\\
        &\quad-2\gamma\sum_{i=1}^N\langle v_i(t),x_i(t)-x_i^\ast\rangle+2\gamma\sum_{i=1}^N \left\langle v_i(t),x_i(t)-x_i^\ast\right\rangle\\
        &=\dfrac{2}{N}\sum_{i,j=1}^N\phi(\lVert x_i(t)-x_j(t)\rVert) \langle v_j(t)-v_i(t),v_i(t)\rangle-2\alpha\sum_{i=1}^N\lVert v_i(t)\rVert^2.
    \end{align*}
    Using the identity $\langle a_j-a_i,a_i\rangle=\dfrac{1}{2}\left(\lVert a_j\rVert^2-\lVert a_i\rVert^2-\lVert a_j-a_i\rVert^2\right)$, and, as the communication kernel $r\mapsto \phi(r)$ in~\eqref{decr_f} is even,
    \begin{equation*}
        \dfrac{\di}{\di t}\mathcal{U}^{\gamma,N}_{\mathrm{tot}}(t)=-\dfrac{1}{N}\sum_{i,j=1}^N\phi(\lVert x_i(t)-x_j(t)\rVert)\lVert v_j(t)-v_i(t)\rVert^2-2\alpha\sum_{i=1}^N\lVert v_i(t)\rVert^2.
    \end{equation*}
\end{proof}
We immediately deduce the following corollary.
\begin{corollaire}
    \label{lemm_do_inv}
    For any $R>0$ the domain
    \begin{equation}
        \label{Domaine_inv_classique}
        \ensinva{R}:=\{\vv{z}:= (x_i,v_i)_{i\leq N}\in \reel^{Nd}\times \reel^{Nd}:\,\mathcal{U}^{\gamma}_{\mathrm{tot}}(\vv{z})=\Ss\left\lVert v_i\right\rVert^2+\gamma\Ss\left\lVert x_i-x_i^\ast\right\rVert^2<R \}.
    \end{equation}
     is stable for the system \eqref{C-S-conf_1}.
\end{corollaire}
\subsection{A refined Lyapunov function and asymptotic stability}\label{sec:RefinedLyapunov}\leavevmode\par\vspace{0.5em}
Lemma  \ref{Lyapou_cible_naive} provides a first result regarding the dissipation involved in~\eqref{C-S-conf_1}. However, this dissipation is \textit{a priori} not sufficient to establish the asymptotic stability of the equilibrium $\pointeq$ since the agglomeration part of $\mathcal{U}^\gamma_{\mathrm{tot}}$ is missing in its derivative. To overcome this,  we modify the energy  $\mathcal{U}^\gamma_{\mathrm{tot}}$ by introducing a cross-energy term that captures the interplay between velocities and positions. Let us consider
 \begin{equation}\label{energie-modifiee}
 	\mathcal{U}^{\sigma,\gamma}(\vv{z})=\mathcal{U}^\gamma_{\mathrm{tot}}(\vv{z})+\sigma\sum_{i=1}^N\left\langle v_i,x_i-x_i^\ast\right\rangle\quad\mathrm{and}\quad\mathcal{U}^{\sigma,\gamma,N}(t)=\mathcal{U}^{\sigma,\gamma}([(x_i(t),v_i(t))_{1\le i\le N}]).
    \end{equation}
 Provided that $\sigma$ is under an explicit threshold, $\mathcal{U}^{\sigma,\gamma,N}(t)$ and   $\mathcal{U}^{\gamma,N}_{\mathrm{tot}}(t)$ are equivalent in the following sense.
\begin{corollaire}
    \label{dissipatif_part}
    For all $0\leq \sigma\leq\sqrt{\gamma}$ and all $t\geq 0$, we have:
    \begin{equation}
       \dfrac{1}{2}\mathcal{U}^{\gamma,N}_{\mathrm{tot}}(t)\leq \mathcal{U}^{\sigma,\gamma,N}(t)\leq \dfrac{3}{2}\mathcal{U}^{\gamma,N}_{\mathrm{tot}}(t).
    \end{equation}
\end{corollaire}
\begin{proof}   The claim is immediate from Cauchy-Schwarz and Young inequalities, for which we have
    \begin{align*}
        \left\lvert \sigma\sum_{i=1}^N\left\langle v_i(t),x_i(t)-x_i^\ast\right\rangle\right\rvert&\leq\sum_{i=1}^N\lVert v_i(t)\rVert \sigma\left\lVert x_i(t) -x_i^\ast\right\rVert\leq\dfrac{1}{2}\left[\sum_{i=1}^N\lVert v_i(t)\rVert^2+\dfrac{\sigma^2}{\gamma}\gamma\sum_{i=1}^N\left\lVert x_i(t) -x_i^\ast\right\rVert^2\right].
    \end{align*}
\end{proof}
\noindent
We state the main result of this section, asserting the exponential decay of the total energy $ \mathcal{U}^{\gamma,N}_{\mathrm{tot}}(t)$  
via $\mathcal{U}^{\sigma,\gamma,N}(t)$.
 \begin{theoreme}
    \label{theoreme_Ineq_Diff_1}
    Given $\theta<\gamma$, and  
    \begin{align*}
    	&\sigma_0:=\min\Bigg(\sqrt{\gamma},\dfrac{2\alpha}{1+\frac{(\phi^+)^2}{\gamma}+\frac{\alpha^2}{2\theta}}\Bigg)\quad\mathrm{and}\quad\delta:=\min\bigg(2\alpha-\bigg(\dfrac{(\phi^{+})^2\ }{\gamma}+\dfrac{\alpha^2}{2\theta}+1\bigg)\sigma,\dfrac{\sigma(\gamma-\theta)}{2\gamma}\bigg),
    \end{align*}
    it holds that: for every $\sigma\in (0,\sigma_0)$, $t\ge 0$, \begin{equation}
        \label{energie_mod_U}
            \mathcal{U}^{\sigma,\gamma,N}(t)\leq \mathcal{U}^{\sigma,\gamma,N}(0)e^{-\frac{2}{3}\delta t} \quad\mathrm{and}\quad \mathcal{U}^{\gamma,N}_{\mathrm{tot}}(t)\leq 3\mathcal{U}^{\gamma,N}_{\mathrm{tot}}(0)e^{-\frac{2}{3}\delta t}.
    \end{equation}
\end{theoreme}
\begin{proof} Owing to Corollary~\ref{dissipatif_part}, it is enough to estimate $\mathcal{U}^{\sigma,\gamma,N}(t)$. Using \eqref{derivee_th3}, we have, for all $t\geq 0$ :
    \begin{align*}
        \dfrac{\di}{\di t}\mathcal{U}^{\sigma,\gamma,N}(t)&=\dfrac{\di}{\di t}\mathcal{U}^{\gamma,N}_{\mathrm{tot}}(t)+\sigma\sum_{i=1}^N\left\langle \dot{v}_i(t),x_i(t)-x_i^\ast\right\rangle+\sigma\sum_{i=1}^N\lVert v_i(t)\rVert^2\\
        &=-\dfrac{1}{N}\sum_{i,j=1}^N\phi(\lVert x_i(t)-x_j(t)\rVert)\lVert v_j(t)-v_i(t)\rVert^2+(\sigma-2\alpha)\sum_{i=1}^N\lVert v_i(t)\rVert^2\\
        &\quad+\dfrac{\sigma}{N}\sum_{i,j=1}^N\phi(\lVert x_i(t)-x_j(t)\rVert)\left\langle x_i(t)-x_i^\ast,v_j(t)-v_i(t)\right\rangle(=:\mathcal U_1(t))\\
        &\quad-\sigma\alpha\sum_{i=1}^N\left\langle x_i(t)-x_i^\ast,v_i(t)\right\rangle(=:\mathcal U_2(t))-\sigma\gamma\sum_{i=1}^N\left\langle x_i(t)-x_i^\ast,x_i(t)-x_i^\ast\right\rangle\\
        &\le (\sigma-2\alpha)\sum_{i=1}^N\lVert v_i(t)\rVert^2+\mathcal{U}_1(t)+\mathcal{U}_2(t)-\sigma\gamma\sum_{i=1}^N\|x_i(t)-x^*_i\|^2.
\end{align*}
\textbf{For $\mathcal{U}_1$.} As $\phi$ given by \eqref{decr_f} is bounded from above by $\phi^+:=\tfrac{\lambda}{\kappa^\beta}$, Cauchy-Schwarz and Young inequalities yield  
\begin{align}
    \label{U4}
    \mathcal{U}_1(t) 
    &\leq \dfrac{1}{N}\sum_{i,j=1}^N\sqrt{\sigma\gamma}\left\lVert x_i(t) -x_i^\ast\right\rVert\dfrac{\phi^+\sqrt{\sigma}}{\sqrt{\gamma}}\lVert v_j(t)-v_i(t)\rVert\notag\\
    &\leq \dfrac{1}{N}\sum_{i,j=1}^N\dfrac{\sigma\gamma}{2}\left\lVert x_i(t) -x_i^\ast\right\rVert^2+\dfrac{1}{N}\sum_{i,j=1}^N\dfrac{(\phi^+)^2\sigma}{2\gamma}\lVert v_j(t)-v_i(t)\rVert^2\notag\\
    &=\dfrac{\sigma\gamma}{2}\sum_{i=1}^N\left\lVert x_i(t) -x_i^\ast\right\rVert^2+\dfrac{(\phi^+)^2\sigma}{\gamma}\sum_{i=1}^N\lVert v_i(t)\rVert^2-\dfrac{(\phi^+)^2\sigma}{N\gamma}N^2\lVert \overline{v}(t)\rVert^2\notag\\
    &\leq\dfrac{\sigma\gamma}{2}\sum_{i=1}^N\left\lVert x_i(t) -x_i^\ast\right\rVert^2+\dfrac{(\phi^+)^2\sigma}{\gamma}\sum_{i=1}^N\lVert v_i(t)\rVert^2.
\end{align}
\textbf{For $\mathcal{U}_2$.} Let $\theta>0$. Then, applying (again) Cauchy-Schwarz and Young inequalities, we get
\begin{align}
    \label{U5}
    \mathcal{U}_2(t) 
    &\leq\dfrac{\sigma\theta}{2}\sum_{i=1}^N\left\lVert x_i(t) -x_i^\ast\right\rVert^2+\dfrac{\alpha^2\sigma}{2\theta}\sum_{i=1}^N\lVert v_i(t)\rVert^2.
\end{align}
Combining the upper-bounds \eqref{U4} and \eqref{U5} yields 
\begin{align*}
    \dfrac{\di}{\di t}\mathcal{U}^{\sigma,\gamma,N}(t)&\leq-\left(2\alpha-\sigma-\dfrac{(\phi^+)^2\sigma}{\gamma}-\dfrac{\alpha^2\sigma}{2\theta}\right)\sum_{i=1}^N\lVert v_i(t)\rVert^2-\left(\dfrac{\sigma(\gamma-\theta)}{2}\right)\sum_{i=1}^N\left\lVert x_i(t) -x_i^\ast\right\rVert^2.
\end{align*}
Choosing $\theta<\gamma$, $\sigma\in(0,\sigma_0)$, and owing to the values of $\sigma_0$ and $\delta$ and to Corollary~\ref{dissipatif_part}, we conclude the claim with
\begin{equation}
	\label{ineq_vitesse_moy}
	\dfrac{\di}{\di t}\mathcal{U}^{\sigma,\gamma,N}(t)\leq -\delta \mathcal{U}^{\sigma,\gamma,N}_{\mathrm{tot}}(t)\leq -\dfrac{2\delta}{3} \mathcal{U}^{\sigma,\gamma,N}(t).
\end{equation}
\end{proof}
\section{Small-noise consistency and quantitative \texorpdfstring{$L^1$}{L1} stability }\label{sec2}
In this section, we return to the stochastic Cucker-Smale model \eqref{Cucker-Smale-bruite} discussed \textit{ab initio} from the outset in the introduction. For notation convenience, we put 
\begin{equation}
Z^{i,N,\varepsilon}_t=\left(X_t^{i,N,\varepsilon},V_t^{i,N,\varepsilon}\right),\qquad	\vv{Z}^{N,\varepsilon}_t=(\vv{X}_t^{N,\varepsilon},\vv{V}_t^{N,\varepsilon}).
\end{equation}
We also denote $\left(X_t^{i,N,0},V_t^{i,N,0}\right)= \left(x^{i,N}(t),v^{i,N}(t)\right)$ the solution of the zero-noise system \eqref{C-S-conf_1} and by $\left(\vv{x_0},\vv{v_0}\right):=\left[(x_0^i)_{i\leq N},(v_0^i)_{i\leq N}\right]$ the common initial conditions of the systems. As (again) the driving force fields are locally Lipschitz and have at most linear growth, the SDE \eqref{Cucker-Smale-bruite} admits a unique strong solution such that, for any finite time horizon $T>0$, 
\begin{equation}\label{MomentCarre}
\mathbb E\Big(\max_{t\in[0,T]}\|\vv{Z}^{N,\varepsilon}_t\|^2\Big)=\mathbb E\Big(\max_{t\in[0,T]}\big(\|\vv{X}^{N,\varepsilon}_t\|^2+\|\vv{V}^{N,\varepsilon}_t\|^2\big)\Big)\le C<\infty.
\end{equation}
(see e.g. \cite[Theorem 3.1.1]{PreRock-07} and \cite[Lemma 3.2]{Mao-11}). Below, in the small-noise regime, we show that for a fixed time horizon $T>0$, the deterministic $(\vv{x}^N,\vv{v}^N)$ is a $L^1$ approximation of $\grandeva{}$ for the uniform topology. This provides a heuristic guarantee that  $\grandeva{}$ will spend a significant amount of time near the equilibrium $\pointeq$ when noise is small. While this convergence is very standard when the drift of the SDE \eqref{overdamped} is globally Lipschitz continuous (see \cite[Th~1.2,p.~30]{Freidlin_Wentzell}), in our context it is only locally Lipschitz continuous. This will eventually not be a strong limitation to overcome since we exploit the boundedness of $\mathcal{U}^{\gamma,N}_{\mathrm{tot}}$ independently of $t$ (see Lemma~\ref{Lyapou_cible_naive}). Let us define the following functions:
\begin{align*}
    b^N(\vv{z}=\initias):=\left(\left[\Psi^N_i(\vv{z})-\alpha v_i-\gamma(x_i-x_i^\ast)\right]_{i\leq N}\right),\;
    \mathrm{with}\;\;\Psi^N_i(\vv{z})= \frac1{N}\sum_{j=1}^N \phi(\norme{x_i-x_j})(v_j-v_i).
\end{align*}

We now state the main theorem of this section  showing, for the uniform topology on the path space, the  $L^1$ convergence of the system \eqref{Cucker-Smale-bruite} to its deterministic counterpart as $\varepsilon$ vanishes.

\begin{theoreme}
    \label{convL1}
    For all $T>0$, we have, for some constant $\widehat C=\widehat C\big(\alpha,\phi_+,\beta,\gamma,\norme{\vv{x}_0},\norme{\vv{v}_0}\big)$ (for $\phi^+=\lambda/\kappa^\beta$),
    \begin{equation}\label{constante_phi}
    	\esp{\max_{0\le t\le T}\norme{\vv{V}_t^{N,\varepsilon}-\vv{v}^{N}(t)}^2}\leq 2 NCdT\varepsilon^2\exp{\left(\widehat CT\right)},
    \end{equation}
for some positive $C$ and
\begin{equation}\label{constante_phi_bis}
	 \lim_{\varepsilon\to 0}\left[\esp{\max_{0\le t\le T}\norme{\vv{X}_t^{N,\varepsilon}-\vv{x}^{N}(t)}}\vee\esp{\max_{0\le t\le T}\norme{\vv{V}_t^{N,\varepsilon}-\vv{v}^{N}(t)} }\right]=0. 
\end{equation}
 In particular, for every compact set $\mathcal{K}$ of $\reel^{2Nd}$ and every $\rho,T>0$, there exists $\varepsilon_0>0$ such that, for every  $\vv{z}_0=(\vv{x_0},\vv{v_0})\in \mathcal{K}$ and $\varepsilon\leq \varepsilon_0$, we have 
\begin{equation}
\label{convunifo}
	\max\left(\mathbb E_{\vv{z_0}}\left[\supre{T}\norme{\vv{V}_t^{N,\varepsilon}-\vv{v}^{N}(t)}\right],\mathbb E_{\vv{z_0}}\left[\supre{T}\norme{\vv{X}_t^{N,\varepsilon}-\vv{x}^{N}(t)}\right]\right)\leq \rho.
\end{equation}
\end{theoreme}
	\begin{proof} According to Lemma \ref{Lyapou_cible_naive} we have $\underset{t\geq 0}{\sup}\,\|\vv{v}^N(t)\|\le \sqrt{\gamma\|\vv{x_0}-\vv{x}^\ast_{\rm eq}\|^2+\|\vv{v}_0\|^2}=C\initia=:C_0$. Using the convexity inequality $(a+b)^2\leq 2(a^2+b^2)$, $a,b\in\reel^+$, we have, for all $0\le t\le T$,
		\begin{align*}
			&\esp{\underset{0\leq s\leq t }{\max}\norme{\vv{V}^{N,\varepsilon}_s-\vv{v}^{N}(s)}^2}\\
			&\le 2\esp{\underset{0\leq s\leq t }{\max}\norme{\inte{s}\left[b^{N}\left(\vv{X}_r^{N,\varepsilon},\vv{V}_r^{N,\varepsilon}\right)-b^{N}\left(\vv{x}^{N}(r),\vv{v}^{N}(r)\right)\right]\deriv r}^2}+ 2 \varepsilon^2\esp{\underset{0\leq s\leq t }{\max}\norme{W_s}^2},
		\end{align*}
		where $W:=(W^i)_i$ is a $Nd$-dimensional Brownian motion. Moreover, by standard Brownian estimates we obtain  $\esp{\underset{0\leq s\leq t }{\max}\norme{W_s}^2}\leq NCdT$ for some positive $C$.\\ Although the ``interaction part'' $\Psi_i^N(\vv{x},\vv{v})=\frac1{N}\sum_{j=1}^N \phi(\norme{x_i-x_j})(v_j-v_i)$ in $b^N$ is not globally Lipschitz continuous, one can easily check that, for all $\vv{z}=(\vv{x},\vv{v})$,  $\vv{z}'=(\vv{x}',\vv{v}')$ such that $\|\vv{v}'\|\le C_0 $ and for all $i\in \ens$  
		\begin{align*}
		\|\Psi_i^N(\vv{z})-\Psi_i^N(\vv{z}')\|&=\frac1{N}\left\|\sum_{j=1}^N\Big( \phi(\norme{x_i-x_j})(v_j-v_i)-\phi(\norme{x'_i-x'_j})(v'_j-v'_i)\pm \phi(\norme{x_i-x_j})(v'_j-v'_i)\Big)\right\|
		\\&\le \frac1{N}\left\|\sum_{j=1}^N \phi(\norme{x_i-x_j})\Big((v_j-v_i)-(v'_j-v'_i)\Big)\right\|\\
        &\quad+\frac1{N}\left\|\sum_{j=1}^N (v'_j-v'_i)\Big(\phi(\norme{x_i-x_j})-\phi(\norme{x'_i-x'_j})\Big)\right\|\\
		&\le 2\phi^+\|\vv{v}-\vv{v}'\|+4 C_0\phi_{\mathrm{Lip}}\|\vv{x}-\vv{x}'\|,
		\end{align*}
		with $\phi$ being $\phi_{\mathrm{Lip}}:=2\lambda\beta(2\beta+1)^{\beta-\frac{1}  {2}}\kappa^{\frac{1}{2}-\beta}(2\beta+2)^{-\beta}$ Lipschitz continuous. The friction and confining force composing the rest of $b^N$ being linear, using Cauchy-Schwarz's inequality and the convexity inequality we have
		\begin{align*}
		&\underset{0\leq s\leq t }{\max}\norme{\inte{s}\left[b^{N}\left(\vv{X}_r^{N,\varepsilon},\vv{V}_r^{N,\varepsilon}\right)-b^{N}\left(\vv{x}^{N}(r),\vv{v}^{N}(r)\right)\right]\deriv r}^2\\
		&\le 2TN \int_0^t\left[\big(\gamma+4 C_0\phi_{\mathrm{Lip}}\big)^2\underset{0\leq r\leq s}{\max}\norme{\vv{X}_r^{N,\varepsilon}-\vv{x}^{N}(r)}^2+\big(2\phi^++\alpha\big)^2\underset{0\leq r\leq s}{\max}\norme{\vv{V}_r^{N,\varepsilon}-\vv{v}^{N}(r)}^2\right]\,\di s.
		\end{align*}
		Since
		\begin{equation}\label{331}
			\supre{T}\norme{\vv{X}_t^{N,\varepsilon}-\vv{x}^{N}(t)}\leq T\supre{T}\norme{\vv{V}_t^{i,N,\varepsilon}-\vv{v}^{N}(t)}.
		\end{equation}
		we eventually deduce that, for $\widehat{C}=4\max\left(\big(\gamma+4 C_0\phi_{\mathrm{Lip}}\big)^2,(2\phi^++\alpha)^2\right)NT(1+T^2)$
		\begin{equation*}
			\esp{\underset{0\leq s\leq t }{\max}\norme{\vv{V}^{N,\varepsilon}_s-\vv{v}^{N}(s)}^2}
			\le \widehat C\int_0^t\esp{\underset{0\leq r\leq s }{\max}\norme{\vv{V}^{N,\varepsilon}_r-\vv{v}^{N}(r)}^2}
			\deriv r+ 2 \varepsilon^2NCdT.
		\end{equation*}
Applying Gr\"onwall lemma yields \eqref{constante_phi}. The second claim \eqref{constante_phi_bis} follows from \eqref{331} by taking in the constant $\widehat C$ the bound $C_0=\underset{\vv{z_0}\in\mathcal K}{\max}\big(\gamma\|\vv{x_0}-\vv{x}^\ast_{\rm eq}\|^2+\|\vv{v}_0\|^2\big)^{1/2}$.
\end{proof}

\section{Arrhenius' law and Kramers' type law for \texorpdfstring{\eqref{Cucker-Smale-bruite}}{eq}}\label{sec3}
In this section, we establish the Arrhenius' law (see Theorem \ref{loi_arrhenius}) and a Kramers'type law (see Theorem \ref{l_kram}) for the exit-time:
\begin{equation}
	\tau^{\alpha,\gamma,N}_{\varepsilon}(R)=\inf\{t>0:\,\left(\varia[X]{t}{i}{\varepsilon},\varia[V]{t}{i}{\varepsilon}\right)_i\notin \mathcal{O}_N^{\gamma}(R) \},
\end{equation}
as discussed in the introduction. We have shown with Corollary \ref{lemm_do_inv} that starting in $\ensinva{R}$, the deterministic counterpart of \eqref{Cucker-Smale-bruite} remains stuck in $\ensinva{R}$ and converges to equilibrium $\vv{z}^*_{\mathrm{eq}}:=\pointeq$ (see Theorem \ref{theoreme_Ineq_Diff_1}). Note that the domain $\ensinva{R}$ represents a neighborhood of the equilibrium in the phase space for the system \eqref{C-S-conf_1}, and has been constructed due to the kinetic nature of our problem. Exiting $\ensinva{R}$ does not imply, for instance, a large deviation of the positions from their equilibria even with a reasonable increase of the threshold $R$. Indeed, in presence of small-noise, its boundary $\partial \mathcal{O}_N^{\gamma}(R)$ could be reached with a substantial increase of the kinetic part of the energy associated to the particles system (see $\mathcal{E}_t$ below) while the positions remains in a small vicinity of their equilibria at the exit. Thus, the stopping-time $\ta$ represents the persistence time of the coherent (both in velocity and position) equilibrium in the phase space under stochastic perturbations.

Proposition \ref{propestimgro} first shows us the scarcity of  the excursions out of $\mathcal{O}_N^{\gamma}(R)$ by establishing rough estimates of $\espa{\tau^{\alpha,\gamma,N}_{\varepsilon}(R)}$. To achieve it, consider the energy functional: 
 \begin{equation*}
       \mathcal{E}_t^{\alpha,N,\varepsilon} =\Ss\norme{\varia[V]{t}{i}{\varepsilon}}^2+\gamma\Ss\norme{\varia[X]{t}{i}{\varepsilon}-x_i^\ast}^2,
\end{equation*}
 which is the stochastic version of $\mathcal{U}^\gamma_{\mathrm{tot}} $ defined in \eqref{ener_u}. (For the notation convenience, we drop the explicit writing of the parameter $\gamma,\alpha,N$ and $\varepsilon$.) 
 \begin{proposition}
    \label{propestimgro}
    Consider the initial condition of~\eqref{Cucker-Smale-bruite}  with $\vv{z_0}=\initia\in \mathcal{O}_N^{\gamma}(R)$. 
    Then we have 
    $$
    \espa{\tau^{\alpha,\gamma,N}_{\varepsilon}(R)}\geq \dfrac{R-\mathcal{E}_0}{Nd\varepsilon^2},\quad \forall\varepsilon\text{ sufficiently small}. 
    $$
    Moreover, we have the estimate
    \begin{align}
    \label{estim_temps_sortie_1}
    &\espa{\tau^{\alpha,\gamma,N}_{\varepsilon}(R)}\leq \dfrac{R-\mathcal{E}_0}{Nd\varepsilon^2-2(\phi^++\alpha)R},\quad\forall \varepsilon>\sqrt{\dfrac{2(\phi^++\alpha)R}{Nd}},\quad\text{with}\quad \phi^+=\dfrac{\lambda}{\kappa^\beta}.
    \end{align}
\end{proposition}
\begin{proof}
  Invoking the Itô formula, we get
  \begin{align}
    \label{itoener}
    \di \mathcal{E}_t
    &=-\dfrac{\lambda}{N}\Sd \dfrac{1}{(\kappa+\lVert \varia[X]{t}{i}{\varepsilon}-\varia[X]{t}{j}{\varepsilon}\rVert^2)^\beta}\norme{\varia[V]{t}{i}{\varepsilon}-\varia[V]{t}{j}{\varepsilon}}^2\di t-2\alpha \Ss\norme{\varia[V]{t}{i}{\varepsilon}}^2\di t\notag\\
    &\quad+N d\varepsilon^2\di t+\varepsilon \Ss \langle \varia[V]{t}{i}{\varepsilon}; \di W_t^i \rangle .
\end{align}
Since the stopped process $\displaystyle\left(\inte{\tau^{\alpha,\gamma,N}_{\varepsilon}(R)\land t}\Ss \langle \varia[V]{s}{i}{\varepsilon}; \di W_s^i \rangle\right)_{t\geq0}$ is a martingale, taking the expectation of $\mathcal{E}_t$ yields 
\begin{equation}\label{411}
    \esp{\mathcal{E}_{\tau^{\alpha,\gamma,N}_{\varepsilon}(R)\land t}}\leq \mathcal{E}_0+\varepsilon^2dN\espa{\tau^{\alpha,\gamma,N}_{\varepsilon}(R)\land t},
\end{equation}
 for all $t\ge 0$.
Let us now suppose that there is $\varepsilon_0>0$ such that 
    \begin{equation}
    \label{cond_petit}
        \espa{\tau^{\alpha,\gamma,N}_{\varepsilon}(R)}<+\infty,\forall\varepsilon<\varepsilon_0,
    \end{equation}
    otherwise the lower-bound is trivial. (We will address this in Remark \ref{RemValid} below.) Applying Lebesgue dominated convergence theorem to the left-hand side and  monotone convergence theorem for the right-hand side of \eqref{411}, for every $\varepsilon\leq\varepsilon_0$, we obtain
\begin{equation}
    \esp{\mathcal{E}_{\tau^{\alpha,\gamma,N}_{\varepsilon}(R)}}\leq \mathcal{E}_0+\varepsilon^2dN\espa{\tau^{\alpha,\gamma,N}_{\varepsilon}(R)}.
\end{equation}
Since  $\left(\varia[X]{}{i}{\varepsilon},\varia[V]{}{i}{\varepsilon}\right)_i$ are a.s. continuous,  $\mathcal{E}_{\tau^{\alpha,\gamma,N}_{\varepsilon}(R)}=R$. This implies the lower-bound.\\
For the upper-bound in \eqref{estim_temps_sortie_1}, we reconsider \eqref{itoener} at the stopping-time $\tau^{\alpha,\gamma,N}_{\varepsilon}(R)\wedge t$ and derive (recalling that $\phi\le \phi^+$) 
\begin{align*}
    \esp{\mathcal{E}_{\tau^{\alpha,\gamma,N}_{\varepsilon}(R)\land t}}
    &\geq \mathcal{E}_0+\varepsilon^2dN\esp{\tau^{\alpha,\gamma,N}_{\varepsilon}(R)\land t}-2\alpha\esp{\inte{\tau^{\alpha,\gamma,N}_{\varepsilon}(R)\land t} \Ss\norme{\varia[V]{s}{i}{\varepsilon}}^2\di s}\\
    &\quad-\dfrac{\phi^+}{N}\esp{\inte{\tau^{\alpha,\gamma,N}_{\varepsilon}(R)\land t}\Sd \norme{\varia[V]{s}{i}{\varepsilon}-\varia[V]{s}{j}{\varepsilon}}^2\di s}\\
    &\geq \mathcal{E}_0+\varepsilon^2dN\esp{\tau^{\alpha,\gamma,N}_{\varepsilon}(R)\land t}-2\alpha\esp{\inte{\tau^{\alpha,\gamma,N}_{\varepsilon}(R)\land t} \Ss\norme{\varia[V]{s}{i}{\varepsilon}}^2\di s}\\
    &\quad-2\phi^+\esp{\inte{\tau^{\alpha,\gamma,N}_{\varepsilon}(R)\land t}\Ss \norme{\varia[V]{s}{i}{\varepsilon}}^2\di s}+2\phi^+N\esp{\inte{\tau^{\alpha,\gamma,N}_{\varepsilon}(R)\land t}\norme{\overline{V}_s^{N,\varepsilon}}^2\di s}\\
    &\geq \mathcal{E}_0+\varepsilon^2dN\esp{\tau^{\alpha,\gamma,N}_{\varepsilon}(R)\land t}-2\left(\alpha+\phi^+\right)\esp{\inte{\tau^{\alpha,\gamma,N}_{\varepsilon}(R)\land t} \Ss\norme{\varia[V]{s}{i}{\varepsilon}}^2\di s}.
\end{align*}
 Before the stopping-time $\tau^{\alpha,\gamma,N}_{\varepsilon}(R)$, we have $\Ss\norme{\varia[V]{s}{i}{\varepsilon}}^2\leq R$ so we deduce 
\begin{equation}\label{412}
    \esp{\mathcal{E}_{\tau^{\alpha,\gamma,N}_{\varepsilon}(R)\land t}}\geq \mathcal{E}_0+\left[\varepsilon^2dN-2\left(\alpha+\phi^+\right)R\right]\espa{\tau^{\alpha,\gamma,N}_{\varepsilon}(R)\land t}.
\end{equation}
Taking $\varepsilon>\sqrt{2(\phi^++\alpha)R/(Nd)} $ yields the almost-sure finiteness of $\tau^{\alpha,\gamma,N}_{\varepsilon}(R)$ with
\begin{equation*}
    \proba{\ta=+\infty}\leq \dfrac{R-\mathcal{E}_0}{\left[\varepsilon^2dN-2\left(\alpha+\phi^+\right)R\right]t},\quad\forall t>0,
\end{equation*}
letting $t$ going to infinity. Finally, we apply Lebesgue dominated convergence theorem for the left-hand side and monotone convergence theorem  for the right-hand side of \eqref{412} to conclude.
\end{proof}
 \begin{remarque}
\label{RemValid} In the elliptic diffusion case, the mean exit-time from an open stable domain $D$ can be estimated via PDE analysis, yielding further its finiteness for all $\varepsilon>0$ (see \cite[Th.~5.7.3,p.~222]{Dembo_Zeitouni}). The extension of this kind of estimate to kinetic dynamics is not so immediate and is not necessary for our scope. For small $\varepsilon$, we may state the validity of \eqref{cond_petit} with the tail-estimate inequality \eqref{majo_kram} below. Further investigations would be needed to understand what happens when the noise intensity is between $\varepsilon_0$ and $\sqrt{\dfrac{2(\phi^++\alpha)R}{Nd}}$ albeit it is, again, not necessary for our scope.
 \end{remarque}
 We next track sharper estimates than those of Proposition \ref{propestimgro}. To do so, we exploit the properties of the noisy particles system relative to those of its deterministic counterpart and the following LDP, further detailed in Annex \ref{annexe_grd_dev}. Define the following kinetic constraint path set
 \begin{align*}
 &\mathcal{H}_T:=\{(\varphi_1,\varphi_2)\in \espacecont:\,\dot{\varphi}_1=\varphi_2 \text{ and }\varphi_1\in H_2\left([0;T]\right)  \},
 \end{align*}
 $\varphi_1$ denoting the position's coordinates and $\varphi_2$ the velocity's ones. Then, we state the following proposition whose detailed proofs can be found in \cite[Prop.~2.10,~p.~192]{Azencott}, \cite[Th.~7, p.~190]{Priouret} and \cite[Th.~2.1, p.~317]{Baldi_Maurel}.
 \begin{proposition}
    \label{ppdeva}
    The family of processes $\grandeva{}_{\varepsilon>0}$ verifies a LDP , with a good rate function defined by: if $(\varphi_1,\varphi_2)\notin \mathcal{H} _T$, $ \mathcal{J}_{0;T}^N\left(\varphi_1,\varphi_2\right)=+\infty $, and if $(\varphi_1,\varphi_2)\in \mathcal{H}_T $ with $\varphi_1=\left(\varphi_1^i\right)_{i\leq N}$,
    \begin{align*}
        \mathcal{J}_{0;T}^N(\varphi_1,\varphi_2)=\mathcal{J}^N_{0;T}(\varphi_1,\dot{\varphi}_1)=\frac{1}{2}\sum_{i=1}^N\int_0^T\left\lVert \ddot{\varphi}_1^i(t)-b^{i,N}(\varphi_1(t),\dot{\varphi}_1(t))+\alpha\dot{\varphi}_1^i(t)+\gamma\left(\varphi^i_1(t)-x_i^\ast\right)\right\rVert^2\deriv t.
    \end{align*}
\end{proposition}
\begin{remarque}
    To simplify the notation of the rate function, we shall write in the sequel $\mathcal{J}^N_{0;T}(\varphi)$ instead of $\mathcal{J}^N_{0;T}(\varphi,\dot{\varphi})$ on $\mathcal{H}_T$, keeping in mind that $\mathcal{J}^N_{0;T}$ is a function mapping $\espacecont$ to $\reel^+\cup\{+\infty\}$, even though we identify (isometrically) the subspace $E=\{(f_1,f_2),\, f_2=\dot{f}_1\}$
    of $C^0\left([0;T],\mathbb{R}^{2Nd}\right)$ to the space $C^1\left([0;T],\mathbb{R}^{Nd}\right)$.
\end{remarque}

 Let $R>0$ and for every  $0<r\leq R$, we consider the closure and boundary of the exit-domain $\mathcal{O}_N^{\gamma}(r)$ with 
\begin{align*}
     \overline{\mathcal{O}_N^{\gamma}(r)}&=\{ (x_i,v_i)_{i\leq N}\in \reel^{Nd}\times \reel^{Nd}:\,\Ss\lVert v_i\rVert^2+\gamma\Ss\left\lVert x_i-x_i^\ast\right\rVert^2\leq r \},\\
     \partial\mathcal{O}_N^{\gamma}(r)&=\{ (x_i,v_i)_{i\leq N}\in \reel^{Nd}\times \reel^{Nd}:\,\Ss\left\lVert v_i\right\rVert^2+\gamma\Ss\left\lVert x_i-x_i^\ast\right\rVert^2=r\}.
\end{align*} 
We define the quasi-potential associated to the particles system \eqref{Cucker-Smale-bruite} by 
  \begin{equation}
 \label{quasipo}
 \begin{aligned}
  &\mathcal{V}_N\left[(\vv{x_0},\vv{v_0}),(\vv{x}_1,\vv{v}_1)\right]=\underset{\substack{T>0,\\(\varphi,\dot{\varphi})\in K_T}}{\inf}\mathcal{J}^N_{0,T}(\varphi),\text{ where }\\ 
  &K_T=\{(\varphi,\dot{\varphi}) : \varphi \in H_2([0;T]),\big(\varphi(0),\dot{\varphi}(0)\big)=(\vv{x_0},\vv{v_0})\mbox{ and }\big(\varphi(T),\dot{\varphi}(T)\big)=(\vv{x}_1,\vv{v}_1)\}. 
 \end{aligned} 
\end{equation}
This quasi-potential helps us to build the exit-cost from our stable domain $\ensinva{R}$ defined in \eqref{Domaine_inv_classique}, given by 
\begin{equation}
    \label{csortie}
    \mathcal{V}_N(R,\alpha,\gamma):=\underset{(\vv{y},\vv{u})\in \partial\mathcal{O}_N^{\gamma}(R)}{\inf}\mathcal{V}_N\left[\pointeq,(\vv{y},\vv{u})\right].
\end{equation}
Equivalently, the exit-cost writes as
\begin{equation}\label{quasipo-bis}
    \mathcal{V}_N(R,\alpha,\gamma)=\underset{\substack{(\vv{y},\vv{u})\in \partial\mathcal{O}_N^{\gamma}(R),\\T>0,\,(\varphi,\dot{\varphi})\in\mathcal{K}_T^\ast(\vv{y},\vv{u})}}{\inf}\mathcal{J}^N_{0,T}(\varphi),
\end{equation}
where the constraint set $\mathcal{K}_T^\ast$ is defined by  
\begin{equation*}
\mathcal{K}_T^\ast(\vv{y},\vv{u}):=\{\varphi \in H_2\left([0;T]\right)\,:\,\left(\varphi(0),\dot{\varphi}(0)\right)=\vv{z}^*_{\mathrm{eq}},
       \big(\varphi(T),\dot{\varphi}(T)\big)=(\vv{y},\vv{u})\,\text{and}\,\left(\varphi,\dot{\varphi}\right)(]0;T[)\subset\mathcal{O}_N^{\gamma}(R)\}.
\end{equation*}
  When noise is small, one can guess that, most of the time, the noisy particles system shall oscillate around a neighborhood of its equilibrium, remaining in $\ensinva{R}$ before attempting an  excursion outside the latter. Hence, we introduce the following sequences of stopping-times 
 \begin{equation}
    \label{temps_osci}
     \begin{aligned}
      &\zeta_0^\varepsilon=0,\quad \chi_n^\varepsilon=\inf\{t\geq\zeta_n^\varepsilon,\,\left(\vv{X}^{N,\varepsilon}_t,\vv{V}^{N,\varepsilon}_t\right)\in \partial\mathcal{O}_N^{\gamma}(r) \},\\
      &\zeta_n^\varepsilon=\inf\{t>\chi_{n-1}^\varepsilon,\,\left(\vv{X}^{N,\varepsilon}_t,\vv{V}^{N,\varepsilon}_t\right)\in \partial\mathcal{O}_N^{\gamma}(R)\cup \partial\mathcal{O}_N^{\gamma}\left(\dfrac{r}{4}\right) \}.
     \end{aligned}
 \end{equation}
 \vspace{-4mm}
   \begin{figure}[H]
    \centering
    \includegraphics[width=0.6\textwidth]{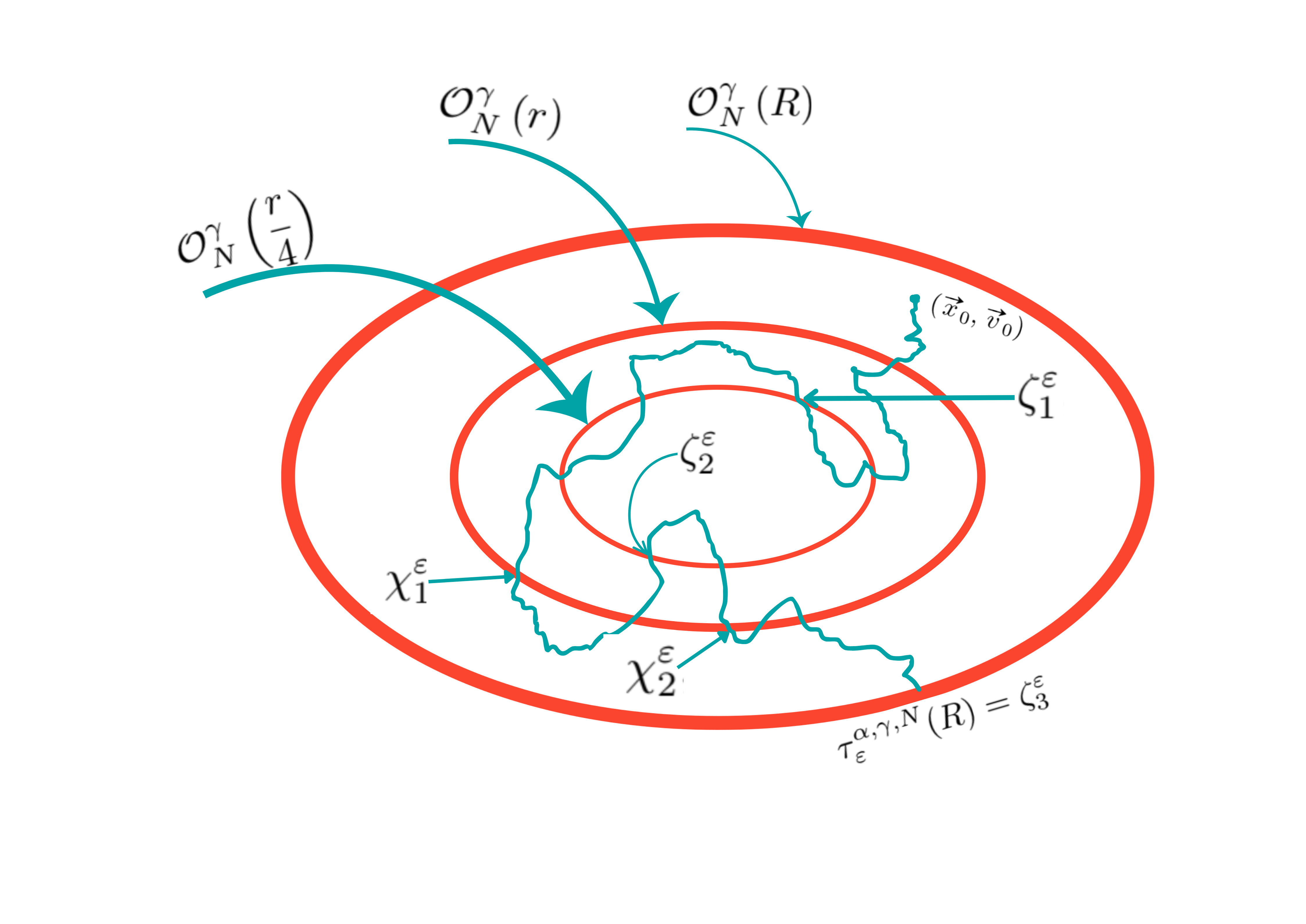}
    \vspace{-12mm}
    \caption{Wandering around the equilibrium before exiting $\ensinva{R}$.}
    \label{fig:oscillations}
\end{figure}
 Let us also define, for every $n\geq 1$, the Markov chain $\vv{\mathcal Z}^{N,\varepsilon}:=\left(\vv{\mathcal Z}_n^{N,\varepsilon}\right)_{n\geq 1}$  with 
 \begin{equation}
    \label{etat_osci}
     \vv{\mathcal Z}_n^{N,\varepsilon}=\left(\vv{X}^{N,\varepsilon}_{\zeta_n^\varepsilon},\vv{V}^{N,\varepsilon}_{\zeta_n^\varepsilon}\right),\qquad n\ge 1,
 \end{equation}
taking values a.s. in $\partial\mathcal{O}_N^{\gamma}(R)\cup \partial\mathcal{O}_N^{\gamma}(r/4)$.
 With these preliminaries, we have the following theorem.
 \begin{theoreme}[Arrhenius' Law]
 \label{loi_arrhenius}
      For every $R,\alpha,\gamma>0$, and all $\vv{z_0}=(\vv{x_0},\vv{v_0})\in \mathcal{O}_N^{\gamma}(R)$, we have
     \begin{equation}\label{arrhenius-asymptotique}
         \lim_{\varepsilon\to 0}\varepsilon^2\log\left(\esper{\left[\tau^{\alpha,\gamma,N}_{\varepsilon}(R)\right]}\right)=\mathcal{V}_N(R,\alpha,\gamma).
     \end{equation}
 \end{theoreme}
 \begin{proof}
  Given $R>0$, we shall show that, for every $\delta>0$, there exists $\varepsilon_0>0$ such that for all $\varepsilon<\varepsilon_0$ we have 
  \begin{equation}
      \label{def_limit}
      \mathcal{V}_N(R,\alpha,\gamma)-\delta<\varepsilon^2\log\left(\esper{\left[\tau^{\alpha,\gamma,N}_{\varepsilon}(R)\right]}\right)<\mathcal{V}_N(R,\alpha,\gamma)+\delta.
  \end{equation}
 \textbf{Upper-bound part of \eqref{def_limit}.} Let $\delta>0$ and $ (\vv{x_0},\vv{v_0})\in \mathcal{O}_N^{\gamma}(R)$. By Lemma \ref{sortiefini} (Appendix \ref{annexe_grd_dev}), there exists a deterministic time  $T(\delta)$ and a real number $\varepsilon_0>0$ not depending on $(\vv{x_0},\vv{v_0})$ such that for every $\varepsilon<\varepsilon_0$
  \begin{equation}
    \label{estim_petit}
     \underset{(\vv{x_0},\vv{v_0})\in \mathcal{O}_N^{\gamma}(R)}{\inf}\probabi\left[\tau^{\alpha,\gamma,N}_{\varepsilon}(R)<T(\delta)\right]\geq \dfrac{1}{2} \exp{\Bigg(-\dfrac{1}{\varepsilon^2}\bigg(\mathcal{V}_N(R,\alpha,\gamma)+\dfrac{\delta}{2}\bigg)\Bigg)}.
 \end{equation}
 Since $\tau^{\alpha,\gamma,N}_{\varepsilon}(R)$ is a positive random variable, we have 
 \begin{align}
    \label{estim_espe}
     \esper{\left[\tau^{\alpha,\gamma,N}_{\varepsilon}(R)\right]}&=\inte{+\infty}\probabi\left(\tau^{\alpha,\gamma,N}_{\varepsilon}(R)>t\right)\di t=\sum_{n=0}^{+\infty}\int_{nT(\delta)}^{(n+1)T(\delta)}\probabi\left(\tau^{\alpha,\gamma,N}_{\varepsilon}(R)>t\right)\di t\notag\\
     &\leq T(\delta) \sum_{n=0}^{+\infty}\probabi\left(\tau^{\alpha,\gamma,N}_{\varepsilon}(R)>nT(\delta)\right).
 \end{align}
 For any $n\in \mathbb{N}$, applying the Markov property with respect to the $\sigma$-algebra $\mathcal{F}_{n T(\delta)}$ immediately yields
 \begin{align*}
     \probabi\left[\tau^{\alpha,\gamma,N}_{\varepsilon}(R)>(n+1)T(\delta)\right]
     &=\esper{\left(\mathbbm{1}_{\{\tau^{\alpha,\gamma,N}_{\varepsilon}(R)>(n+1)T(\delta)\}\cap\{\tau^{\alpha,\gamma,N}_{\varepsilon}(R)>nT(\delta)\}}\right)}\\
     &=\esper{\left(\mathbbm{1}_{\{\tau^{\alpha,\gamma,N}_{\varepsilon}(R)>nT(\delta)\}}\mathbb{P}_{\vv{Z}^{N,\varepsilon}_{nT(\delta)}}\left[\tau^{\alpha,\gamma,N}_{\varepsilon}(R)>T(\delta)\right]\right)}\\
     &\leq \probabi{\left[\tau^{\alpha,\gamma,N}_{\varepsilon}(R)>nT(\delta)\right]}\underset{\vv{z}\in \mathcal{O}_N^{\gamma}(R)}{\max}\mathbb{P}_{\vv{z}}\left[\tau^{\alpha,\gamma,N}_{\varepsilon}(R)>T(\delta)\right].
 \end{align*}
 By induction, it follows 
 \begin{align*}
     \probabi\left[\tau^{\alpha,\gamma,N}_{\varepsilon}(R)>nT(\delta)\right]&\leq\left(\underset{(\vv{x},\vv{y})\in \mathcal{O}_N^{\gamma}(R)}{\max}\mathbb{P}_{(\vv{x},\vv{y})}\left[\tau^{\alpha,\gamma,N}_{\varepsilon}(R)>T(\delta)\right]\right)^n\\
     &=\left(1-\underset{(\vv{x},\vv{y})\in \mathcal{O}_N^{\gamma}(R)}{\inf}\mathbb{P}_{(\vv{x},\vv{y})}\left[\tau^{\alpha,\gamma,N}_{\varepsilon}(R)\leq T(\delta)\right]\right)^n.
 \end{align*}
 Hence, combined with \eqref{estim_petit}, for $\varepsilon<\varepsilon_0$, the inequality \eqref{estim_espe} yields 
 \begin{align}
    \label{majo_kram}
     \esper{\left[\tau^{\alpha,\gamma,N}_{\varepsilon}(R)\right]}&\leq T(\delta) \sum_{n=0}^{+\infty}\left(1-\underset{(\vv{x},\vv{y})\in \mathcal{O}_N^{\gamma}(R)}{\inf}\mathbb{P}_{(\vv{x},\vv{y})}\left[\tau^{\alpha,\gamma,N}_{\varepsilon}(R)\leq T(\delta)\right]\right)^n\notag\\
     &\leq 2T(\delta)\exp{\Bigg(\dfrac{1}{\varepsilon^2}\bigg(\mathcal{V}_N(R,\alpha,\gamma)+\dfrac{\delta}{2}\bigg)\Bigg)}.
 \end{align}
 Thus, there exists $\varepsilon_1>0$ such that for every $\varepsilon<\varepsilon_1$, we have the upper-bound in \eqref{def_limit}
 \begin{align*}
     \varepsilon^2\log\left(\esper{\left[\tau^{\alpha,\gamma,N}_{\varepsilon}(R)\right]}\right)<\mathcal{V}_N(R,\alpha,\gamma)+\delta.
 \end{align*}
 \textbf{Lower-bound part of \eqref{def_limit}.} \textbf{Step 1.} Let us consider the Markov chain $\vv{\mathcal Z}^{N,\varepsilon}$ as in \eqref{etat_osci} and define
 \begin{equation}
     \varrho_\varepsilon:=\inf\{n\in \mathbb{N},\, \vv{\mathcal Z}_n^{N,\varepsilon}\in\partial\mathcal{O}_N^{\gamma}(R)\}.
 \end{equation}
 corresponding to the first time $\vv{\mathcal Z}^{N,\varepsilon}$ hits $\partial\mathcal{O}_N^{\gamma}(R)$. Let  $\initia\in \partial\mathcal{O}_N^{\gamma}(r/4)$ and $h<\mathcal{V}_N(R,\alpha,\gamma)$. Applying Lemma~\ref{ev_impossi} (Appendix \ref{annexe_grd_dev}) yields, for $\varepsilon<\varepsilon_0$ with $\varepsilon_0>0$ independent of $ \initia$, the following upper-bound
   \begin{align}
    \label{sortie_v1}
     \probabi{\Big(\vv{\mathcal Z}_1^{N,\varepsilon}\in\partial\mathcal{O}_N^{\gamma}(R)\Big)}&=\probabi{\left(\varrho_\varepsilon=1\right)}=\probabi{\left(\zeta_1^\varepsilon=\tau^{\alpha,\gamma,N}_{\varepsilon}(R)\right)}\leq \exp{\left(-\frac{1}{\varepsilon^2}\left(\mathcal{V}_N(R,\alpha,\gamma)-h\right)\right)}.
 \end{align}
 For $n\in \mathbb{N}$ and $z:=(z_n)_n\subset\reel^{2dN}$, define the functional 
  \begin{equation*}
      \varrho(z):=\inf\{n\in \mathbb{N}^\ast,\, z_n\in  \partial\mathcal{O}_N^{\gamma}(R)\}.
  \end{equation*}
 Recalling the time-shift operator $\theta_n((\vv{z}(t))_t)=(\vv{z}(t+n))_t$, the (strong) Markov property conditioning by $\mathcal{F}_{\zeta_n^\varepsilon}$ gives
  \begin{align*}
       \probabi{(\varrho_\varepsilon>n+1)}&=\esper{\left(\mathbbm{1}_{\{\varrho_\varepsilon>n+1\}}\mathbbm{1}_{\{\varrho_\varepsilon>n\}}\right)}=\esper{\left(\mathbbm{1}_{\{n+\varrho\left(\theta_n(\vv{\mathcal Z}^{N,\varepsilon})\right)>n+1\}}\mathbbm{1}_{\{\varrho_\varepsilon>n\}}\right)}\\
 &=\esper{\left(\mathbbm{1}_{\{\varrho_\varepsilon>n\}}\esper{\left(\mathbbm{1}_{\{\varrho\left(\theta_n(\vv{\mathcal Z}^{N,\varepsilon})\right)>1\}}\middle\vert\mathcal{F}_{\zeta_n^\varepsilon}\right)}\right)}\\
 &=\esper{\left(\mathbbm{1}_{\{\varrho_\varepsilon>n\}}\mathbb{E}_{\vv{\mathcal  Z}_n^{N,\varepsilon}}\left(\mathbbm{1}_{\{\varrho_\varepsilon>1\}}\right)\right)}.
  \end{align*}
  On the event $\{\varrho_\varepsilon>n\}$, $\vv{\mathcal Z}_n^{N,\varepsilon}\notin \partial\mathcal{O}_N^{\gamma}(R)$ so  $\vv{\mathcal Z}_n^{N,\varepsilon}\in \partial\mathcal{O}_N^{\gamma}\left(r/4\right)$. Thus, we deduce 
  \begin{equation}
      \probabi{(\varrho_\varepsilon>n+1)}\geq \left(\underset{\initias\in \partial\mathcal{O}_N^{\gamma}\left(r/4\right)}{\inf}\probabis{(\varrho_\varepsilon>1)}\right)\times \probabi{(\varrho_\varepsilon>n)}.
  \end{equation}
For all  $\initia\in \partial\mathcal{O}_N^{\gamma}(\frac{r}{4})$, by induction and thanks to \eqref{sortie_v1}, for every $\varepsilon<\varepsilon_0$, we get
\begin{equation}
    \label{mino_ten_sort}
    \probabi{(\varrho_\varepsilon>n)}\geq \Bigg[1-\exp{\left(-\frac{1}{\varepsilon^2}\bigg(\mathcal{V}_N(R,\alpha,\gamma)-h\bigg)\right)}\Bigg]^n. 
\end{equation}
 \textbf{Step 2.} Given the identities 
 \begin{align*}
     \tau^{\alpha,\gamma,N}_{\varepsilon}(R)=\zeta^\varepsilon_{\varrho_\varepsilon}=\sum_{n=1}^{\varrho_\varepsilon} (\zeta_{n}^\varepsilon-\zeta_{n-1}^\varepsilon)=\sum_{n=1}^{+\infty}\mathbbm{1}_{\{\varrho_\varepsilon\geq n\}}(\zeta_{n}^\varepsilon-\zeta_{n-1}^\varepsilon),
 \end{align*}
applying the Beppo-Levi theorem yields 
 \begin{align*}
     \esper{\left[\tau^{\alpha,\gamma,N}_{\varepsilon}(R)\right]}=\sum_{n=1}^{+\infty}\esper{\left[\mathbbm{1}_{\{\varrho_\varepsilon\geq n\}}(\zeta_{n}^\varepsilon-\zeta_{n-1}^\varepsilon)\right]}.
 \end{align*}
 By definition, we have  $\chi_{n-1}^\varepsilon\geq \zeta_{n-1}^\varepsilon$. By applying the Markov property at the stopping-time $\chi_{n-1}^\varepsilon$, for every  $\vv{z_0}=\initia\in \partial\mathcal{O}_N^{\gamma}\left(r/4\right)$, we obtain
 \begin{equation}
    \label{mino_esper}
     \esper{\left[\tau^{\alpha,\gamma,N}_{\varepsilon}(R)\right]}\geq\underset{\initias\in \partial\mathcal{O}_N^{\gamma}(r)}{\inf}\esperd{\left[\tau_{\mathcal{A}(r)}\left(\vv{X}^{N,\varepsilon},\vv{V}^{N,\varepsilon}\right)\right]}\sum_{n=1}^{+\infty} \probabi{(\varrho_\varepsilon\geq n)},
 \end{equation}
  where $\tau_A$ is as in \eqref{fonctionnelle}  and  $\mathcal{A}(r):=\partial\mathcal{O}_N^{\gamma}(R)\cup \partial\mathcal{O}_N^{\gamma}\left(r/4\right)$. For simplicity we denote $\tau_{\mathcal{A}(r)}\left(\vv{X}^{N,\varepsilon},\vv{V}^{N,\varepsilon}\right)=: \tau_{\mathcal{A}(r)}$.

 Applying Lemma \ref{temp_min}, one can find $\varepsilon_1>0$ and a deterministic time $T(r)$ such that, for every $\initia\in\ensinvab{r}$  and all $\varepsilon<\varepsilon_1$, we have
 \begin{equation}
     \probabi{\left[\tau_{\mathcal{A}(r)}\geq T(r)\right]}\geq\dfrac{1}{2}.
 \end{equation}
 Moreover, using the Markov inequality, for $\varepsilon<\varepsilon_1$ and all $\initia\in\ensinvab{r}$, we deduce the lower-bound $$\esper{\left[\tau_{\mathcal{A}(r)}\right]}\geq T(r) \probabi{\left[\tau_{\mathcal{A}(r)}\geq T(r)\right]}\geq\dfrac{T(r)}{2}.$$
 Hence, for all $\initia\in \partial\mathcal{O}_N^{\gamma}\left(\dfrac{r}{4}\right)$ and all $\varepsilon<\min(\varepsilon_0,\varepsilon_1)$, thanks to \eqref{mino_ten_sort}, inequality \eqref{mino_esper} implies
 \begin{equation}
    \label{estim_transi}
     \esper{\left[\tau^{\alpha,\gamma,N}_{\varepsilon}(R)\right]}\geq\dfrac{T(r)}{2}\exp{\left(\frac{1}{\varepsilon^2}\bigg[\mathcal{V}_N(R,\alpha,\gamma)-h\bigg]\right)},
 \end{equation}
so that for $\varepsilon<\varepsilon_2$ with $\varepsilon_2>0$, we deduce 
\begin{equation}
    \label{mino_Arrhenius}
    \varepsilon^2\log\left(\esper{\left[\tau^{\alpha,\gamma,N}_{\varepsilon}(R)\right]}\right)>\mathcal{V}_N(R,\alpha,\gamma)-2h.
\end{equation}
It suffices to take $h=\delta/2$ to conclude, yielding the  Arrhenius' law when the system \eqref{Cucker-Smale-bruite} starts with $\vv{z_0}\in\partial\mathcal{O}_N^{\gamma}\left(r/4\right)$.\\
\textbf{Step 3.} It remains to extend \eqref{mino_Arrhenius} to all $\vv{z}_0=\initia\in \mathcal{O}_N^{\gamma}(R)$. To this aim, observe
\begin{align}
           \esper{\left[\tau^{\alpha,\gamma,N}_{\varepsilon}(R)\right]}&\geq\esper{\left(\mathbbm{1}_{\{\tau^{\alpha,\gamma,N}_{\varepsilon}(R)>\zeta_1^\varepsilon\}}\tau^{\alpha,\gamma,N}_{\varepsilon}(R)\right)}\notag\\
           &\geq \esper{\left[\mathbbm{1}_{\{\tau^{\alpha,\gamma,N}_{\varepsilon}(R)>\zeta_1^\varepsilon\}}\tau_{\partial\mathcal{O}_N^{\gamma}(R)}\left(\theta_{\zeta_1^\varepsilon}\grandeva{}\right)\right]},
\end{align}
where $\theta$ is the shift-time operator and $\tau_A$  for a Borel set $A$ is the functional defined in \eqref{fonctionnelle}. Applying the strong Markov property by considering the stopping-time $\zeta_1^\varepsilon$, we obtain
\begin{equation}
    \label{estim_temps_general}
    \esper{\left[\tau^{\alpha,\gamma,N}_{\varepsilon}(R)\right]}\geq \probabi{\bigg(\tau^{\alpha,\gamma,N}_{\varepsilon}(R)>\zeta_1^\varepsilon\bigg)}\underset{\initias\in \partial\mathcal{O}_N^{\gamma}\left(r/4\right)}{\inf}\esperd{\left[\tau^{\alpha,\gamma,N}_{\varepsilon}(R)\right]}.
\end{equation}
Since the noisy particles system is close to its deterministic counterpart with high probability, thanks to Lemma \ref{exit_b_ball} for $r$ very small, one can find $\varepsilon_1>0$ (depending on initial conditions) such that for every $\varepsilon<\varepsilon_1$ $$\probabi{\left(\tau^{\alpha,\gamma,N}_{\varepsilon}(R)>\zeta_1^\varepsilon\right)}\geq \dfrac{1}{2}.$$ 
Merging the estimates \eqref{estim_transi} and \eqref{estim_temps_general}, we deduce the validity of inequality \eqref{mino_Arrhenius} for $\varepsilon<\min(\varepsilon_0,\varepsilon_1)$ when the system \eqref{Cucker-Smale-bruite} starts in $\mathcal{O}_N^{\gamma}(R)$, which completes the proof of the Arrhenius' law.
 \end{proof}
 We now state our second main result, the Kramers's type law for the system $\grandeva{}$.
\begin{theoreme}[Kramers' type law]
\label{l_kram}
    For all $\initia\in {\mathcal{O}}_N^{\alpha,\gamma}(R)$ and for all $\delta>0$, we have the limit
    \begin{equation*}
        \underset{\varepsilon\to 0}{\lim}\,\probabi{\left(\Norme{\varepsilon^2\log\Big(\tau^{\alpha,\gamma,N}_{\varepsilon}(R)\Big)-\mathcal{V}_N(R,\alpha,\gamma)}<\delta\right)}=1.
    \end{equation*}
\end{theoreme}
\begin{proof}
\textbf{Upper-tail.}  Let $\delta >0$ and $\initia\in \ensinva{R}$.   Applying the Markov inequality,  
\begin{align*}
    \probabi{\left[\tau^{\alpha,\gamma,N}_{\varepsilon}(R)\geq \exp{\left(\dfrac{\mathcal{V}_N(R,\alpha,\gamma)+\delta}{\varepsilon^2}\right)}\right]}&\leq\exp{\left(-\dfrac{\mathcal{V}_N(R,\alpha,\gamma)+\delta}{\varepsilon^2}\right)}\esper{\left[\tau^{\alpha,\gamma,N}_{\varepsilon}(R)\right]}\\
   & \leq  2 T(\delta)\exp{\left(-\dfrac{\delta}{2\varepsilon^2}\right)},
\end{align*}
where the last line follows from \eqref{majo_kram} and choosing $0<\varepsilon<\varepsilon_0$. Thus
\begin{equation*}
    \underset{\varepsilon\to 0}{\lim}\,
    \probabi{\left(\varepsilon^2\log\Big(\tau^{\alpha,\gamma,N}_{\varepsilon}(R)\Big)-\mathcal{V}_N(R,\alpha,\gamma)<\delta\right)}=1.
\end{equation*}
\textbf{Lower-tail}. Our argument is inspired by \cite{Dembo_Zeitouni}. Without loss of generality, let us suppose that $\delta<\mathcal{V}_N(R,\alpha,\gamma)$. Let $k\in \mathbb{N}^\ast$, $T_0>0$ and $r$ to be chosen in the sequel. Let us take $\initia\in \ensinva{R}$. We have
\begin{align*}
    \probabi{\left(\ta\leq kT_0\right)}&\leq \probabi{\left(\{\zeta_1^\varepsilon<\ta\}\cap\{\ta\leq kT_0\}\right)}\\
    &\quad+\probabi{\left(\zeta_1^\varepsilon\geq\ta\right)}.
\end{align*}
The quantity $\probabi{\left(\zeta_1^\varepsilon\geq\ta\right)}$ goes to $0$ when $\varepsilon$ goes to $0$ for $r$ sufficiently small, thanks to Lemma~\ref{exit_b_ball}. We shall show that for a certain integer $k$,
\begin{equation}
    \mathcal{P}:=\probabi{\left(\{\zeta_1^\varepsilon<\ta\}\cap\{\ta\leq kT_0\}\right)}\underset{\varepsilon\to 0}{\longrightarrow} 0.
\end{equation}
 Applying the strong Markov property and considering the stopping-time $\zeta_1^\varepsilon$, we get 
\begin{align} \label{ineq_usef}
    \mathcal{P}&\leq  \probabi{\left(\{\zeta_1^\varepsilon<\ta\}\cap\{\tau_{\ensinvab{R}}\left(\theta_{\zeta_1^\varepsilon}\grandeva{}\right)\leq kT_0\}\right)}\notag\\
    &\leq\probabi{\left(\zeta_1^\varepsilon<\ta\right)}\underset{\initias\in \ensinvab{\frac{r}{4}}}{\max}\probabis{\left(\ta\leq kT_0\right)}\notag\\
    &\leq \underset{\initias\in \ensinvab{\frac{r}{4}} }{\max}\probabis{\left(\ta\leq kT_0\right)},
\end{align}
with $\tau_A$ as in \eqref{fonctionnelle}.
Henceforth, let us now take $\initia\in \ensinvab{\dfrac{r}{4}}$ and define the following events (recalling \eqref{temps_osci} for the setting of $\zeta_n^\varepsilon$ and $\chi_{n}^\varepsilon$): 
\begin{align}
    \mathcal{E}_k:=\bigcup_{j=1}^k\{\ta =\zeta_j^\varepsilon\}\quad\mathrm{ and}\quad\mathcal{B}_k(T_0):=\{\underset{j\in [\![1;k]\!]}{\min}\left(\zeta_j^\varepsilon-\chi_{j-1}^\varepsilon\right)\leq T_0\}.
\end{align}
Let $\omega \in \prescript{}{}{\mathcal{E}_k^c}\cap\prescript{}{}{\mathcal{B}_k^c(T_0)}$. Recalling that  $\chi_j^\varepsilon\geq \zeta_j^\varepsilon$ for all $j\in \mathbb{N}^\ast$, we obtain the lower-bound
\begin{align*}
    \ta(\omega)>\zeta_k^\varepsilon(\omega)=\sum_{j=1}^k\left(\zeta_j^\varepsilon(\omega)-\zeta_{j-1}^\varepsilon(\omega)\right)> k\underset{j\in [\![1;k]\!]}{\min}\left(\zeta_j^\varepsilon(\omega)-\chi_{j-1}^\varepsilon(\omega)\right)>kT_0.
\end{align*}
It follows
\begin{align}
    \label{sort_b_T}
    \probabi{\left(\ta\leq kT_0\right)}&\leq \probabi{\left(\mathcal{E}_k\right)}+\probabi{\left(\mathcal{B}_k(T_0)\right)}\notag\\
    &\leq \sum_{j=1}^k\probabi{\left(\ta= \zeta_j^\varepsilon\right)}+\sum_{j=1}^k\probabi{\left(\zeta_j^\varepsilon-\chi_{j-1}^\varepsilon\leq T_0\right)}.
\end{align}
Setting next $\mathcal{A}(r):=\partial\mathcal{O}_N^{\gamma}(R)\cup \partial\mathcal{O}_N^{\gamma}\left(\dfrac{r}{4}\right)$ and $\tau_{\mathcal{A}(r)}\left(\vv{X}^{N,\varepsilon},\vv{V}^{N,\varepsilon}\right)=: \tau_{\mathcal{A}(r)}$, by the strong Markov property, 
\begin{align*}
\probabi{\left(\zeta_j^\varepsilon-\chi_{j-1}^\varepsilon\leq T_0\right)}&=\probabi{\left(\tau_{\mathcal{A}(r)}\left(\theta_{\chi_{j-1}^\varepsilon }\grandeva{}\right)\leq T_0\right)}\\
&\leq \underset{\initias\in \ensinvab{r}}{\max}\probabis{\left(\tau_{\mathcal{A}(r)}\leq T_0\right)}.
\end{align*}
Moreover, for every $\initias\in \ensinvab{r}$, we have  \begin{equation*}
    \{\tau_{\mathcal{A}(r)}\leq T_0\}\subset \{\supre{T_0}\sqrt{\norme{\vv{X}^{N,\varepsilon}_t-\vv{x}}^2+\norme{\vv{V}^{N,\varepsilon}_t-\vv{v}}^2}\geq \mu_r\},
\end{equation*} 
with $\mu_r:=\mu_r^2\land\mu_r^1$ where the positive real numbers $\mu_r^1$ and $\mu_r^2$ are defined by
\begin{align*}
    &\mu_r^1:=\inf\{\norme{x-y}: (x,y)\in \ensinvab{r}\times \ensinvab{R}\},\quad\mu_r^2:=\inf\{\norme{x-y}: (x,y)\in \ensinvab{r}\times \ensinvab{\frac{r}{4}}\}.
\end{align*}
We deduce that
\begin{align*}
    \probabis{\left(\tau_{\mathcal{A}(r)}\leq T_0\right)}&\leq \probabis{\left(\supre{T_0}\norme{\vv{X}^{N,\varepsilon}_t-\vv{x}}\geq\dfrac{\mu_r}{2}\right)}+\probabis{\left(\supre{T_0}\norme{\vv{V}^{N,\varepsilon}_t-\vv{v}}\geq\dfrac{\mu_r}{2}\right)}.
\end{align*}
By virtue of Lemma \ref{proche_ini} with $c=\mathcal{V}_N(R,\alpha,\gamma)-\dfrac{\delta}{2}$ and Remark~\ref{extent_impo}, one can find a deterministic time $T_0:=T(\mu_r,c)>0$ for $r$ sufficiently small and depending on $\delta$, a real number $\varepsilon_1>0$ not depending on initial conditions such that for every $\initia \in \ensinvab{\dfrac{r}{4}}$ and every $\varepsilon<\varepsilon_1$ , the inequality \eqref{sort_b_T} yields 
\begin{align*}
    \probabi{\left(\ta\leq kT_0\right)}&\leq \left[1+8Nd\right]k\exp{\left(-\dfrac{1}{\varepsilon^2}\left(\mathcal{V}_N(R,\alpha,\gamma)-\dfrac{\delta}{2}\right)\right)}.
\end{align*}
Taking $k:=\left\lfloor\exp{\left(\dfrac{1}{\varepsilon^2}\left(\mathcal{V}_N(R,\alpha,\gamma)-\delta\right)\right)}\right\rfloor+1$, for $\varepsilon\leq \varepsilon_2$, the positive real number $\varepsilon_2$ depending on the initial conditions, the inequality \eqref{ineq_usef} induces 
\begin{align*}
    \mathcal{P}\leq \left[1+8Nd\right]\left(\exp{\left(-\dfrac{\delta}{2\varepsilon^2}\right)}+\exp{\left(-\dfrac{1}{\varepsilon^2}\left[\mathcal{V}_N(R,\alpha,\gamma)-\dfrac{\delta}{2}\right]\right)}\right).
\end{align*}
We complete the proof by noticing that  $$\probabi{\left(\ta\leq \exp{\left(\dfrac{1}{\varepsilon^2}\left(\mathcal{V}_N(R,\alpha,\gamma)-\delta\right)\right)} \right)}\leq\probabi{\left(\ta\leq k T_0\right)}.$$
\end{proof}

\section{Numerical illustrations}\label{sec5}
In this section, we illustrate our findings, mainly the Arrhenius' law given in Theorem \ref{loi_arrhenius}, by simulations. We choose to investigate cases in which the particle target is common, namely $x_i^\ast=x^*$ for all $i$. In this case, it is possible to improve the upper-bound given at Lemma \ref{posi_cou}, for the exit-cost $\mathcal{V}_N(R,\alpha,\gamma)$, making it only dependent on the friction coefficient $\alpha$ and the positive constant $R$. In fact, it is easy to derive the bound $\mathcal{V}_N(R,\alpha,\gamma)\leq \widehat{\mathcal{V}}_N(R,\alpha,\gamma)$, considering that $\widehat{\mathcal{V}}_N(R,\alpha,\gamma)$ is the exit-cost under the restricted constraint set \eqref{quasipo-bis} where the $i$-th components of $\varphi$ are indistinguishable paths (i.e. $\varphi:=(\phi,\cdots,\phi)$ with $\phi\in C^1([0;T],\reel^{d})$); that is 
\begin{equation*}
    \widehat{\mathcal{V}}_N(R,\alpha,\gamma):=N\underset{\substack {\left(x,v\right)\in \mathcal{D}(R/N)\\T>0,\,\phi\in \widehat{K}_T}}{\inf}\dfrac{1}{2}\inte{T}\!\norme{\ddot{\phi}(t)+\alpha\dot{\phi}(t)+\gamma(\phi(t)-x^\ast)}^2\di t,
\end{equation*}
where $\widehat{K}_T:=  \{\phi \in H_2([0;T],\reel^d):\big(\phi(0),\dot{\phi}(0)\big)=(x^\ast,0_{\reel^d}),\big(\phi(T),\,\dot{\phi}(T)\big)=\left(x,v\right)\}$ and where the end-point $(\vv{y},\vv{u})\in\partial \mathcal O^\gamma_N(R)$ reduces to 
$$
\mathcal{D}(R/N):=\{\left(x,v\right)\in\reel^{2d}: \norme{v}^2+\gamma\norme{x-x^\ast}^2=\frac{R}{N}\}.
$$  
As the corresponding value function further reduces to  
\begin{equation}
    \tag{$\mathcal{P}$}
    \underset{\substack {T>0,\,\phi\in \widehat{K}_T}}{\inf}\inte{T}\norme{\ddot{\phi}(t)+\alpha\dot{\phi}(t)+\gamma(\phi(t)-x^\ast)}^2\di t=\alpha\norme{v}^2+\gamma\alpha\norme{x-x^\ast}^2,
\end{equation}
Thus, $\widehat{\mathcal{V}}_N(R,\alpha,\gamma)=\alpha R$ and Lemma~\ref{posi_cou} becomes $\dfrac{\alpha R}{2}<\mathcal{V}_N(R,\alpha,\gamma)\leq\alpha R $. 

With Theorem \ref{loi_arrhenius}, we have shown that, in the small-noise regime, $\esper{\left[\tau^{\alpha,\gamma,N}_{\varepsilon}(R)\right]}\approx e^{\frac{\mathcal{V}_N(R,\alpha,\gamma)}{\varepsilon^2}}$, so that the importance of the refined bound lies in making the exit numerically feasible, since we employ a basic (and naive) Monte-Carlo approach to estimate the latter expectation by emphasizing the statistical relevance of our estimations. In addition, the Arrhenius' law is given asymptotically so that we do not clearly grasp under which threshold it stands. Thus, we overall take reasonable small noises for our illustration mostly under the bound $\sqrt{\dfrac{2(\lambda\phi^++\alpha)R}{Nd}}$ (see Proposition \ref{propestimgro}) in order to avoid a fast exit from the domain $\ensinva{R}$ . Finally, we use an Euler-Maruyama  scheme $\grandevae{}$ to approximate $\grandeva{}$, the discretization step being  taken homogeneous with $h=0.01$. Thus, the exit-time is approximated by  
\begin{equation*}
    \ta\approx\tau_\varepsilon^h(R):=\inf\{n\in \mathbb{N},\,\grandevae{nh}\notin \mathcal{O}_N^{\gamma}(R) \}\times h.
\end{equation*}
and its mean value by simple Monte-Carlo samples. 
Note that the discretization scheme may be too rough in view of the fact that an exit occurs at very large time. To overcome this limitation, we verify \textit{a posteriori} its quality by choosing an acceptable error threshold $\varepsilon_a$, namely: $$\norme{\vv{V}_{\tau_\varepsilon^h(R)}^{N,\varepsilon,h}}^2+\gamma \norme{\vv{X}_{\tau_\varepsilon^h(R)}^{N,\varepsilon,h}-\vv{x}_{\mathrm{eq}}^\ast}^2-R\leq\varepsilon_a.$$

\begin{figure}[H]
    \centering
    \includegraphics[width=0.975\textwidth]{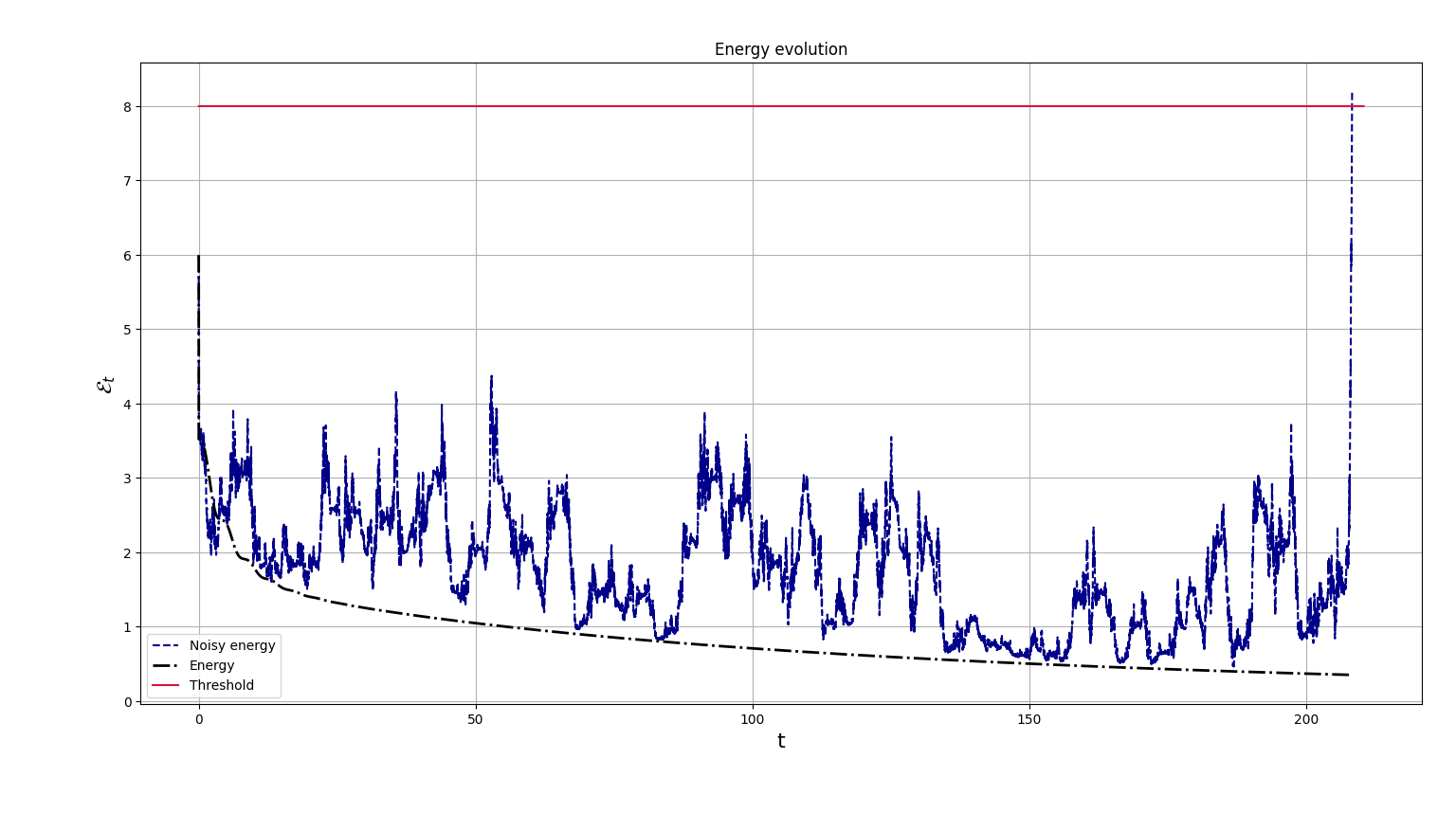}
    \vspace{-8mm}
    \captionsetup{font=small}
    \caption{\footnotesize{In continuous line: Sample of the energy $\mathcal{E}^\varepsilon$ for $R=8, \alpha=0.2,\gamma=0.6,d=1,N=2,\lambda=300,\kappa=1,\beta=0.5,x^\ast=-5$ and $\varepsilon=0.4$. $\mathcal{E}^\varepsilon$ overcomes the threshold $R=8$ in a long time and most of the time is near $0$. In dash-dot: Plot of $\mathcal{E}^0$.}}
    \label{fig:oscillations_R_8}
\end{figure}

\begin{figure}[H]
    \centering
    \includegraphics[width=0.975\textwidth]{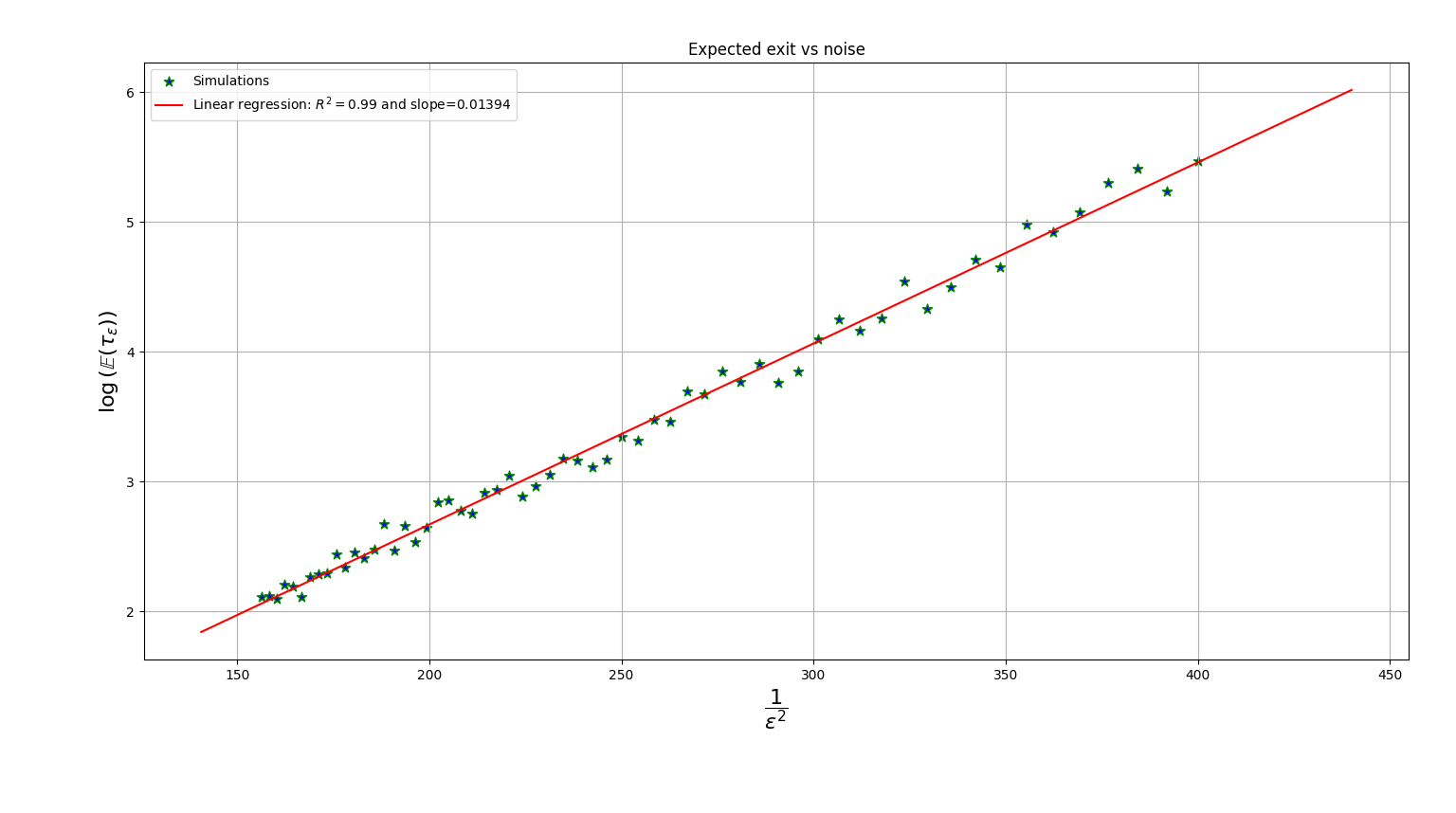}
    \vspace{-8mm}
    \captionsetup{font=small}
    \caption{\footnotesize{Numerical simulations of $\log\left(\esp{\tau_\varepsilon}\right)$ on a $1/\varepsilon^2$ scale for $\varepsilon\in [0.05;0.08]$. The parameters are set to $R=1/8, \alpha=0.2,\gamma=0.2,d=3,N=5,\lambda=5,\kappa=1,\beta=0.5$ and $x^\ast=(-5;7;9)$. $100$ Monte-Carlo simulations have been performed.}}
    \label{fig:oscillations_R_5}
\end{figure}

\noindent
Figure \ref{fig:oscillations_R_8} shows the evolution of $ \mathcal{E}^\varepsilon_t =\norme{\vv{V}_t^{N,\varepsilon}}^2+\gamma \norme{\vv{X}_t^{N,\varepsilon}-\vv{x}_{\mathrm{eq}}^\ast}^2$ and $\mathcal{E}^0_t=\mathcal{U}^{\gamma,N}_{\mathrm{tot}}(t)$ (see \eqref{ener_u}) for $R=8$ . Most of the time, $\mathcal{E}_t^\varepsilon$ wanders around $0$ attempting  to overcome the threshold $R$. At each exit failure, the particles system forgets its past and its energy is attracted by $0$. Figure \ref{fig:oscillations_R_8} also highlights the proximity between  $\mathcal{E}^\varepsilon_t$ and $\mathcal{E}^0_t$ at least during the first stage of the motion, corroborating the proximity  between the noisy system and its deterministic counterpart at least in probability.
Arrhenius' law is illustrated in Figure \ref{fig:oscillations_R_5} with $R\!=\!1/{8}$ and $\alpha\!=\!0.2$. The figure displays $\log\left(\esp{\tau_\varepsilon}\right)$ versus $1/\varepsilon^2$ for $\varepsilon\in [0.05;0.08]$, as well as the corresponding linear regression. We obtain a coefficient of determination equal to $0.99$ and the slope gives us an estimation of the exit-cost corroborating the refined bounds obtained previously.

\vfill

\appendix
\numberwithin{equation}{section}
\section{Large deviations principle}\label{annexe_grd_dev}
\numberwithin{equation}{subsection}
\newtheorem{theoremea}{Theorem}[subsection]
\newtheorem{remarquea}[theoremea]{Remark}
\newtheorem{propositiona}[theoremea]{Proposition}
\newtheorem{definitiona}[theoremea]{Definition}
\newtheorem{lemmea}[theoremea]{Lemma}
\newtheorem{corollairea}[theoremea]{Corollary}
\newtheorem{hypothesea}[theoremea]{Hypothesis}
\newtheorem{hypothesesa}[theoremea]{Hypotheses}
\subsection{Classical results}\leavevmode\par\vspace{0.5em}

We recall (\cite{Freidlin_Wentzell,Dembo_Zeitouni}) that the family of processes $\{\grandeva{}\}_{\varepsilon>0}$ satisfies a LDP (large deviation principle) if
there exists a lower semi-continuous functional $\mathcal{J}^N_{0,T} :C^0_{\vv{z}}([0;T],\mathbb{R}^{2Nd})\longmapsto \mathbb{R}^+\cup\{+\infty\}$ (called rate function) such that  
    \begin{align*}  
    &\bullet\text{for any open subset }O\text{ of }C^0_{\vv{z}}\left([0;T],\mathbb{R}^{2Nd}\right),\quad -\underset{(\varphi,\dot{\varphi})\in O}{\inf}\mathcal{J}_{0,T}^N(\varphi)\leq \underset{\varepsilon\to 0}{\liminf}\,\varepsilon^2\log\mathbb{P}_{\vv{z}}\left(\grandeva{}\in O\right),\\
    &\bullet\text{for any closed subset }F\text{ of }C^0_{\vv{z}}\left([0;T],\mathbb{R}^{2Nd}\right),\quad \underset{\varepsilon\to 0}{\limsup}\,\varepsilon^2\log\mathbb{P}_{\vv{z}}\left(\grandeva{}\in F\right)
        \leq-\underset{(\varphi,\dot{\varphi})\in F}{\inf}\mathcal{J}_{0,T}^N(\varphi).
      \end{align*}
We may also recall the following class of rate functions, ensuring the existence of a minimizer when optimizing $\mathcal{J}_{0,T}^N(\varphi)$ on non-empty closed set of $\espacecont$: given $T>0$, we say that $\mathcal{J}_{0,T}^N(\varphi)$ is a \textit{good rate function} if, for all compact set $\mathcal{K}$ of $\reel^{2Nd}$ and every  $M>0$, the level set
    \begin{equation}
        \label{ens_niv}
        \mathcal{S}_T(\mathcal{K},M)=\{(\varphi,\dot{\varphi}):(\varphi(0),\dot{\varphi}(0))\in \mathcal{K},\, \varphi \in C^1([0;T],\mathbb{R}^{Nd}) \text{ and }\mathcal{J}_{0;T}^N(\varphi)\leq M\},
    \end{equation}
is a compact set of $\espacecont$. 

As previously discussed (see Proposition~\ref{ppdeva}),
given the constraint set 
$$
\mathcal{H}_T:=\{(\varphi_1,\varphi_2)\in \espacecont:\,\dot{\varphi}_1=\varphi_2 \text{ and }\varphi_1\in H_2\left([0;T]\right)  \},
$$ 
the family $\{\grandeva{}\}_{\varepsilon>0}$ satisfies a LDP with the good rate function:
\begin{equation}
    \label{ratefe}
    \mathcal{J}_{0,T}^N(\varphi)=\frac{1}{2}\sum_{i=1}^N\int_0^T\!\left\lVert \ddot{\varphi}^i(t)-b^{i,N}(\varphi(t),\dot{\varphi}(t))+\alpha\dot{\varphi}_1^i(t)+\gamma\left(\varphi^i(t)-x_i^\ast\right)\right\rVert^2\deriv t.
\end{equation}
We recall that this rate function should be understood as a measure of similarity to the deterministic counterpart of the system \eqref{Cucker-Smale-bruite}. Heuristically, the form can be obtained as follows. Consider on $[0;T]$ the unique solution of 
\begin{equation}
    \label{systeme_squelette}
    \dot{g}(t)=G^N(g(t))+\Gamma\dot{f}(t),\qquad 
    g(0)=\initia,
\end{equation}
for $\displaystyle f\!\in \!H^2\left([0;T];\reel^{2Nd}\right):=\{f\in\espacecontabs:f(0)=0 \text{ and } \int_0^T\lVert\dot{f}(t)\rVert^2 \deriv t<+\infty \}$ and 
\begin{align*}
    G^{N}\initias=
    \begin{pmatrix}
        \vv{v}=(v_i)_{i\leq N}\\
         \left(\Psi^{N}_i\initias\right)_i-\alpha \vv v-\gamma\left(\vv x-x^\ast\right)\\
    \end{pmatrix}, x^*:=\left((x^*_j)_{j \in \ens}\right)\quad\mathrm{and}\quad\Gamma=\begin{pmatrix}
        0&0\\
        0&\mathrm{I}_{Nd}
    \end{pmatrix}.
\end{align*}
Owing to the strong wellposedness of \eqref{Cucker-Smale-bruite},
 $\grandeva{}$ can be formulated as a transformation $\mathcal{G}\left(\varepsilon W\right)$ of a driving $\reel^{Nd}\times\reel^{Nd}$-Brownian motion.
Supposing that Brownian motions can be approximated by $f\in \!H^2\left([0;T];\reel^{2Nd}\right)$, \eqref{systeme_squelette} can be similarly viewed as $g=\mathcal{G}\left(\varepsilon f\right)$. Based on Schilder's theorem (see \cite[Th.~2.1,Th.~2.2, p.~58]{Freidlin_Wentzell}, \cite[Th.~5.2.3, p.~185]{Dembo_Zeitouni}) which provides the rate function of $\left(\varepsilon W\right)_{\varepsilon>0}$, and, provided that $\mathcal{G}$ is continuous, the contraction principle \cite[Th.~3.1,p.~64]{Freidlin_Wentzell} yields the following rate function for $\grandeva{}$: for all $(\varphi_1,\varphi_2)\in \espacecont$
\begin{equation*}
    \mathcal{J}_{0;T}(\varphi_1,\varphi_2)=\underset{\{f\in H^2:\mathcal{G}(f)=(\varphi_1,\varphi_2)\}}{\inf}\dfrac{1}{2}\int_0^T\lVert \dot{f}(t)\rVert^2 \deriv t,
\end{equation*}
with the convention $\inf\emptyset =+\infty$. Hence, explicitating the minimizing space eventually yields
$$\{f\in H^2:\mathcal{G}(f)=(\varphi_1,\varphi_2)\}=\mathcal H_T,$$ and so justifies the form \eqref{ratefe}. Rigorously, the above derivation involves a few technical issues which can eventually be carried out (we again refer to \cite{Azencott,Priouret,Baldi_Maurel} for details).

 We finally give two useful propositions for the rest of this section, stating some uniform control when $\vv{Z}^{N,\varepsilon}\!=\grandeva{}$ starts from a compact set or remains in a compact set. Their proofs can be found in \cite[p.~96-98]{Azencott}.
\begin{propositiona}
\label{devunifmin}
Given a compact set $\mathcal{K}$ of $\reel^{2Nd}$ and $T$, $a$, $\rho$ and $\eta$ positive, there exist $\varepsilon_0>0$ and $r>0$ such that if $\vv{z}=(\vv{x},\vv{v})\in \mathcal{K}$ and  $(\varphi,\dot{\varphi})\in C^0([0;T],\mathbb{R}^{2Nd})$ is such that  $(\varphi,\dot{\varphi})\left([0;T]\right)\subset \mathcal{K}$, $\norme{(\varphi(0),\dot{\varphi}(0))-(\vv{x},\vv{v})}<r$ and $\mathcal{J}_{0;T}(\varphi)\leq a$, then we have: for every $\varepsilon<\varepsilon_0$, 
\begin{equation*}
    -\mathcal{J}_{0;T}(\varphi)-\eta\leq \varepsilon^2\log\left(\mathbb{P}_{\vv{z}}\left(\mathcal{D}^{T}\left((\varphi,\dot{\varphi}),\vv{Z}^{N,\epsilon}
    \right)< \rho\right)\right).
\end{equation*}
\end{propositiona}
\begin{propositiona}
\label{devunifmax}
Consider a compact set $K$ of $\reel^{2Nd}$. Then, for every positive $T$,$\rho$, $\eta$ and $A$, one can find $\varepsilon_0>0$ such that if $(\vv{x},\vv{v})\in K$, for all $\varepsilon<\varepsilon_0$, we have, for $\mathcal{S}_{T}$ the level set in \eqref{ens_niv},
\begin{align}
    \varepsilon^2\log\left(\mathbb{P}_{\vv{z}}\left(\mathcal{D}^{T}\left[\vv{Z}^{N,\varepsilon},\mathcal{S}_{T}\big(\{\vv{z}\},A\big)\right]\geq \rho\right)\right)\leq -A+\eta.
\end{align}
\end{propositiona}

\subsection{Controllability and properties of the exit-cost }\leavevmode\par\vspace{0.5em}\label{sec4}
 In this section, we state  two key lemmas aiming to build a suitable exit-trajectory of the domain $\mathcal{O}_N^{\gamma}(R)$ to derive essential properties of the rate function $\mathcal{J}_{0,T}^N(\varphi)$. The first lemma is a tool to connect points of $\overline{\mathcal{O}_N^{\gamma}(r)}$ to the equilibrium  $\pointeq$ in finite time; the connecting trajectory having its rate function negligible when $r$ is sufficiently small. We also show that this interpolating function remains inside  $\mathcal{O}_N^{\gamma}(R)$. In the ``overdamped'' cases, the authors in~\cite{Freidlin_Wentzell} define a uniform rectilinear motion, connecting two distinct points $x,y$  with a finite time $T=\norme{x-y}$, given by $\varphi(t)=x+\dfrac{t}{\lVert x-y\lVert}(y-x)$. This idea does not work in our context because all interpolating functions should respect the kinematic constraint $\dot{x}(t)=v(t)$.\smallbreak
 To overcome this difficulty, we propose a motion with ``bang-bang'' acceleration on  time interval $[0;T]$, that is to say the motion has two main phases. On the first phase, particles undergo an acceleration equal to $a_1$ on interval $[0;\tau]$ and on the second phase, a deceleration $a_2$ on $[\tau;T]$, until they hit the desired target $\vv{X}_{\mathrm{eq}}^\ast$ with null velocities where $\tau$ is some time. 
 \begin{lemmea}
 \label{chemin_au_pteq}
There exist two functions $g^\gamma,h^{\lambda,\beta,\kappa}_{\gamma,\alpha}:[0,\infty)\longmapsto[0,\infty)$, positive, and for every $0<r\leq R$, and all $(\vv{x_0},\vv{v_0})\in \overline{\mathcal{O}_N^{\gamma}(r)}\setminus ~\{\pointeq\}$, there exists $0<T^*<T^\gamma(r)<\infty$ and a function $\psi\in H_2\left([0;T^*]\right)$, verifying 
\begin{equation}
    \begin{split}
   &\left(\psi(0),\dot{\psi}(0)\right)=(\vv{x_0},\vv{v_0}),\quad \left(\psi(T^\ast),\dot{\psi}(T^\ast)\right)=\pointeq\!, \\
   &\left(\psi,\dot{\psi}\right)\left(\big[0;T^\ast\big]\right)\subset \mathcal{O}_N^{\gamma}(g^\gamma(r))\quad\mathrm{and}\quad  \mathcal{J}_{0,T^\ast}^N(\psi)\leq h^{\lambda,\beta,\kappa}_{\gamma,\alpha}(r),
   \end{split}
\end{equation}
with $\displaystyle\underset{r\to 0}{\lim}\max\left(g^\gamma(r),h^{\lambda,\beta,\kappa}_{\gamma,\alpha}(r)\right)=0$.
 \end{lemmea}
 \begin{proof}
   On the finite time interval $[0;T]$, with $T>0$ to be chosen later,  define the path $\psi$ connecting $\vv{x}_0$ to $\vv{x}^*_{\mathrm{eq}}$ with $\psi(t)=\psi_1(t)=a_1\dfrac{t^2}{2}+\vv{v}_0t+\vv{x_0} \text{ if } t\in [0;\tau]$ and $\psi(t)=\psi_2(t)=a_2\dfrac{(T-t)^2}{2}+\vv{x}_{\mathrm{eq}}^\ast \text{ if } t\in [\tau;T]$ where $\tau:=\frac{T}{2}$. Choosing $a_1=\dfrac{4\Delta x}{T^2}-\dfrac{3\vv{v_0}}{T}$ and $a_2=\dfrac{\vv{v_0}}{T}-\dfrac{4\Delta x}{T^2}$,  where $\Delta x:=\vv{x}_{\mathrm{eq}}^\ast-\vv{x_0}$, ensures that $\psi\in C^1([0;T])$ with $\psi_1(\tau)=\psi_2(\tau)$, $\dot{\psi}_1(\tau)=\dot{\psi}_2(\tau)$ and $\dot{\psi}\in C^1\left([0;T]\setminus\{\tau\};\reel^{Nd}\right)$.\\
   Choosing then  $T^\ast:=\max\left(\norme{\vv{x_0}-\vv{x}_{\mathrm{eq}}^\ast}^{\frac{1}{3}},\norme{v_0}^{\frac{1}{2}}\right)$ and recall that, by assumption,  $\norme{\vv{v_0}}^2+\gamma\norme{\vv{x_0}-\vv{x}_{\mathrm{eq}}^\ast}^2\leq r$.
   We also set $T=T^ \gamma(r):=\left(\dfrac{r}{\gamma}\right)^{\frac{1}{6}}+r^\frac{1}{4}$, and $\tau^*$ accordingly, so that $T^\ast\leq T^\gamma(r)$  $\max(\norme{a_1},\norme{a_2})<4 T^\gamma(r)$ and, for every $t\in[\tau^\ast;T^\ast]$, 
      \begin{equation*}
          \norme{\psi(t)-\vv{x}_{\mathrm{eq}}^\ast}\leq \dfrac{(T^\ast)^2}{8}\norme{a_2}\leq \dfrac{[T^\gamma(r)]^2}{8}\norme{a_2}< \dfrac{[T^\gamma(r)]^3}{2}\quad\mathrm{and}\quad \norme{\dot{\psi}(t)}\leq \norme{a_2}\dfrac{T^\ast}{4}< [T^\gamma(r)]^2.
      \end{equation*}
      For all $t\in [0;\tau^\ast]$, we have $\norme{\dot{\psi}(t)}
      \leq \norme{a_1}\,T^\ast + \norme{\vv{v_0}}
      < 
      5[T^\gamma(r)]^2$ and
     \begin{align*}
\norme{\psi(t)-\vv{x}_{\mathrm{eq}}^\ast}
&\leq \frac{(T^\ast)^2}{8}\,\norme{a_1} + \norme{\vv{v_0}}\,\frac{T^\ast}{2} + \norme{\Delta x} 
<\dfrac{[T^\gamma(r)]^3}{2}+\dfrac{[T^\gamma(r)]^3}{2}+[T^\gamma(r)]^3=2[T^\gamma(r)]^3.
\end{align*}
This implies
\begin{align*}
    \norme{\dot{\psi}(t)}^2+\gamma\norme{\psi(t)-\vv{x}_{\mathrm{eq}}^\ast}^2<25[T^\gamma(r)]^4+4\gamma[T^\gamma(r)]^6,
\end{align*}
and so, defining $g^\gamma(r):=25[T^\gamma(r)]^4+4\gamma[T^\gamma(r)]^6$, we obtain $(\psi,\dot{\psi})\big(\big[0;T^\ast\big]\big)\subset \mathcal{O}_N^{\gamma}(g^\gamma(r))$.\\
It remains to show that the rate function of $\psi$ can be arbitrarily small with respect to $r$ in terms of a bound $h^{\lambda,\beta,\kappa}_{\gamma,\alpha}(r)$. The convexity inequality implies 
\begin{align*}
    \mathcal{J}^N_{0,T^\ast}(\psi)&\leq 2\Bigg(\Ss\inte{T^*}\norme{\ddot{\psi}_i(t)}^2\di t+\dfrac{\lambda^2}{\kappa^{2\beta} N}\Sd\inte{T^\ast}\norme{\dot{\psi}_i(t)-\dot{\psi}_j(t)}^2\di t\Bigg)\\
    &\quad+2\Bigg(\Ss\alpha^2\inte{T^\ast}\norme{\dot{\psi}_i(t)}^2\di t+\gamma^2\Ss\inte{T^\ast}\norme{{\psi_i}(t)-x_i^\ast}^2\di t\Bigg)\\
    &\leq 2\Bigg(\inte{T^\ast}\norme{\ddot{\psi}(t)}^2\di t+\left(\dfrac{4\lambda^2}{\kappa^{2\beta}}+\alpha^2\right)\inte{T^\ast}\norme{\dot{\psi}(t)}^2\di t+\gamma^2 \inte{T^\ast}\norme{{\psi}(t)-\vv{x}_{\mathrm{eq}}^\ast}^2\di t\Bigg).
\end{align*}
As we have the following upper-bounds
\begin{flalign*}
    &\inte{T^\ast}\norme{\ddot{\psi}(t)}^2\di t\leq \dfrac{T^\ast}{2}(\norme{a_1}^2+\norme{a_2}^2)\leq 16[T^\gamma(r)]^3,&\\
    &\inte{T^\ast}\norme{\dot{\psi}(t)}^2\di t\leq \dfrac{T^\ast}{2}\bigg(\underset{t\in[0;\tau^\ast]}{\max}\norme{\dot{\psi_1}(t)}^2+\underset{t\in[\tau^\ast;T^\ast]}{\max}\norme{\dot{\psi_2}(t)}^2\bigg)\leq 13[T^\gamma(r)]^5,&\\
    \mathrm{and}\quad&\inte{T^\ast}\norme{{\psi}(t)-\vv{x}_{\mathrm{eq}}^\ast}^2\di t\leq \dfrac{T^*}{2}\bigg(\underset{t\in[0;\tau^\ast]}{\max}\norme{\psi_1(t)-\vv{x}_{\mathrm{eq}}^\ast}^2+\underset{t\in[\tau^\ast;T^\ast]}{\max}\norme{\psi_2(t)-\vv{x}_{\mathrm{eq}}^\ast}^2\bigg)&\\
    &\quad\quad\quad\quad\quad\quad\quad\quad\quad\quad\leq\dfrac{T^\gamma(r)}{2}\bigg[4[T^\gamma(r)]^6+\dfrac{1}{4}[T^\gamma(r)]^6\bigg]=\dfrac{17}{8}[T^\gamma(r)]^7,
\end{flalign*}
we get $\mathcal{J}_{0,T^\ast}^N(\psi)\leq h^{\lambda,\beta,\kappa}_{\gamma,\alpha}(r)$ for
\begin{equation}
    \label{majo_quasipo}
h^{\lambda,\beta,\kappa}_{\gamma,\alpha}=    2\Bigg(16[T^\gamma(r)]^2+13\bigg(\dfrac{4\lambda^2}{\kappa^{2\beta}}+\alpha^2\bigg)[T^\gamma(r)]^5+\dfrac{17\gamma^2}{8}[T^\gamma(r)]^7\Bigg).
\end{equation}

 \end{proof}
 We now give a second lemma stating that one can construct a particles system with starting points in $\partial\mathcal{O}_N^{\gamma}(R)$ leaving outside $\overline{\mathcal{O}_N^{\gamma}(R)}$. Such particles system is given by $\vv \Psi^N=(\Psi_i)_{i\leq N}$ solution of the ODEs 
 \begin{equation}
    \label{dynamique_inverse}
     \ddot{\Psi}_i(t)=-b^{i,N}\big((\Psi_j(t),\dot{\Psi}_j(t))_{j\in \ens}\big)+\alpha\dot{\Psi}_i(t)-\gamma(\Psi_i(t)-x_i^\ast),\qquad \big(\Psi(0),\dot{\Psi}(0)\big)=(x_0^i,v_0^i)\in\ensinvab{R},
 \end{equation}
 where $b^{i,N}$ is defined in \eqref{fonction_communis}. This system admits a unique solution (since the $b^{i,N}$ are locally Lipschitz and with at most linear growth). Let us suppose that  $(\vv{x_0},\vv{v_0})\in \ensinvab{R}$ and introduce the energy function (following \eqref{ener_u}):
 \begin{equation*}
     \mathcal{E}(t):=\mathcal{U}^{\gamma}_{\mathrm{tot}}\big([\Psi_i(t),\dot{\Psi}_i(t)]_{1\le i\leq N}\big)=\Ss\norme{\dot{\Psi}_i(t)}^2 +\gamma\Ss\norme{\Psi_i(t)-x_i^\ast}^2,\qquad t\ge 0.
 \end{equation*}
 \begin{lemmea}
 \label{trajectoire_sortie} $\mathcal{E}$ is increasing on $\reel^+$, and
 there exists $\delta_c>0$ such that, for all $t\geq 0$, $\mathcal{E}(t)\geq \dfrac{1}{3}\mathcal{E}(0)e^{\frac{2}{3}\delta_c t}$. 
 \end{lemmea}
 \begin{proof} For every $t\geq 0$, we have
   \begin{align}
        \label{derive_ener}
       \dfrac{\di}{\di t}\mathcal{E}(t)
       &=\dfrac{2}{N}\Sd \phi\left(\norme{\Psi_i(t)-\Psi_j(t)}\right)\left\langle\dot{\Psi}_i(t);\dot{\Psi}_i(t)-\dot{\Psi}_j(t)\right\rangle+2\alpha\Ss\langle\dot{\Psi}_i(t),\dot{\Psi}_i(t)\rangle\notag\\
       &\quad-2\gamma\Ss\langle\dot{\Psi}_i(t);\Psi_i(t)-x_i^\ast\rangle+2\gamma\Ss\left\langle\dot{\Psi}_i(t);\Psi_i(t)-x_i^\ast\right\rangle\notag\\
       &=\dfrac{1}{N}\Ss\phi(\norme{\Psi_i(t)-\Psi_j(t)})\norme{\dot{\Psi}_i(t)-\dot{\Psi}_j(t)}^2+2\alpha\Ss\norme{\dot{\Psi}_i(t)}^2\geq 0.
   \end{align}
   The above yields
   \begin{equation}
   	\{t\geq 0:\,\dfrac{\di}{\di t}\mathcal{E}(t)=0\}=\{t\geq 0:\dot{\Psi}_i(t)=0,\forall i\in\ens\}.
   \end{equation}
   Thus if $\mathcal{E}$ was not increasing over a non-trivial interval $[T_1;T_2]$, in view of the left part of \eqref{dynamique_inverse} necessarily $\Psi_i(t)=x^*_i$ within $[T_1;T_2]$ which is impossible since $(\vv{x}^*_{\mathrm{eq}},0_{\reel^{dN}})$ lies outside $\ensinvab{R}$. Hence $\mathcal{E}$ is necessarily increasing.
   \\
   To obtain the lower-bound $\mathcal{E}(t)\ge \dfrac{1}{3}\mathcal{E}(0)e^{\frac{2}{3}\delta_c t}$, as in Section \ref{sec:RefinedLyapunov}, we introduce the energy 
   \begin{align}
   	\mathcal{E}^\sigma(t)=\mathcal{U}^{\sigma,\gamma}([(\Psi_i(t),\dot{\Psi}_i(t))_{1\le i\le N}])=\mathcal{E}(t)+\sigma\Ss\left\langle\Psi_i(t)-x_i^\ast;\dot{\Psi}_i(t)\right\rangle, 
   \end{align}
   with $\sigma\in \reel$. Following \textit{in extenso} the proof of Corollary \ref{dissipatif_part}, for all $-\sqrt{\gamma}\leq\sigma\leq \sqrt{\gamma}$ and all $t\geq 0$, we have $\dfrac{1}{2}\mathcal{E}(t)\leq\mathcal{E}^\sigma(t)\leq \dfrac{3}{2}\mathcal{E}(t)$. Studying further the growth of $\mathcal{E}^\sigma$ when $\sigma$ is in a suitable interval $(-\sigma_c;0)$ will show the exponential explosion to $+\infty$ of $\mathcal{E}^\sigma$ and eventually the claim.
   
For $\sigma<0$, and, $t\geq 0$, we have
   \begin{align*}
       \dfrac{\di}{\di t}\mathcal{E}^\sigma(t)
       &=\dfrac{1}{N}\Ss\phi(\norme{\Psi_i(t)-\Psi_j(t)})\norme{\dot{\Psi}_i(t)-\dot{\Psi}_j(t)}^2\quad \\
       &\quad+\dfrac{\sigma}{N}\Sd\phi(\norme{\Psi_i(t)-\Psi_j(t)})\left\langle \Psi_i(t)-x_i^\ast;\dot{\Psi}_i(t)-\dot{\Psi}_j(t)\right\rangle\quad\left(=:\sigma\mathcal{I}_1(t)\right)\\
       &\quad+(2\alpha+\sigma)\Ss\norme{\dot{\Psi}_i(t)}^2-\sigma\gamma\Ss \left\langle\Psi_i(t)-x_i^\ast;\Psi_i(t)-x_i^\ast\right\rangle+\sigma\alpha \Ss\left\langle\Psi_i(t)-x_i^\ast;\dot{\Psi}_i(t)\right\rangle(=:\sigma\mathcal{I}_2(t)).
 \end{align*}
 
For $\mathcal{I}_1(t)$, using Cauchy-Schwarz and Young inequalities and recalling $\phi\le \phi^+=\dfrac{\lambda}{\kappa^\beta}$, we have 
 \begin{align}\label{ener_1}
     \mathcal{I}_1(t) &\leq\dfrac{\phi^+}{N}\Sd\norme{\Psi_i(t)-x_i^\ast}\norme{\dot{\Psi}_i(t)-\dot{\Psi}_j(t)}\notag\\
     &\leq \dfrac{\gamma}{2N}\Sd\norme{\Psi_i(t)-x_i^\ast}^2+\dfrac{(\phi^+)^2}{2N\gamma}\Sd\norme{\dot{\Psi}_i(t)-\dot{\Psi}_j(t)}^2\notag\\
     &=\dfrac{\gamma}{2}\Ss\norme{\Psi_i(t)-x_i^\ast}^2+\dfrac{(\phi^+)^2}{\gamma}\Ss\norme{\dot{\Psi}_i(t)}^2-\dfrac{(\phi^+)^2N}{\gamma}\norme{\frac1{N}\sum_{j=1}^N\Psi_j(t)}^2.
 \end{align}
 For $\mathcal{I}_2^\sigma(t)$, let $\eta>0$ and applying again thanks to Cauchy-Schwarz and Young inequalities, we get
 \begin{equation}\label{ener_2}
     \mathcal{I}_2(t)\leq \alpha\Ss\sqrt{\eta}\norme{\Psi_i(t)-x_i^\ast}\dfrac{1}{\sqrt{\eta}}\norme{\dot{\Psi}_i(t)}
     \leq \dfrac{\eta}{2}\Ss \norme{\Psi_i(t)-x_i^\ast}^2+\dfrac{\alpha^2}{2\eta}\Ss\norme{\dot{\Psi}_i(t)}^2.
 \end{equation}
 Combining \eqref{ener_1} and \eqref{ener_2} (taking into account that the factor $\sigma<0$ inverts the inequality), it follows that 
 \begin{align}
    \label{ineq_diff_explo}
     \dfrac{\di}{\di t}\mathcal{E}^\sigma(t)&\geq \left(2\alpha+\sigma\left(1+\dfrac{(\phi^+)^2}{\gamma}+\dfrac{\alpha^2}{2\eta}\right)\right)\Ss\norme{\dot{\Psi}_i(t)}^2-\sigma\Big(\dfrac{\gamma-\eta}{2}\Big)\Ss \norme{\Psi_i(t)-x_i^\ast}^2.
 \end{align}
  To get positive factors in the right-hand side of the inequality
 we choose $\eta<\gamma$, and next $-\sigma_c<\sigma<0$ for
 \begin{equation*}
     \sigma_c=\min\left(\sqrt{\gamma};2\alpha\Big(1+\frac{(\phi^+)^2}{\gamma}+\frac{\alpha^2}{2\eta}\Big)^{-1}\right).
 \end{equation*}
 Then, we define
 \begin{equation*}
     \delta_c=\min\left(2\alpha+\sigma\left(1+\dfrac{(\phi^+)^2}{\gamma}+\dfrac{\alpha^2}{2\eta}\right);-\sigma\Big(\dfrac{\gamma-\eta}{2}\Big)\right),
 \end{equation*}
 and the differential inequality \eqref{ineq_diff_explo} implies $\dfrac{\di}{\di t}\mathcal{E}^\sigma(t)\geq \delta_c \mathcal{E}(t)\geq \dfrac{2}{3}\delta_c\mathcal{E}^\sigma(t)$.
Thus, we have $\mathcal{E}(t)\geq \dfrac{1}{3}\mathcal{E}(0)e^{\frac{2}{3}\delta_c t}$.
 \end{proof}
Now, we establish some properties of the exit-cost \eqref{csortie}; essentially its positiveness which is crucial for the validity of the main Theorems \ref{loi_arrhenius} and \ref{l_kram}, and justifies its exponential scale. 
\begin{lemmea}
    \label{posi_cou}
    Let us consider $h^{\lambda,\beta,\kappa}_{\gamma,\alpha}$ given at \eqref{majo_quasipo} and the exit-cost $\mathcal{V}_N(R,\alpha,\gamma)$. Then, we have
    \begin{equation*}
        \dfrac{\alpha R}{2}<\mathcal{V}_N(R,\alpha,\gamma)\leq h^{\lambda,\beta,\kappa}_{\gamma,\alpha}\left(R\right).
    \end{equation*}
    In particular, $\displaystyle\lim_{\alpha\to+\infty}\mathcal{V}_N(R,\alpha,\gamma)=\lim_{R\to+\infty}\mathcal{V}_N(R,\alpha,\gamma)=+\infty$ and $\displaystyle\lim_{R\to 0}\mathcal{V}_N(R,\alpha,\gamma)=0$.
\end{lemmea}
\begin{proof}
Let $T>0$ and $\varphi=(\varphi_i)_i $ be a function in $H_2\left([0;T]\right)$, such that $\left(\varphi(0),\dot{\varphi}(0)\right)=\pointeq$ and $\left(\varphi(T),\dot{\varphi}(T)\right)\in \ensinvab{R}$. Recall that $ \phi_{ij}(t):=\phi\left(\norme{\varphi_i(t)-\varphi_j(t)}\right)$ where $\phi(x):=\lambda\left(\kappa+x^2\right)^{-\beta}$. Write $\displaystyle\mathcal{J}_{0;T}^N(\varphi)=\dfrac{1}{2}\sum_{i=1}^N\inte{T}A_i(t)\di t$ where
$$A_i(t):=\Big\|\ddot{\varphi}_i(t)+\alpha\dot{\varphi}_i(t)+\gamma\left(\varphi_i(t)-x_i^{\ast}\right)+\dfrac{1}{N}\sum_{j=1}^N\phi_{ij}(t)\left(\dot{\varphi}_i(t)-\dot{\varphi}_j(t)\right)\Big\|^2.$$
Observing that
\begin{align*}
    A_i(t)
    &=\norme{\ddot{\varphi}_i(t)+\gamma\left(\varphi_i(t)-x_i^\ast\right)}^2+\alpha^2\norme{\dot{\varphi}_i(t)}^2+\norme{\dfrac{1}{N}\sum_{j=1}^N\phi_{ij}(t)\left(\dot{\varphi}_i(t)-\dot{\varphi}_j(t)\right)}^2\\
    &\quad+2\alpha\scalp[\dot{\varphi}_i(t)]{\ddot{\varphi}_i(t)}+2\alpha\gamma\scalp[\dot{\varphi}_i(t)]{\varphi_i(t)-x_i^\ast}\\
    &\quad+2\scalp[\dfrac{1}{N}\sum_{j=1}^N\phi_{ij}(t)\left(\dot{\varphi}_i(t)-\dot{\varphi}_j(t)\right)]{\ddot{\varphi}_i(t)+\gamma\left(\varphi_i(t)-x_i^\ast\right)}\\
    &\quad+\dfrac{2\alpha}{N}\sum_{j=1}^N\phi_{ij}(t)\scalp[\dot{\varphi}_i(t)]{\dot{\varphi}_i(t)-\dot{\varphi}_j(t)},
\end{align*}
using the identity $\displaystyle 2\sum_{j=1}^N\phi_{ij}(t)\scalp[\dot{\varphi}_i(t)]{\dot{\varphi}_i(t)-\dot{\varphi}_j(t)}=\sum_{j=1}^N\phi_{ij}(t)\norme{\dot{\varphi}_i(t)-\dot{\varphi}_j(t)}^2$ and the inequality $2\scalp[a]{b}\geq -\norme{a}^2-\norme{b}^2$, we deduce 
\begin{align}
    \label{ineq_bon_fct}
    \dfrac{1}{2}\Ss\inte{T}A_i(t)\di t&\geq \alpha\Ss\inte{T}\scalp[\dot{\varphi}_i(t)]{\ddot{\varphi}_i(t)} \di t+\alpha\gamma\Ss\inte{T}\scalp[\dot{\varphi}_i(t)]{\varphi_i(t)-x_i^\ast}\di t+\alpha^2\Ss\inte{T}\norme{\dot{\varphi}_i(t)}^2\di t\notag\\
    &=\dfrac{1}{2}\left(\Ss\norme{\dot{\varphi}_i(T)}^2+\gamma\Ss \norme{{\varphi}_i(T)-x_i^*}^2\right)+\alpha^2\Ss\inte{T}\norme{\dot{\varphi}_i(t)}^2\di t\notag\\
    &=\dfrac{\alpha R}{2}+\alpha^2\Ss\inte{T}\norme{\dot{\varphi}_i(t)}^2\di t .
\end{align}
The last term of the right-hand side is necessarily positive, otherwise $\dot{\varphi}(t)=0_{}$ for all $t\in[0;T]$ and so  $\varphi(t)=\vv{x}_{\mathrm{eq}}^\ast$ meanwhile $(\vv{x}_{\mathrm{eq}}^\ast,0_{\reel^{Nd}}) \notin\ensinvab{R}$. For the upper-bound, it suffices to take the path connecting $\ensinvab{R}$ to $\pointeq$ given by Lemma \ref{chemin_au_pteq} and reverse it in time; $\mathcal{V}_N(R,\alpha,\gamma)\leq h^{\lambda,\beta,\kappa}_{\gamma,\alpha}\left(R\right)$ immediately follows.
\end{proof}
\subsection{Technical lemmas}\label{ssec:Appendice-lemmestech}\leavevmode\par\vspace{0.5em}\label{sec6}
In this section, we give a series of important lemmas used in the proof of Theorem \ref{loi_arrhenius} and Theorem \ref{l_kram}. Although the overall strategy is inspired by the arguments of \cite{Freidlin_Wentzell}, several technical building blocks do not directly carry over. In particular, as the diffusion of \eqref{Cucker-Smale-bruite} is degenerate, several auxiliary lemmas therein must be revisited and extended. \\
The first one states that the system of particles \eqref{Cucker-Smale-bruite} can not remain \textit{ad infinitum} in what we call the ``$r$-annulus'' $\overline{\mathcal{O}_N^{\gamma}(R)}\setminus\mathcal{O}_N^{\gamma}(r)$, when starting from the latter. Thus, there is a dichotomy: either the particles hit $\ensinvad{r}$, or  they  shall exit the domain, and hit $\ensinvab{R}$ at some time. We define the following stopping-time:
 \begin{equation}
     \Lambda_{\varepsilon}^{\alpha,\gamma,N}(R,r)=\inf\{t\geq 0,\, \left(\vv{X}^{N,\varepsilon}_t,\vv{V}^{N,\varepsilon}_t\right)\in \partial\mathcal{O}_N^{\gamma}(R)\cup \overline{\mathcal{O}}_N^{\alpha,\gamma}(r)\}.
 \end{equation}
\begin{lemmea}[Time spent in the annulus]
 \label{temps_dans_annneau}

 For every $0<r<R$, we have the following statements:
 \begin{enumerate}
     \item There exist $\upsilon(r)>0$ and a deterministic finite time $T(r)>0$ such that 
     \begin{equation*}
         \mathcal{J}_{0;T}^N(\varphi)\geq \upsilon(r)\left[T-T(r)\right],
     \end{equation*}
      for all $T\geq T(r)$ and all $(\varphi,\dot{\varphi})\in C^0\left([0;T],\mathbb{R}^{2Nd}\right) $ such that $(\varphi,\dot{\varphi})([0,T])\subset \overline{\mathcal{O}_N^{\gamma}(R)}\setminus\mathcal{O}_N^{\gamma}(r)$.
     \item There exists a real number $\rho(r)>0$, a time $T(r)>0$ and $\varepsilon_0>0$ such that for all  $\varepsilon<\varepsilon_0$, all $\vv{z_0}=(\vv{x_0},\vv{v_0})\in \overline{\mathcal{O}_N^{\gamma}(R)}\setminus\mathcal{O}_N^{\gamma}(r)$, and every $T\geq T(r)$, we have 
     \begin{equation}
         \probabi\left(\Lambda_{\varepsilon}^{\alpha,\gamma,N}(R,r)>T\right)\leq \exp{\bigg(-\dfrac{\rho(r)}{\varepsilon^2}\big(T-T(r)\big)\bigg)}.
     \end{equation}
 \end{enumerate}
 \end{lemmea}
 
 \begin{proof}\textbf{Step 1.} Let us consider $\vv{z_0}=(\vv{x_0},\vv{v_0})$ in $\overline{\mathcal{O}_N^{\alpha,\gamma}(R)}$. Then, one can find a  time (uniform with respect to initial conditions) such that the noiseless particles system  $\big(\vv{x}^{N},\vv{v}^{N}\big)$ hits the set $\mathcal{O}_N^{\gamma}\left(\dfrac{r}{2}\right)$. Thanks to Theorem~\ref{theoreme_Ineq_Diff_1}, there exists  $\delta>0$ such that 
 \begin{equation}
     \mathcal{U}^{\gamma,N}_{\mathrm{tot}}(t):=\sum_{i=1}^N\lVert v^{i,N}(t)\rVert^2+\gamma\sum_{i=1}^N\lVert x^{i,N}(t)-x_i^\ast\rVert^2\leq 3Re^{-\frac{2}{3}\delta t},
 \end{equation}
 Hence, defining the time $T(r)=\dfrac{3}{\delta} \log(6R/r)$, for all $t\geq T(r)$, we get $
 \mathcal{U}^{\gamma,N}_{\mathrm{tot}}(t)<r/2$. We next introduce 
 \begin{equation*}
     \mathcal{C}(r)=\{(\varphi,\dot{\varphi}):\varphi \in C^1\left([0;T(r)],\mathbb{R}^{Nd}\right),\,\left(\varphi,\dot{\varphi}\right)([0;T(r)])\subset \overline{\mathcal{O}_N^{\gamma}(R)}\setminus\mathcal{O}_N^{\gamma}(r)\}.
 \end{equation*}
 (which is non-empty and closed) and put   $c(r)=\underset{(\varphi,\dot{\varphi})\in \mathcal{C}(r)}{\inf}\mathcal{J}_{0;T(r)}(\varphi)$. The latter is positive since $\mathcal{J}_{0;T}^N$ is a good rate function (see Proposition~\ref{ppdeva}), and attains its minimum on $\mathcal{C}(r)$ at $(\varphi^*,\dot{\varphi}^*)$; $c(r)$ so can not be equal to $0$, otherwise $(\vv{x}^N,\vv{z}^N)$ should be at the same time in $\overline{\mathcal{O}_N^{\gamma}(R)}\setminus\mathcal{O}_N^{\gamma}(r)$ and in $\mathcal{O}_N^{\gamma}(r/2)$ which is impossible. \\
 \textbf{Step 2.} Let $T\geq T(r)$ and let $\varphi$ be a path such that $(\varphi,\dot{\varphi})([0;T])\subset \overline{\mathcal{O}_N^{\gamma}(R)}\setminus\mathcal{O}_N^{\gamma}(r)$. Put $\Delta_1(r):=\left\lfloor \frac{T}{T(r)}\right\rfloor$ and
 \begin{align*}
     \varphi^k(t)=\varphi\left(t+kT(r)\right), \forall (t,k)\in [0;T(r)]\times [\![0;\Delta_1(r)-1]\!].
 \end{align*}
 For any $k\in [\![0;\Delta_1(r)-1]\!]$, the function $(\varphi^k,\dot{\varphi}^k)$ evolves inside $\overline{\mathcal{O}_N^{\gamma}(R)}\setminus\mathcal{O}_N^{\gamma}(r)$ and we deduce 
 \begin{align*}
     \mathcal{J}_{0;T(r)}^N(\varphi)\geq \sum_{k=0}^{\Delta_1(r)-1}\mathcal{J}_{0;T(r)}^N(\varphi^k)\geq c(r)\Delta_1(r)\geq \dfrac{c(r)}{T(r)}\left(T-T(r)\right),
 \end{align*}
 yielding the first point of the lemma.\\
 \textbf{Step 3.} Let us now take $\eta>0$ and consider the enlargement  $\overline{\mathcal{O}_N^{\gamma}(R+\eta)}$ of $\overline{\mathcal{O}_N^{\gamma}(R)}.$
 We consider the (compact) level set $\mathcal{S}_{T(r)}\big(\{(\vv{x_0},\vv{v_0})\},c(r)\big)$  given at \eqref{ens_niv} and define 
 \begin{align*}
    \mathcal{K}_1&=\{(\varphi,\dot{\varphi})\in C^0\left([0;T(r)],\mathbb{R}^{2Nd}\right):(\varphi,\dot{\varphi})([0;T(r)])\subset\overline{\mathcal{O}_N^{\gamma}(R+\eta)}\}.
 \end{align*}
 Hence, we get the decomposition $\mathcal{S}_{T(r)}\big(\{(\vv{x_0},\vv{v_0})\},c(r)\big) =:\mathcal{V}_1\cup\mathcal{V}_2$ with 
 \begin{equation}
     \mathcal{V}_1=\mathcal{S}_{T(r)}\Big(\{(\vv{x_0},\vv{v_0})\},c(r)\Big)\bigcap \mathcal{K}_1,\quad \mathcal{V}_2:=\mathcal{S}_{T(r)}\Big(\{(\vv{x_0},\vv{v_0})\},c(r)\Big)\bigcap \prescript{}{}{\mathcal{K}_1^c}.
 \end{equation}
 We then set $ \sigma_r:=\delta_r\land \alpha_R$, where the scalars $ \delta_r$ and $\alpha_R$ are defined as follow
 \begin{align*}
     &\delta_r:=\inf\{\norme{\vv{z}-\vv{z}'}:(\vv{z},\vv{z'})\in \overline{\mathcal{O}_N^{\gamma}\left(\dfrac{r}{2}\right)}\times \left(\overline{\mathcal{O}_N^{\gamma}(R)}\setminus\mathcal{O}_N^{\gamma}(r)\right)\},\\
     &\alpha_R:= \inf\{\norme{\vv{z}-\vv{z}'}:(\vv{z},\vv{z'})\in \overline{\mathcal{O}_N^{\gamma}\left(R\right)}\times \prescript{}{}{\left(\overline{\mathcal{O}_N^{\gamma}(R+\eta)}\right)^c}\}.
 \end{align*}
The sets $\overline{\mathcal{O}_N^{\gamma}(R)}\setminus\mathcal{O}_N^{\gamma}(r)$ and $\overline{\mathcal{O}_N^{\gamma}(r/2)}$ are disjoint compact sets of $\reel^{2Nd}$ so that $\delta_r>0$ ; likewise $\alpha_R>0$ by compactness of $\overline{\mathcal{O}_N^{\gamma}\left(R\right)}$. From the first point of the lemma, there exists $c(r)>0$ and $T(r)>0$ such that $\mathcal{J}_{0;T(r)}(\phi)>c(r)$, for all $\phi \in C^1([0;T(r)],\mathbb{R}^{Nd}) $ and $(\phi,\dot{\phi})([0;T(r)])\subset \overline{\mathcal{O}_N^{\gamma}(R+\eta)}\setminus\mathcal{O}_N^{\gamma}(r/2)$. \\
Let us now consider a path $\psi \in C^1([0,T(r)],\mathbb{R}^{Nd})$ such that $(\psi(0),\dot{\psi}(0))=(\vv{x_0},\vv{v_0})\in  \overline{\mathcal{O}_N^{\gamma}(R)}\setminus\mathcal{O}_N^{\gamma}(r)$ and $(\psi,\dot{\psi})([0;T(r)]\subset \overline{\mathcal{O}_N^{\gamma}(R)}\setminus\mathcal{O}_N^{\gamma}(r)$. If we state that $ (\varphi,\dot{\varphi})\in \mathcal{V}_1$ then this function hits $\mathcal{O}_N^{\gamma}\left(r/2\right)$ before $T(r)$, otherwise it should remain in $\overline{\mathcal{O}_N^{\gamma}(R+\eta)}\setminus\mathcal{O}_N^{\gamma}\left(\dfrac{r}{2}\right)$ leading to the contradiction $c(r)<c(r)$. Thus, we have $\mathcal{D}^{T(r)}\big((\varphi,\dot{\varphi}),(\psi,\dot{\psi})\big)\geq \delta_r\geq \sigma_r$. If we suppose that $(\varphi,\dot{\varphi})\in \mathcal{V}_2 $, since there exists a time $t^\ast$ before $T(r)$ where $(\varphi,\dot{\varphi})$ is outside $\overline{\mathcal{O}_N^{\gamma}(R+\eta)} $, we have $\mathcal{D}^{T(r)}\big((\varphi,\dot{\varphi}),(\psi,\dot{\psi})\big)\geq \alpha_R\geq \sigma_r$. Consequently $$
\mathcal{D}^{T(r)}\big((\psi,\dot{\psi}),\mathcal{S}_{T(r)}\big(\{(\vv{x_0},\vv{v_0})\},c(r)\big)\big)\geq \sigma_r,
$$
 for all $C^1$-path $\psi$ such that  $(\psi,\dot{\psi})$ starts from $\initia$ and remains in $ \overline{\mathcal{O}_N^{\gamma}(R)}\setminus\mathcal{O}_N^{\gamma}(r)$ over $[0;T(r)]$. Applying Proposition \ref{devunifmax} with 
   $\vv{z_0}=\left(\vv{x_0},\vv{v_0}\right)\in\overline{\mathcal{O}_N^{\gamma}(R)}\setminus\mathcal{O}_N^{\gamma}(r)$, there is  $\varepsilon_0>0$ (depending on $r$ and not on $\vv{z_0}$) such that for every $\varepsilon\leq \varepsilon_0$, and  
 \begin{align}
    \label{estim_prob}
     \probabi\left(\Lambda_{\varepsilon}^{\alpha,\gamma,N}(R,r)>T(r)\right)&\leq\probabi\left(\left(\vv{X}^{N,\varepsilon},\vv{V}^{N,\varepsilon}\right)\in F_r\right)\leq \e^{-\frac{c(r)}{2\varepsilon^2}},
 \end{align}
 for
 \begin{equation*}
 	F_r:=\{(\varphi,\dot{\varphi})\in C^0\left([0,T(r)],\mathbb{R}^{2Nd}\right):\mathcal{D}^{T(r)}\big((\varphi,\dot{\varphi}),\mathcal{S}_{T(r)}\big(\{(\vv{x_0},\vv{v_0})\},c(r)\big)\big)\geq \sigma_r\}.
 \end{equation*}
 \textbf{Step 4.} For the extension of \eqref{estim_prob} to times larger than $T(r)$, we use the Markov property of $\left(\vv{X}^{N,\varepsilon},\vv{V}^{N,\varepsilon}\right)$. Recall the time-shift operator $\theta_t$ in \eqref{fonctionnelle} and define the functional on the path space $C([0,\infty);\reel^{2Nd})$:
 \begin{align}
    \label{op_transla}
     &\Lambda((x_s,y_s)_{s\geq 0})=\inf\{t\geq 0,(x_t,y_t)\in \partial\mathcal{O}_N^{\gamma}(R)\cup \overline{\mathcal{O}_N^{\gamma}(r)}\},
 \end{align}
     so that $\Lambda_{\varepsilon}^{\alpha,\gamma,N}(R,r)=\Lambda(\vv{X}^{N,\varepsilon},\vv{V}^{N,\varepsilon})$. For simplicity, we denote $\Lambda_{\varepsilon}^{\alpha,\gamma,N}\left(R,r\right)=\Lambda_{\varepsilon}^{\alpha,\gamma,N}$ and again write $\vv{Z}^{N,\varepsilon}=(\vv{X}^{N,\varepsilon},\vv{V}^{N,\varepsilon})$. On the event $\{\Lambda_{\varepsilon}^{\alpha,\gamma,N}>nT(r)\}$, we have the identity $\Lambda_{\varepsilon}^{\alpha,\gamma,N}=nT(r)+\Lambda\big(\theta_{nT(r)}(\vv{Z}^{N,\varepsilon})\big)$. Thus
     \begin{align*}
         \probabi\left(\Lambda_{\varepsilon}^{\alpha,\gamma,N}>(n+1)T(r)\right)&=\mathbb{E}_{\vv{z_0}}\left(\mathbbm{1}_{\{\Lambda\left(\theta_{nT(r)}\left(\vv{X}^{N,\varepsilon},\vv{V}^{N,\varepsilon}\right)\right)> T(r)\}\cap \{\Lambda_{\varepsilon}^{\alpha,\gamma,N}>nT(r)\}}\right)\\
         &=\mathbb{E}_{\vv{z_0}}\left(\mathbbm{1}_{\{\Lambda_{\varepsilon}^{\alpha,\gamma,N}>nT(r)\}}Y\right),
     \end{align*}
     where, by the Markov property,
     \begin{align*}
         Y&=\mathbb{E}_{\vv{Z}^{N,\varepsilon}_{nT(r)}}\left(\mathbbm{1}_{\{\Lambda(\vv{Z}^{N,\varepsilon})> T(r)\} }\right)=\mathbb{E}_{\vv{Z}^{N,\varepsilon}_{nT(r)}}\left(\mathbbm{1}_{\{\Lambda_{\varepsilon}^{\alpha,\gamma,N}> T(r)\}}\right).
     \end{align*}
     Hence, since on  $\{\Lambda_{\varepsilon}^{\alpha,\gamma,N}>nT(r)\}$, $\vv{Z}^{N,\varepsilon}_{nT(r)}$ is in the annulus $\overline{\mathcal{O}_N^{\gamma}(R)}\setminus\mathcal{O}_N^{\gamma}(r)$, we deduce
     \begin{align*}
         \probabi\left(\Lambda_{\varepsilon}^{\alpha,\gamma,N}>(n+1)T(r)\right)&\leq \underset{(\vv{x},\vv{y})\in\overline{\mathcal{O}_N^{\gamma}(R)}\setminus\mathcal{O}_N^{\gamma}(r)}{\max}\!\!\!\!\!\!\!\!\!\!\!\!\!\!\!\!\mathbb{P}_{(\vv{x},\vv{y})}\left(\Lambda_{\varepsilon}^{\alpha,\gamma,N}>T(r)\right)\probabi\left(\Lambda_{\varepsilon}^{\alpha,\gamma,N}> nT(r)\right)\\
         &\leq \exp{\big(\!\!\!-c(r)/(2\varepsilon^2)\big)}\probabi\left(\Lambda_{\varepsilon}^{\alpha,\gamma,N}> nT(r)\right)
     \end{align*}
     using \eqref{estim_prob} for the last estimate. By induction it follows that, for $\varepsilon\leq \varepsilon_0$, $$\probabi\left(\Lambda_{\varepsilon}^{\alpha,\gamma,N}>n T(r)\right)\leq \exp{\Bigg(\!\!\!-\dfrac{nc(r)}{2\varepsilon^2}\Bigg)}.$$ 
     Hence, for every $T\geq T(r)$,  the proof is completed with
     \begin{align*}
         \probabi\left(\Lambda_{\varepsilon}^{\alpha,\gamma,N}(R,r)>T\right)&\leq \probabi\left(\Lambda_{\varepsilon}^{\alpha,\gamma,N}(R,r)>T(r)\Bigg\lfloor \dfrac{T}{T(r)}\Bigg\rfloor\right)\\
         &\leq \exp{\Bigg(\!\!\!-\Bigg\lfloor \dfrac{T}{T(r)}\Bigg\rfloor\dfrac{c(r)}{2\varepsilon^2}\!\Bigg)}\leq \exp{\Bigg(\!\!\!-\dfrac{c(r)\big[T-T(r)\big]}{2T(r)\varepsilon^2}\!\Bigg)}.
     \end{align*}
     \end{proof}
 We state a lemma showing the construction  of a given deterministic time such that the probability of exiting $\mathcal{O}_N^{\gamma}(R)$ before this given time has a lower-bound proportional to $\exp{\left(-\frac{C}{\varepsilon^2}\right)}$ with $C$ arbitrarily close to the exit-cost.
 \begin{lemmea}[Exiting before a finite time]
\label{sortiefini}
 For every $\delta>0$, there exist a deterministic time $T(\delta)>0$ and $\varepsilon_0>0$ such that for all $\vv{z_0}=(\vv{x_0},\vv{v_0})\in \mathcal{O}_N^{\gamma}(R)$ and for all $\varepsilon<\varepsilon_0$, we have
  \begin{equation}
     \probabi\left(\tau^{\alpha,\gamma,N}_{\varepsilon}(R)<T(\delta)\right)\geq \dfrac{1}{2} \exp{\Bigg(\!\!\!-\dfrac{1}{\varepsilon^2}\bigg(\mathcal{V}_N(R,\alpha,\gamma)+\dfrac{\delta}{2}\bigg)\!\Bigg)}.
 \end{equation}
\end{lemmea}
\begin{proof}

  Let us take $\delta>0$ and  $0<r<R$ which shall be calibrated later on.\\
  \textbf{Step 1.} We set $\widehat{T}_{1,r}=\dfrac{3}{\delta} \log\big(18R/r\big)$
such that $\grandevad{} $ remains in $\mathcal{O}_N^{\gamma}(r/6)$ after this time. For now,  we construct a deterministic time $\widehat{T}_{2,r}$ such that starting from $\overline{\mathcal{O}_N^{\gamma}(r)}$ a particles system exits $\mathcal{O}_N^{\gamma}(R)$ before $\widehat{T}_{2,r}$, by reaching the stable point $\pointeq$ by means of a smooth path $\psi$, having $\mathcal{J}_{0;T(r)}^N(\psi,\dot{\psi})$ close to the exit-cost $\mathcal{V}_N(R,\alpha,\gamma)$.\\By the characterization of the exit-cost, there exist $T^{\delta}>0$ and a path $\psi^{(1)}\in H^2([0;T^\delta])$ verifying: 
\begin{align*}
    &(\psi^{(1)}(0),\dot{\psi}^{(1)}(0))=\pointeq ,\quad (\psi^{(1)}(T^\delta),\dot{\psi}^{(1)}(T^\delta))=(x_{\mathrm{out}}^\delta,v_{\mathrm{out}}^\delta)\in \partial\mathcal{O}_N^{\gamma}(R)\\
    &\left(\psi^{(1)},\dot{\psi}^{(1)}\right)\left(]0;T^\delta[\right)\subset\mathcal{O}_N^{\gamma}(R)\quad \mathrm{and}\quad\mathcal{J}_{0,T^\delta}^N\left(\psi^{(1)}\right)\leq \mathcal{V}_N(R,\alpha,\gamma)+\dfrac{\delta}{12}.
\end{align*}
Let us consider $(\vv{x_0},\vv{v_0})\in\overline{\mathcal{O}_N^{\gamma}(r)}$. Thanks to Lemma~\ref{chemin_au_pteq}, we connect $\initia$ to the equilibrium point $\pointeq$ via a function $\psi^{(2)}$ such that there exists a finite time $T^*=T^*(\vv{z_0})\leq T(r)$ for which we have 
\begin{equation}
    \begin{split}
   &(\psi^{(2)}(0),\dot{\psi}^{(2)}(0))=(\vv{x_0},\vv{v_0}),\quad (\psi^{(2)}(T^\ast),\dot{\psi}^{(2)}(T^\ast))=\pointeq \\
   &(\psi^{(2)},\dot{\psi}^{(2)})\big(\big[0;T^\ast\big]\big)\subset \mathcal{O}_N^{\gamma}(g(r))\quad \mathrm{and} \quad \mathcal{J}_{0,T^\ast}^N\left(\psi^{(2)}\right)\leq h(r),
   \end{split}
\end{equation}
with $g(r)$ and $h(r)$ going to $0$ when $r\to 0$. So, we can choose  $r$ small enough so that 
\begin{equation*}
   (\psi^{(2)},\dot{\psi}^{(2)})\big(\big[0;T^\ast\big]\big)\subset \mathcal{O}_N^{\gamma}(R)\quad\mathrm{and}\quad  \mathcal{J}_{0,T^\ast}\left(\psi^{(2)}\right)\leq \dfrac{\delta}{12}.
\end{equation*}
We next connect $(x_{\mathrm{out}}^\delta,v_{\mathrm{ou}t}^\delta)$ to a point out of  $\overline{\mathcal{O}_N^{\gamma}(R)}$. To do so, put $\psi^{(3)}:=\vv{\Psi}^N\big|_{[0;\widehat{T}_c]}$ the restriction to $[0;\widehat{T}_c]$ of $\vv{\psi}^N$ the particles system \eqref{dynamique_inverse} starting from $(x_{\mathrm{out}}^\delta,v_{\mathrm{out}}^\delta)$ at time $0$. The time $\widehat{T}_c$ is chosen so that the rate function of $\mathcal{J}_{0,\widehat{T}(\mu)}^N$ at $\psi^{(3)}$ is small for some positive $\widehat{T}(\mu)$. To achieve this, let us take $\mu,\eta>0$, consider $\partial\mathcal{O}_N^{\gamma}(R+\eta)$ and define
\begin{equation}
    T_\eta^\delta=\inf\{t>0:\big(\vv{\psi}^N(0),\vv{\dot{\psi}}^N(0)\big)=(x_{\mathrm{out}}^\delta,v_{\mathrm{out}}^\delta)\text{ and }\big(\vv{\Psi}^N(t),\vv{\dot{\Psi}}^N(t)\big)\in\partial\mathcal{O}_N^{\gamma}(R+\eta)\}.
\end{equation}
By virtue of Lemma~\ref{trajectoire_sortie}, this hitting-time is positive and finite with  $T_\eta^\delta\leq T(R,\eta)=\dfrac{3}{2\delta_c}\log\left(\dfrac{3(R+\eta)}{R}\right)$.
From this, we define $\widehat{T}(\mu)=\mu T_\eta^\delta$ with $0<\mu<1$ such that $\mathcal{J}_{0,\widehat{T}(\mu)}^N\left(\psi^{(3)}\right)\leq \dfrac{\delta}{12}$. This is obtained by observing that  
\begin{align*}
    \mathcal{J}_{0,\widehat{T}(\mu)}^N\left(\psi^{(3)}\right)&\leq \Ss\inte{\widehat{T}(\mu)}\norme{\dfrac{2}{N}\sum_{j=1}^N\phi\left(\norme{\psi^{(3)}_i(t)-\psi^{(3)}_j(t)}\right)\left(\dot{\psi}^{(3)}_i(t)-\dot{\psi}^{(3)}_j(t)\right)+2\alpha \dot{\psi}^{(3)}_i(t)}^2\di t\notag\\
    &\leq \inte{\widehat{T}(\mu)}\left(\left(\frac{4\lambda^2}{N\kappa^{2\beta}}+\frac{4\lambda\alpha}{N\kappa^\beta}\right)\Sd \norme{\dot{\psi}^{(3)}_i(t)-\dot{\psi}^{(3)}_j(t)}^2+4\alpha^2\Ss \norme{\dot{\psi}^{(3)}_i(t)}^2\right)\di t\notag\\
    &\leq \inte{\widehat{T}(\mu)}\left[\frac{8\lambda^2}{\kappa^{2\beta}}+\frac{8\lambda\alpha}{\kappa^\beta}+4\alpha^2\right]\Ss \norme{\dot{\psi}^{(3)}_i(t)}^2\di t.
\end{align*}
Since for all $t\leq \widehat{T}(\mu)$, we have $ \Ss \norme{\dot{\psi}^{(3)}_i(t)}^2<R+\eta$, it follows that
\begin{equation*}
    \mathcal{J}_{0,\widehat{T}(\mu)}^N\left(\psi^{(3)}\right)\leq \widehat{T}(\mu)(R+\eta)\left[\frac{8\lambda^2}{\kappa^{2\beta}}+\frac{8\lambda\alpha}{\kappa^\beta}+4\alpha^2\right].
\end{equation*}
The bound $\mathcal{J}_{0,\widehat{T}(\mu)}^N\left(\psi^{(3)}\right)\leq \dfrac{\delta}{12}$ is so granted by choosing $\mu$ and $\widehat{T}_c$ with $$\mu<\mu_0:=\min\left(1,\delta\left(12T_\eta^\delta\left[\frac{8\lambda^2}{\kappa^{2\beta}}+\frac{8\lambda\alpha}{\kappa^\beta}+4\alpha^2\right](R+\eta)\right)^{-1}\right)\quad\mathrm{and}\quad \widehat{T}_c:=\widehat{T}\left(\dfrac{\mu_0}{2}\right).
$$
 From $\psi^{(i)}$, we construct a path connecting $\overline{\mathcal{O}_N^{\gamma}(r)}$ to the set  outside  $\overline{\mathcal{O}_N^{\alpha,\gamma}(R)}$, and reaching the equilibrium, with
 \begin{equation}\label{lebonchemin}
     \psi(t)=\begin{cases}
        &\!\!\!\!\!\!\psi^{(2)}(t),\text{ if } t\in \big[0;T^\ast\big],\\
       & \!\!\!\!\!\!\psi^{(1)}(t-T^\ast),\text{ if } t\in \big[T^\ast;T^\ast+T^\delta\big],\\
     &\!\!\!\!\!\!\psi^{(3)}\big(t-(T^\ast+T^\delta)\big), \text{ if } t\in \big[T^\ast+T^\delta;T^\ast+T^\delta+\widehat{T}_c\big],
     \end{cases}
 \end{equation}
 $\psi$ being defined up to $T((\vv{x_0},\vv{v_0}))=T_0:=T^\ast+T^\delta+\widehat{T}_c\leq \widehat{T}_{2,r}:=T(r)+T^\delta+\widehat{T}_c$. Along $\psi$,  $\mathcal{J}_{0,T_0}^N$ satisfies
 \begin{equation}
    \label{Minoration_taux}
     \mathcal{J}_{0,T_0}^N(\psi)\leq \mathcal{V}_N(R,\alpha,\gamma)+\dfrac{\delta}{4}.
 \end{equation}
 \textbf{Step 2.} We suppose $\vv{z_0}=(\vv{x_0},\vv{v_0})\in \overline{\mathcal{O}_N^{\gamma}(r)}$. Then, we know that $\psi$  as in \eqref{lebonchemin} verifies $
\left(\psi(T_0),\dot{\psi}\left(T_0\right)\right) \notin \overline{\mathcal{O}_N^{\gamma}(R)}
$. Thus, there exists $h>0$ such that $B\left(\left(\psi(T_0),\dot{\psi}(T_0)\right),h\right)\cap \overline{\mathcal{O}_N^{\gamma}(R)}=\emptyset$. Since, we have $\left(\psi,\dot{\psi}\right)([0;T_0])\subset\overline{\mathcal{O}_N^{\gamma}(R+\eta)}$, Proposition \ref{devunifmax} implies that there is a $\varepsilon_0>0$, independent of $\initia$, and such that, 
for every $\varepsilon<\varepsilon_0$,
 \begin{equation*}
     \probabi\bigg(\mathcal{D}^{T_0}\left(\left(\vv{X}^{N,\varepsilon},\vv{V}^{N,\varepsilon}\right),(\psi,\dot{\psi})\right)<h\bigg) \geq \exp{\Bigg(\!\!\!-\dfrac{1}{\varepsilon^2}\bigg(\mathcal{J}_{0,T_0}^N(\psi)+\dfrac{\delta}{4}\bigg)\!\Bigg)}.
 \end{equation*}
 From \eqref{Minoration_taux}, we further deduce the following inequality:
 \begin{equation}
     \probabi\left(\mathcal{D}^{T_0}\left(\left(\vv{X}^{N,\varepsilon},\vv{V}^{N,\varepsilon}\right),(\psi,\dot{\psi})\right)<h\right) \geq \exp{\Bigg(\!\!\!-\dfrac{1}{\varepsilon^2}\bigg(\mathcal{V}_N(R,\alpha,\gamma)+\dfrac{\delta}{2}\bigg)\!\Bigg)}.
 \end{equation}
 For $\mathcal{C}_\varepsilon:=\{\left(\vv{X}^{N,\varepsilon}_{T_0},\vv{V}^{N,\varepsilon}_{T_0}\right) \in B\left(\left(\psi(T_0),\dot{\psi}(T_0)\right),h\right) \}$, we have the inclusions
 \begin{align*}
     \{\mathcal{D}^{T_0}\left(\left(\vv{X}^{N,\varepsilon},\vv{V}^{N,\varepsilon}\right),(\psi,\dot{\psi})\right)<h\}&\subset \mathcal{C}_\varepsilon\subset\{\left(\vv{X}^{N,\varepsilon}_{T_0},\vv{V}^{N,\varepsilon}_{T_0}\right) \notin  \overline{\mathcal{O}_N^{\gamma}(R)}\}\subset \{\tau^{\alpha,\gamma,N}_{\varepsilon}(R)<T_0\},
 \end{align*}
 the last right set lying itself in $\{\tau^{\alpha,\gamma,N}_{\varepsilon}(R)<\widehat{T}_{2,r}\}$.
 Therefore, for all $\varepsilon<\varepsilon_0$  and $\vv{z_0}\in \overline{\mathcal{O}_N^{\gamma}(r)}$, we obtain
 \begin{equation}
    \label{mino_t2}
     \probabi\Big(\tau^{\alpha,\gamma,N}_{\varepsilon}(R)<\widehat{T}_{2,r}\Big)\geq \exp{\Bigg(\!\!\!-\dfrac{1}{\varepsilon^2}\bigg(\mathcal{V}_N(R,\alpha,\gamma)+\dfrac{\delta}{2}\bigg)\!\Bigg)},
 \end{equation}
 and so
  \begin{equation}
    \label{mino_t2a}
     \probabi\Big(\tau^{\alpha,\gamma,N}_{\varepsilon}(R)<\widehat{T}_{1,r}+\widehat{T}_{2,r}\Big)\geq \dfrac{1}{2}\exp{\Bigg(\!\!\!-\dfrac{1}{\varepsilon^2}\bigg(\mathcal{V}_N(R,\alpha,\gamma)+\dfrac{\delta}{2}\bigg)\!\Bigg)}.
 \end{equation}
 \textbf{Step 3.} We extend \eqref{mino_t2a} to all points $\vv{z_0}$ of $\mathcal{O}_N^{\gamma}(R)$. Let us take $\vv{z_0}\in \mathcal{O}_N^{\gamma}(R)\setminus\overline{\mathcal{O}_N^{\gamma}(r)}$ and define 
 \begin{equation}
     \sigma_\varepsilon:=\inf\{t\geq 0, \left(\vv{X}^{N,\varepsilon}_t,\vv{V}^{N,\varepsilon}_t\right)\in\overline{\mathcal{O}_N^{\gamma}(r)}\}=\tau_{\ensinvad{r}}\left(\vv{X}^{N,\varepsilon},\vv{V}^{N,\varepsilon}\right),
 \end{equation}
  the first hitting-time of $\overline{\mathcal{O}_N^{\gamma}(r)}$
(recalling, from \eqref{fonctionnelle}, that $\tau_A\left(\left(x_t,y_t\right)_{t\geq 0}\right)=\inf\{t\geq 0, (x_t,y_t)\in A\}$). For $\theta_t$ the time-shift operator $\theta$ (see \eqref{fonctionnelle}), we have the inclusion $ \mathcal{A}^1_\varepsilon\cap\mathcal{A}^2_\varepsilon\cap \mathcal{A}^3_\varepsilon\subset\{\tau^{\alpha,\gamma,N}_{\varepsilon}(R)<\widehat{T}_{1,r}+\widehat{T}_{2,r}\}
$, for
\begin{align*}
  &\mathcal{A}^1_\varepsilon:=  \{\sigma_\varepsilon<\tau^{\alpha,\gamma,N}_{\varepsilon}(R)\},\quad \mathcal{A}^2_\varepsilon:=\{\sigma_\varepsilon<\widehat{T}_{1,r}\},\quad\mathcal{A}^3_\varepsilon:=\{\tau_{\ensinvab{R}}\left(\theta_{\sigma_\varepsilon}\left(\vv{X}^{N,\varepsilon},\vv{V}^{N,\varepsilon}\right)\right)<\widehat{T}_{2,r}\}.
\end{align*}
 Hence $\esper{\left(\mathbbm{1}_{\mathcal{A}^1_\varepsilon\cap\mathcal{A}^2_\varepsilon\cap \mathcal{A}^3_\varepsilon}\right)}\leq \probabi\big(\tau^{\alpha,\gamma,N}_{\varepsilon}(R)<\widehat{T}_{1,r}+\widehat{T}_{2,r}\big)$.
 Applying the Markov property and \eqref{mino_t2}, with $\varepsilon<\varepsilon_0$, 
 \begin{align}
    \label{A_estim}
     \esper{\left(\mathbbm{1}_{\mathcal{A}^1_\varepsilon\cap\mathcal{A}^2_\varepsilon\cap \mathcal{A}^3_\varepsilon}\right)}
     &=\esper{\left(\mathbbm{1}_{\{\sigma_\varepsilon<\tau^{\alpha,\gamma,N}_{\varepsilon}(R)\}\cap\{\sigma_\varepsilon<\widehat{T}_{1,r}\}\cap\{\tau_{\ensinvab{R}}\left(\theta_{\sigma_\varepsilon}\left(\vv{X}^{N,\varepsilon},\vv{V}^{N,\varepsilon}\right)\right)<\widehat{T}_{2,r}\}}\right)}\notag\\
     &=\esper{\left(\mathbbm{1}_{\{\sigma_\varepsilon<\widehat{T}_{1,r}\}\cap\{0<\tau_{\ensinvab{R}}\left(\theta_{\sigma_\varepsilon}\left(\vv{X}^{N,\varepsilon},\vv{V}^{N,\varepsilon}\right)\right)<\widehat{T}_{2,r}\}}\right)}\notag\\
      &=\esper{\left(\mathbbm{1}_{\{\sigma_\varepsilon<\widehat{T}_{1,r}\}}\esper{\left(\mathbbm{1}_{\{0<\tau_{\ensinvab{R}}\left(\theta_{\sigma_\varepsilon}\left(\vv{X}^{N,\varepsilon},\vv{V}^{N,\varepsilon}\right)\right)<\widehat{T}_{2,r}\}}\middle\vert\mathcal{F}_{\sigma_\varepsilon}\right)}\right)}\notag\\
      &\geq\probabi\left(\sigma_\varepsilon<\widehat{T}_{1,r}\right)\underset{(\vv{x},\vv{y})\in \ensinvad{r}}{\inf}\mathbb{P}_{(\vv{x},\vv{y})}\left(\tau^{\alpha,\gamma,N}_{\varepsilon}(R)<\widehat{T}_{2,r}\right)\notag\\
      &\geq\probabi\left(\sigma_\varepsilon<\widehat{T}_{1,r}\right)\exp{\Bigg(\!\!\!-\dfrac{1}{\varepsilon^2}\bigg(\mathcal{V}_N(R,\alpha,\gamma)+\dfrac{\delta}{2}\bigg)\!\Bigg)}.
 \end{align}
 Since we have the following inequality 
 \begin{align*}
        \gamma\norme{\vv{X}^{N,\varepsilon}_{\widehat{T}_{1,r}}-\vv{x}_{\mathrm{eq}}^\ast}^2+\norme{\vv{V}^{N,\varepsilon}_{\widehat{T}_{1,r}}}^2&\leq 2\gamma\norme{\vv{X}^{N,\varepsilon}_{\widehat{T}_{1,r}}-\vv{x}^{N}(\widehat{T}_{1,r})}^2+2\gamma\norme{\vv{x}^{N}(\widehat{T}_{1,r})-\vv{x}_{\mathrm{eq}}^\ast}^2\\
        &\quad+2\norme{\vv{V}^{N,\varepsilon}_{\widehat{T}_{1,r}}-\vv{v}^{N}(\widehat{T}_{1,r})}^2+2\norme{\vv{v}^{N}(\widehat{T}_{1,r})}^2,
 \end{align*}
and, by definition of $\widehat{T}_{1,r}$ at Step 1, it follows that
 \begin{align*}
     \probabi{\left(\sigma_\varepsilon<\widehat{T}_{1,r}\right)}\geq\probabi{\big(A_\varepsilon\cap B_\varepsilon \big)},
 \end{align*}
 where the sets $A_\varepsilon$ and $B_\varepsilon$ are defined by
 \begin{equation*}
     A_\varepsilon=\{\underset{0\leq t\leq \widehat{T}_{1,r}}{\max}\norme{\vv{X}^{N,\varepsilon}_t-\vv{x}^{N}(t)}<\sqrt{\dfrac{r}{12}}\}\quad\text{and}\quad B_\varepsilon=\{\underset{0\leq t\leq \widehat{T}_{1,r}}{\max}\norme{\vv{V}^{N,\varepsilon}_t-\vv{v}^{N}(t)}<\sqrt{\dfrac{r}{12}}\}.
 \end{equation*}
 Using successively $\probabi{\big(A_\varepsilon\cap B_\varepsilon \big)}\geq \probabi{\big(A_\varepsilon\big)}+\probabi{\big(B_\varepsilon\big)}-1 $,  Theorem~\ref{convL1}, since $\vv{z_0}$ is in the compact $\ensinvad{R}$, one can find $\varepsilon_1>0$ (independently of $\vv{z_0}$), such that for any $\varepsilon<\varepsilon_1$, $\probabi{\left(A_\varepsilon\cap B_\varepsilon \right)}\geq \dfrac{1}{2}$. Consequently, setting $\varepsilon_2:=\min(\varepsilon_0,\varepsilon_1)$, the inequality \eqref{A_estim} implies: for every  $\varepsilon<\varepsilon_2$, 
 \begin{equation}
     \probabi\left(\tau^{\alpha,\gamma,N}_{\varepsilon}(R)<\widehat{T}_{1,r}+\widehat{T}_{2,r}\right)\geq \dfrac{1}{2} \exp{\Bigg(-\dfrac{1}{\varepsilon^2}\bigg(\mathcal{V}_N(R,\alpha,\gamma)+\dfrac{\delta}{2}\bigg)\Bigg)}.
 \end{equation}
Defining $T(\delta):=\widehat{T}_{1,r}+\widehat{T}_{2,r}$ ends the proof.
\end{proof}
\noindent
Recalling $\{\zeta_n^\varepsilon\}_n$ in \eqref{temps_osci}
we show that, as $\varepsilon\rightarrow0$, the probability of $\{\zeta_1^\varepsilon=\tau^{\alpha,\gamma,N}_{\varepsilon}(R)\}$ decreases exponentially fast.
\begin{lemmea}
\label{ev_impossi}
 For any $0<h<\mathcal{V}_N(R,\alpha,\gamma)$, there exist $r,\varepsilon_0>0$ such that, for all $\vv{z_0}\in  \mathcal{O}_N^{\gamma}(R)$ and $\varepsilon<\varepsilon_0$, 
 \begin{equation}
     \probabi{\left(\zeta_1^\varepsilon=\tau^{\alpha,\gamma,N}_{\varepsilon}(R)\right)}\leq \exp{\left(-\frac{1}{\varepsilon^2}\left(\mathcal{V}_N(R,\alpha,\gamma)-h\right)\right)}.
 \end{equation}
 \end{lemmea}
 \begin{remarquea}
 	\label{extent_impo}
 	Replicating \textit{in extenso} the proof of Lemma \ref{ev_impossi} below -- under the same assumptions -- we can deduce that, for all $n\in\mathbb{N}^\ast$, and all $\varepsilon<\varepsilon_0$, with $\varepsilon_0>0$ not depending on the initial conditions, we also have
 	\begin{equation}
 		\probabi{\left(\zeta_n=\tau^{\alpha,\gamma,N}_{\varepsilon}(R)\right)}\leq \exp{\left(\!\!-\frac{1}{\varepsilon^2}\left(\mathcal{V}_N(R,\alpha,\gamma)-h\right)\!\right)}.
 	\end{equation}
 \end{remarquea}
 \begin{proof}
    Let us take $T>0$ arbitrary for the moment, and set $\mathcal{A}(r)=\partial\mathcal{O}_N^{\gamma}(R)\cup \partial\mathcal{O}_N^{\gamma}\left(\dfrac{r}{4}\right)$ and $\mathcal{B}(r)=\partial\mathcal{O}_N^{\gamma}(R)$. Considering the functional $\tau_A$ defined at \eqref{fonctionnelle}, the strong Markov property with stopping-times $\chi_0^\varepsilon$ (see \eqref{temps_osci}) yields
 \begin{align}
    \label{estim_first_exit}
    \begin{aligned}
    \probabi{(\zeta_1^\varepsilon=\tau^{\alpha,\gamma,N}_{\varepsilon}(R))}&=\esper{\left(\mathbbm{1}_{\{\tau_{\mathcal{A}(r)}\left(\theta_{\chi_0}(\vv{Z}^{N,\varepsilon})\right)=\tau_{\mathcal{B}(r)}\left(\theta_{\chi_0}(\vv{Z}^{N,\varepsilon}_t)\right)\}}\right)}\\
     &\leq \underset{\vv{z}\in \partial\mathcal{O}_N^{\gamma}(r)}{\max}\mathbb P_{\vv{z}}{\big(\tau_{\mathcal{A}(r)}(\vv{Z}^{N,\varepsilon})=\tau^{\alpha,\gamma,N}_{\varepsilon}(R)\big)}\\
     &\leq \underset{\vv{z}\in \partial\mathcal{O}_N^{\gamma}(r)}{\max} \left[P_1(T,\vv{z})+P_2(T,\vv{z})\right],
     \end{aligned}
 \end{align}
 where the quantities $P_1$ and $P_2$ are defined as 
 \begin{equation*}
    P_1(T,\vv{z})=\mathbb P_{\vv{z}}{\left(\tau_{\mathcal{A}(r)}(\vv{Z}^{N,\varepsilon})=\tau^{\alpha,\gamma,N}_{\varepsilon}(R)\geq T\right)}\quad\mathrm{and}\quad
     P_2(T,\vv{z})=\mathbb P_{\vv{z}}{\left(\tau_{\mathcal{A}(r)}(\vv{Z}^{N,\varepsilon})=\tau^{\alpha,\gamma,N}_{\varepsilon}(R)<T\right)}.
 \end{equation*}
 \textbf{For $P_1(T,\vv{z})$.} Thanks to Lemma~\ref{temps_dans_annneau}, there exist a constant $\rho(r)>0$, a finite deterministic time $\Delta(r)$ and $ \varepsilon_0>0$ such that for every $\vv{z}\in \partial\mathcal{O}_N^{\gamma}(r)$, $P_1(T,\vv{z})\leq \exp{\Big(\!\!-\frac{\rho(r)}{\varepsilon^2}\left(T-\Delta(r)\right)\!\Big)}$, for all $T\geq \Delta(r) $ and all $\varepsilon<\varepsilon_0$. Thus, for any $T\geq \Delta^\ast_r:= \Delta(r)+\dfrac{\mathcal{V}_N(R,\alpha,\gamma)}{\rho(r)}$,
 for every $\initias\in \partial\mathcal{O}_N^{\gamma}(r)$ and for all $\varepsilon<\varepsilon_0$, we have
 \begin{equation}
    \label{P_2_T}
    P_1(T,\vv{z})= \mathbb P_{\vv{z}}{\big(\tau_{\mathcal{A}(r)}(\vv{Z}^{N,\varepsilon})=\tau^{\alpha,\gamma,N}_{\varepsilon}(R)\geq T \big)}\leq \exp{\Big(\!\!-\frac{1}{\varepsilon^2}\mathcal{V}_N(R,\alpha,\gamma)\!\Big)}.
 \end{equation}
 \textbf{For $P_2(T,\vv{z})$.} Let us take $0<h<\mathcal{V}_N(R,\alpha,\gamma)$ (which is immediate, for $R,\alpha,\gamma>0$ considering Lemma~\ref{posi_cou}), and consider the compact level set $\mathcal{S}_{T}(h,r):=\mathcal{S}_{T}\left(\partial\mathcal{O}_N^{\gamma}(r),\mathcal{V}_N(R,\alpha,\gamma)-\dfrac{h}{2}\right)$ defined at~\eqref{ens_niv} and the closed~set 
 \begin{equation*}
     \mathcal{W}_r^T=\{\big(\phi,\dot{\phi}\big):\,\phi \in C^1([0;T],\reel^{Nd}),\big(\phi(0),\dot{\phi}(0)\big)\in \partial\mathcal{O}_N^{\gamma}(r)\,\text{and}\, \big(\phi(t),\dot{\phi}(t)\big)\in \partial\mathcal{O}_N^{\gamma}(R)\,\text{ for some $t\in [0,T]$}\},
 \end{equation*}
  consisting all paths hitting $\partial\mathcal{O}_N^{\gamma}(R)$ before $T$.  
 We state that there exists a deterministic time $\widehat{T}_{r}$ such that for all $T\geq \widehat{T}_{r}$, $\mathcal{W}_r^T$ is not empty. Indeed, thanks to Lemma \ref{trajectoire_sortie}, we can build a path starting from $\vv{z_0}\in\ensinvab{r}$ and reaching $\ensinvab{R}$ at a time $T(\vv{z_0})$. This time is bounded by $\widehat{T}_{r}:=\dfrac{3}{2\delta_c}\log\left(\dfrac{3R}{r}\right)$, still according to Lemma~\ref{trajectoire_sortie}.\\
 Next defining $\widehat{\Delta}_r:=\Delta^\ast_r\vee \widehat{T}_{r}$, we have $\mathcal{S}_{\widehat{\Delta}_r}(h,r)\cap~\mathcal{W}_r^{\widehat{\Delta}_r}=\emptyset $ for $r$ sufficiently small. Indeed, if the opposite was true, there should exist a path $\left(\varphi,\dot{\varphi}\right)\in\mathcal{W}_r^T$ starting from $\partial\mathcal{O}_N^{\gamma}(r) $ at time $0$, reaching $\partial\mathcal{O}_N^{\gamma}(R)$ at time $t^\ast$ before $\widehat{\Delta}_r$ and such that $\mathcal{J}_{0;T}(\varphi)$ function is below $\mathcal{V}_N(R,\alpha,\gamma)-\dfrac{h}{2} $. Yet, Lemma~\ref{chemin_au_pteq} ensure the existence of $ (\psi,\dot{\psi})$ connecting $ \left(\varphi(0),\dot{\varphi}(0)\right)$ to the equilibrium point $\pointeq$ in a given time $t^\ast\leq t(r)$, staying at $\pointeq$ until  $t(r)$ and remaining in $\mathcal{O}_N^{\gamma}(R) $ after. For $r$ small, on $[0;t(r)]$, $\mathcal{J}_{0;t(r)}(\psi)$ is below $\dfrac{h}{4}$. Combining $\varphi$ and $\psi$, the function 
  \begin{equation}
     \phi(t)=\begin{cases}
        &\!\!\!\!\!\!\psi(t(r)-t),\text{ if } t\in \left[0;t(r)\right],\\
       &\!\!\!\!\!\! \varphi(t-t(r)),\text{ if } t\in \left[t(r);t^\ast+t(r)\right],\\
     \end{cases}
 \end{equation}
 immediately yields to the contradiction: $\mathcal{V}_N(R,\alpha,\gamma)\leq \mathcal{J}_{0;t^\ast+t(r)}(\phi) \leq \mathcal{V}_N(R,\alpha,\gamma)-\dfrac{h}{4}$. So $\left(\varphi,\dot{\varphi}\right)$ can not exist.\\
As $\mathcal{S}_{\widehat{\Delta}_r}(h,r)\cap~\mathcal{W}_r^{\widehat{\Delta}_r}=\emptyset $, for $r$ small, we can define the positive scalar $\delta^\prime(r):=\mathcal{D}^{\widehat{\Delta}_r}\left(\mathcal{S}_{\widehat{\Delta}_r}(h,r),\mathcal{W}_r^{\widehat{\Delta}_r}\right)$.
Using Proposition~\ref{devunifmax}, there exists $\varepsilon_1>0$ such that for every $\varepsilon<\varepsilon_1$ and for all $\vv{z}=\initias\in \partial\mathcal{O}_N^{\gamma}(r)$: 
 \begin{align}
    \label{P_1_T}
    \begin{aligned}
     P_2(\widehat{\Delta}_r,\vv{z})&\leq \mathbb P_{\vv{z}}\bigg( \mathcal{D}^{\widehat{\Delta}_r}\Big(\vv{Z}^{N,\varepsilon},\mathcal{S}_{\widehat{\Delta}_r}\big(\{\vv{z}\},\mathcal{V}_N(R,\alpha,\gamma)-\dfrac{h}{2}\big)\Big)\geq \delta^\prime(r)\bigg)\\
     &\leq \exp{\left(-\frac{1}{\varepsilon^2}\bigg(\mathcal{V}_N(R,\alpha,\gamma)-\dfrac{h}{2}-\dfrac{h}{6}\bigg)\right)}= \exp{\left(-\frac{1}{\varepsilon^2}\bigg(\mathcal{V}_N(R,\alpha,\gamma)-\dfrac{2h}{3}\bigg)\right)}.
     \end{aligned}
 \end{align}
Consequently, taking $\varepsilon_2<\min(\varepsilon_0,\varepsilon_1)$, from \eqref{estim_first_exit}, \eqref{P_2_T} with $T:=\widehat{\Delta}_r$ and \eqref{P_1_T}, we conclude, for  $\varepsilon<\varepsilon_2$:
 \begin{align}\label{estimation_lieu_sortie}
      \probabi{\left(\zeta_1^\varepsilon=\tau^{\alpha,\gamma,N}_{\varepsilon}(R)\right)}&\leq \exp{\left(\!\!-\frac{1}{\varepsilon^2}\mathcal{V}_N(R,\alpha,\gamma)\!\right)}+\exp{\left(\!\!-\frac{1}{\varepsilon^2}\bigg(\mathcal{V}_N(R,\alpha,\gamma)-\dfrac{2h}{3}\bigg)\!\right)} \notag\\
&\leq \exp{\left(\!\!-\frac{1}{\varepsilon^2}\bigg(\mathcal{V}_N(R,\alpha,\gamma)-h\bigg)\!\right)}.
 \end{align}
 \end{proof}
 The following lemma states that with high probability, for small $\varepsilon$ and  starting from $\ensinvab{r}$, $r<R$, the first hitting-time of  $ \grandeva{}$ to the set $\partial\mathcal{O}_N^{\gamma}(R)\cup \partial\mathcal{O}_N^{\gamma}\left(r/4\right)$ starting from $\ensinvab{r}$, admits a deterministic lower-bound. 
 \begin{lemmea}
 \label{temp_min}
 There exist $r\in (0,R)$, $\varepsilon_0>0$ and $T(r)>0$ such that, for every $\vv{z_0}=\initia\in \ensinvab{r}$ and all $\varepsilon<\varepsilon_0$, we have, for $\tau_{\mathcal A(r)}=\inf\{t\ge 0\,:\,\vv{Z}^{N,\varepsilon}_t\in \mathcal A(r)\}$, with $\mathcal{A}(r):=\partial\mathcal{O}_N^{\gamma}(R)\cup \partial\mathcal{O}_N^{\gamma}\left(r/4\right)$,
\begin{equation}
     \probabi{\left(\tau_{\mathcal{A}(r)}\geq T(r)\right)}\geq\eta.
 \end{equation}
 \end{lemmea}
 \begin{proof}
      According to Lemma \ref{Lyapou_cible_naive}, moreover, the energy $\mathcal{U}^{\gamma,N}_{\mathrm{tot}}(t)=\mathcal{U}^{\gamma}_{\mathrm{tot}}(\left[x_i(t),v_i(t)\right]_{1\le i\leq N})=\sum_{i=1}^N\lVert v_i(t)\rVert^2+\gamma\sum_{i=1}^N\left\lVert x_i(t) -x_i^\ast\right\rVert^2$ associated to the (zero-noise) system \eqref{C-S-conf_1} is nonincreasing. Hence, for every  $\vv{z_0}\in \partial\mathcal{O}_N^{\gamma}(r)$, there exists a  deterministic minimal time $T(\vv{z_0})$ such that $\grandevad{T(\vv{z_0})}\in \mathcal{O}_N^{\gamma}\left(3r/4\right)$ and remains within this domain after. We then define the following deterministic time
     \begin{equation}
     	\label{temps_deter_min}
     	T(r):=\underset{\vv{z_0}\in \partial\mathcal{O}_N^{\gamma}(r)}{\inf}T(\vv{z_0})>0,
     \end{equation}
    for which we have, for every  $t\in\left[0;T(r)\right]$, $\dfrac{3r}{4}\leq\mathcal{U}^{\gamma,N}_{\mathrm{tot}}(t)\leq r$. The positiveness of $T(r)$ comes from the fact that if it was not the case, it should exist a sequence $\left(\vv{x_0}^{n},\vv{v_0}^{n}\right)_n$ in $ \partial\mathcal{O}_N^{\gamma}(r)$ such that the sequence with general term $T_n(r):=T\left(\vv{x_0}^{n},\vv{v_0}^{n}\right)$ goes to $0$. The compactness of $\partial\mathcal{O}_N^{\gamma}(r)$ yields an extraction of a subsequence (denoted also by $\left(\vv{x_0}^{n},\vv{v_0}^{n}\right)_n$) converging toward a point  $\initias$ of $ \partial\mathcal{O}_N^{\gamma}(r)$. By continuity of the flow of $\grandevad{}$ we should have $r\leq \dfrac{3r}{4}$ which leads to a contradiction.
     Then we set 
     \begin{equation}
     	\mathcal{K}_r:=\{ (x_i,v_i)_{i\leq N}\in \reel^{Nd}\times \reel^{Nd},\,\dfrac{3r}{4}\leq\mathcal U^\gamma(\left[x_i,v_i\right]_{1\le i\leq N})
     	\leq r \}.
     \end{equation}
     Since $\partial\mathcal{O}_N^{\gamma}(R)\cap\mathcal{K}_r=\emptyset$, the minimal distance $\delta(r):=\inf\{\norme{\vv{z}-\vv{z}'},\,(\vv{z},\vv{z}')\in \partial\mathcal{O}_N^{\gamma}(R)\times \mathcal{K}_r\}$ is positive. We now choose $r$ sufficiently small such that $ \alpha(r):=\sqrt{\dfrac{r}{17}}+\sqrt{\dfrac{r}{17\gamma}}\leq \delta(r)$. It is easy to show that $\delta(r)$ converges to $\overline{\delta}_R:=d\left(\ensinvab{R},\pointeq\right)>0$ when $r$ goes to 0, so that, for $r$  small, $\delta(r)\geq\dfrac{\overline{\delta}_R}{2}\geq \alpha(r)$. We define the following events
     \begin{equation*}
     	C_\varepsilon^r=\{\underset{0\leq t\leq T(r)}{\max}\norme{\vv{X}^{N,\varepsilon}_t-\vv{x}^{N}(t)}<\sqrt{\dfrac{r}{17\gamma}}\}\quad\mathrm{and}\quad D_\varepsilon^r=\{\underset{0\leq t\leq T(r)}{\max}\norme{\vv{V}^{N,\varepsilon}_t-\vv{v}^{N}(t)}<\sqrt{\dfrac{r}{17}}\}.
     \end{equation*}
     On the event $C_\varepsilon^r\cap D_\varepsilon^r$,  $\underset{0\leq t\leq T(r)}{\max}\left(\norme{\vv{V}^{N,\varepsilon}_t-\vv{v}^{N}(t)}+\norme{\vv{X}^{N,\varepsilon}_t-\vv{x}^{N}(t)}\right)<\alpha(r)$.
     Still on $C_\varepsilon^r\cap D_\varepsilon^r$, let us assume that $\tau_{\mathcal{A}(r)}<T(r)$. By definition of $\tau_{\mathcal{A}(r)}$, two cases are then possible.\\
     \underline{First case}: $\grandeva{\tau_{\mathcal{A}(r)}}\in \partial\mathcal{O}_N^{\gamma}(R)$. If $\tau_{\mathcal{A}(r)}<T(r)$, $\grandevad{\tau_{\mathcal{A}(r)}}\in \mathcal{K}_r$ and so $\norme{\vv{X}^{N,\varepsilon}_{\tau_{\mathcal{A}(r)}}-\vv{x}^{N}(\tau_{\mathcal{A}(r)})} +\norme{\vv{V}^{N,\varepsilon}_{\tau_{\mathcal{A}(r)}}-\vv{v}^{N}(\tau_{\mathcal{A}(r)})}\geq\alpha(r)$ leads to a contradiction.\\
     \underline{Second case}: $\grandeva{\tau_{\mathcal{A}(r)}}\in \partial\mathcal{O}_N^{\gamma}(r/4)$ and so
     \begin{align*}
     	\gamma\norme{\vv{x}^{N}(\tau_{\mathcal{A}(r)})-\vv{x}_{\text{eq}}^{\ast}}^2+\norme{\vv{v}^{N}(\tau_{\mathcal{A}(r)})}^2&\leq 2\gamma\norme{\vv{X}^{N,\varepsilon}_{\tau_{\mathcal{A}(r)}}-\vv{x}^{N}(\tau_{\mathcal{A}(r)})}^2+2\gamma\norme{\vv{X}^{N,\varepsilon}_{\tau_{\mathcal{A}(r)}}-\vv{x}_{\text{eq}}}^2\\
     	&\quad+2\norme{\vv{V}^{N,\varepsilon}_{\tau_{\mathcal{A}(r)}}-\vv{v}^{N}(\tau_{\mathcal{A}(r)})}^2+2\norme{\vv{V}^{N,\varepsilon}_{\tau_{\mathcal{A}(r)}}}^2\\
     &<r/2+2r/17<3r/4.
     \end{align*}
    Hence $\grandevad{\tau_{\mathcal{A}(r)}}\in \mathcal{O}_N^{\gamma}(3r/4)$ and, by definition of $T(r)$, we have the contradiction $\tau_{\mathcal{A}(r)}<T(r)\leq \tau_{\mathcal{A}(r)}$. \\
    Therefore $C_\varepsilon^r\cap D_\varepsilon^r\subset \{\tau_{\mathcal{A}(r)}<T(r)\}$ and so, for every $\initia\in \partial\mathcal{O}_N^{\gamma}(r)$ the following lower-bound holds
     \begin{equation}
     	\probabi{(\tau_{\mathcal{A}(r)}\geq T(r))}\geq \probabi{(C_\varepsilon^r\cap D_\varepsilon^r)}\geq\probabi{(C_\varepsilon^r)}+\probabi{(D_\varepsilon^r)}-1.
     \end{equation}
     By virtue of \eqref{convunifo} in Theorem \ref{convL1}, one can find $\varepsilon_0>0$ such that for all $\initia\in \partial\mathcal{O}_N^{\gamma}(r)$ and for all $\varepsilon\!<\!\varepsilon_0$:
     \begin{equation*}
     	\probabi{\left(\tau_{\mathcal{A}(r)}\geq T(r)\right)}\geq\eta.
     \end{equation*}
 \end{proof}
 We state that when the noise is small in the particles system $\eqref{Cucker-Smale-bruite}$, the event of exiting $\ensinva{R}$ before hitting $\ensinvab{R}$ or $\ensinvab{r/4}$ has its probability close to $0$, the main contribution to this event is due to the attraction to the equilibrium ensuring with high probability the hitting of $\ensinvab{r/4}$.
 \begin{lemmea}
 \label{exit_b_ball}
There exists $r\in(0;R)$  such that, given $\vv{z_0}=\initia\in\ensinva{R}$ and $\zeta_1^\varepsilon$ as in \eqref{temps_osci}, we have 
\begin{equation*}
  \underset{\varepsilon\to 0}{\lim}\,\probabi{\left(\tau^{\alpha,\gamma,N}_{\varepsilon}(R)>\zeta_1^\varepsilon\right)}=1.  
\end{equation*}

 \end{lemmea}
 \begin{proof}
First, suppose that $\vv{z_0}\in \ensinvad{r}$ and consider:
 \begin{equation}
     \varrho_\varepsilon:=\inf\{n\in \mathbb{N},\, \vv{\mathcal Z}_n^{N,\varepsilon}\in\partial\mathcal{O}_N^{\gamma}(R)\},
 \end{equation}
 where $ \vv{\mathcal Z}_n^{N,\varepsilon}$ is the Markov chain \eqref{etat_osci}. We thus have the inclusion $\{\varrho_\varepsilon>1\}\subset \{\tau^{\alpha,\gamma,N}_{\varepsilon}(R)>\zeta_1^\varepsilon\}$. Indeed, on $\{\varrho_\varepsilon>1\} $, the events $\{\tau^{\alpha,\gamma,N}_{\varepsilon}(R)=\zeta_1^\varepsilon\}$ and $\{\tau^{\alpha,\gamma,N}_{\varepsilon}(R)<\zeta_1^\varepsilon\}$ are empty, since when $ \initia\in \ensinvad{r}$, a first exit from $\ensinva{R}$ happens after $\zeta_1^\varepsilon$. We may note that the inequality \eqref{mino_ten_sort} remains true when $\initia\in \ensinvad{r}$\footnote{This can be established by reconsidering the key starting point \eqref{sortie_v1} with $\initia\in \ensinvad{r}$ and adapting the rest of the arguments for \eqref{mino_ten_sort}; details are left aside for the sake of conciseness.} with $h=\frac{\mathcal{V}_N(R,\alpha,\gamma)}{2}$ and $n=1$, it follows 
 \begin{equation*}
     \underset{\varepsilon\to 0}{\lim}\,\probabi{\left(\tau^{\alpha,\gamma,N}_{\varepsilon}(R)>\zeta_1^\varepsilon\right)}=1.
 \end{equation*}
 Now, let us suppose $\initia\in \ensinva{R}\setminus\ensinvad{r}$. We define the deterministic time
 \begin{equation}
    \label{temps_0}
     T_0^\ast:=T\initia=\inf\{T\geq 0: \forall t\geq T ,(\vv{x}^N(t),\vv{v}^N(t))\in \ensinvad{r/4}\}>0.
 \end{equation}
The positiveness of $T_0^\ast$ can be obtained reasoning as for \eqref{temps_deter_min}. We now define the compact set
$$\mathcal{K}_0:=\{\grandevad{t}, t\in[0;T_0^\ast]\},$$
for which  the inclusion $ \mathcal{K}_0\subset \ensinva{R}$ so that $\mathcal{K}_0\cap \ensinvab{R}=\emptyset$. Let us take $\Delta_0:=\big(d(\mathcal{K}_0,\ensinvab{R})\big)^2>0$, and define $\Delta_r^\ast:=\min(r,\Delta_0)$ and the event 
 \begin{equation*}
     \mathcal{A}_N^\varepsilon:=\{\supre{T_0^\ast}\sqrt{\gamma\norme{\vv{X}^{N,\varepsilon}_t-\vv{x}^{N}(t)}^2+\norme{\vv{V}^{N,\varepsilon}_t-\vv{v}^{N}(t)}^2}\leq \dfrac{\min\left(1,\sqrt{\gamma}\right)\sqrt{\Delta_r^\ast}}{2}\}.
 \end{equation*}
In particular, observing that, on  $ \mathcal{A}_N^\varepsilon$, 
  \begin{align*}
        \gamma\norme{\vv{X}^{N,\varepsilon}_{T_0^\ast}-\vv{x}_{\mathrm{eq}}^\ast}^2+\norme{\vv{V}^{N,\varepsilon}_{T_0^\ast}}^2&\leq 2\gamma\Big(\norme{\vv{X}^{N,\varepsilon}_{T_0^\ast}-\vv{x}^{N}(T_0^\ast)}^2+\norme{\vv{x}^{N}(T_0^\ast)-\vv{x}_{\mathrm{eq}}^\ast}^2\Big)\\
        &\quad+2\Big(\norme{\vv{V}^{N,\varepsilon}_{T_0^\ast}-\vv{v}^{N}(T_0^\ast)}^2+\norme{\vv{v}^{N}(T_0^\ast)}^2\Big)\!\!\\
        &\leq \dfrac{\Delta_r^\ast}{2}+\dfrac{r}{2}\leq r,
 \end{align*}
 necessarily
 \begin{equation}
    \label{inc_ev}
    \mathcal{A}_N^\varepsilon\subset \{\grandeva{T_0^\ast}\in \ensinvad{r}\}. 
 \end{equation}
 Observe also the equality
 \begin{align}
    \label{event}
     \{\tau^{\alpha,\gamma,N}_{\varepsilon}(R)\neq\zeta_1^\varepsilon\}\cap\mathcal{A}_N^\varepsilon=\left(\{\tau^{\alpha,\gamma,N}_{\varepsilon}(R)>\zeta_1^\varepsilon\}\cap\mathcal{A}_N^\varepsilon\right)\bigcup\mathcal{B}_\varepsilon(r),
 \end{align}
 with $\mathcal{B}_\varepsilon(r):=\{\tau^{\alpha,\gamma,N}_{\varepsilon}(R)<\zeta_1^\varepsilon\}\cap\mathcal{A}_N^\varepsilon$. We split $\mathcal{B}_\varepsilon(r)$ by $\mathcal{B}_\varepsilon(r)=\Omega_N^{\varepsilon,1}(r)\cup\Omega_N^{\varepsilon,2}(r)$ where $\Omega_N^{\varepsilon,1}$ and $\Omega_N^{\varepsilon,2}$ are
 \begin{equation*}
     \Omega_N^{\varepsilon,1}(r):=\mathcal{A}_N^\varepsilon\bigcap \{T_0^\ast<\tau^{\alpha,\gamma,N}_{\varepsilon}(R)<\zeta_1^\varepsilon\}\quad\mathrm{and}\quad\Omega_N^{\varepsilon,2}(r):=\mathcal{A}_N^\varepsilon\bigcap\{\tau^{\alpha,\gamma,N}_{\varepsilon}(R)<\zeta_1^\varepsilon\}\bigcap \{\tau^{\alpha,\gamma,N}_{\varepsilon}(R)\leq T_0^\ast\}.
 \end{equation*}
 The event $\Omega_N^{\varepsilon,1}(r)$ is indeed empty because $\grandeva{T_0^\ast}\in \ensinvad{r}$ thanks to \eqref{inc_ev}; the noisy system reaches $\ensinvab{R}$ the first time from $\zeta_1^\varepsilon$ onward.\\
 Similarly for $\omega\in\Omega_N^{\varepsilon,2}(r)$, on $\mathcal{A}_N^\varepsilon$, we have, on one hand,
 \begin{align*}
     \sqrt{\gamma\norme{\vv{X}^{N,\varepsilon}_{\tau^{\alpha,\gamma,N}_{\varepsilon}(R)}-\vv{x}^{N}(\tau^{\alpha,\gamma,N}_{\varepsilon}(R))}^2+\norme{\vv{V}^{N,\varepsilon}_{\tau^{\alpha,\gamma,N}_{\varepsilon}(R)}-\vv{v}^{N}(\tau^{\alpha,\gamma,N}_{\varepsilon}(R))}^2}\leq \dfrac{\min\left(1,\sqrt{\gamma}\right)\sqrt{\Delta_r^\ast}}{2},
 \end{align*}
 and, on the other hand, $\grandeva{\tau^{\alpha,\gamma,N}_{\varepsilon}(R)}\in \ensinvab{R}$ and $(\vv{x}^{N}(\tau^{\alpha,\gamma,N}_{\varepsilon}(R)),\vv{v}^{N}(\tau^{\alpha,\gamma,N}_{\varepsilon}(R)))\in \mathcal{K}_0$. By definition of $\Delta_r^\ast $, it follows that
 \begin{align*}
     \sqrt{\Delta_r^\ast}\leq d\left(\mathcal{K}_0,\ensinvab{R}\right)\leq d\left(\grandeva{\tau^{\alpha,\gamma,N}_{\varepsilon}(R)},\grandevad{\tau^{\alpha,\gamma,N}_{\varepsilon}(R)}\right)\leq \dfrac{\sqrt{\Delta_r^\ast}}{2},
 \end{align*}
which leads to a contradiction. This yields $\Omega_N^{\varepsilon,2}(r)=\emptyset$.\\
Hence \eqref{inc_ev} and  \eqref{event} give the inclusion $ \{\tau^{\alpha,\gamma,N}_{\varepsilon}(R)\neq\zeta_1^\varepsilon\}\cap\mathcal{A}_N^\varepsilon\subset \{\tau^{\alpha,\gamma,N}_{\varepsilon}(R)>\zeta_1^\varepsilon\}$. 
Taking $h=\frac{\mathcal{V}_N(R,\alpha,\gamma)}{2}$, according to Lemma~\ref{ev_impossi}, there exist $\varepsilon_0,r>0$ such that for any $\initia\in \ensinva{R}\setminus\ensinvad{r}$ and for all $\varepsilon<\varepsilon_0$ we have
 \begin{align}
     \probabi{\left(\tau^{\alpha,\gamma,N}_{\varepsilon}(R)>\zeta_1^\varepsilon\right)}&\geq \probabi{\left(\{\tau^{\alpha,\gamma,N}_{\varepsilon}(R)\neq\zeta_1^\varepsilon\}\cap\mathcal{A}_N^\varepsilon\right)}\geq \probabi{\left(\tau^{\alpha,\gamma,N}_{\varepsilon}(R)\neq\zeta_1^\varepsilon\right)}+\probabi{\left(\mathcal{A}_N^\varepsilon\right)}-1\notag\\
     &\geq \probabi{\left(\mathcal{A}_N^\varepsilon\right)}-\exp{\left(-\dfrac{\mathcal{V}_N(R,\alpha,\gamma)}{2\varepsilon^2}\right)}.
 \end{align}
 Since, according to Theorem \ref{convL1}, $\displaystyle\lim_{\varepsilon\rightarrow 0}\probabi{\left(\mathcal{A}_N^\varepsilon\right)}=1$, the claim follows.
 
 \end{proof}
 We close this section with a rather classical Gaussian estimate on the proximity between  $\grandeva{}$ and its starting point $\initia\in \ensinva{R}$, during a short time. 
 \begin{lemmea}
 \label{proche_ini}
 For all $\delta,c>0$, there exists a deterministic time $T(\delta,c)>0$, such that for all $\vv{z_0}=\initia\in \ensinva{R}$, $T\leq T(\delta,c)$ and for all $\varepsilon>0$, we have
 \begin{align*}
     \max\Big(\probabi{\left(\supre{T}\norme{\vv{X}^{N,\varepsilon}_t-\vv{x_0}}\geq \delta\right)},\probabi{\left(\supre{T}\norme{\vv{V}^{N,\varepsilon}_t-\vv{v_0}}\geq \delta\right)}\Big)\leq 4Nd\exp{\left(-\dfrac{c}{\varepsilon^2}\right)}.
 \end{align*}
  \end{lemmea}
 \begin{proof}
Let us take $\initia\in \ensinva{R}$. We rewrite the system \eqref{Cucker-Smale-bruite} as
\begin{equation*}
	\deriv\vv{X}^{N,\varepsilon}_t=\vv{V}^{N,\varepsilon}_t\deriv t,\qquad\vv{V}^{N,\varepsilon}_t=b^{N}\left(\grandeva{t}\right)\deriv t -\alpha\vv{V}^{N,\varepsilon}_t\deriv t-\gamma\left(\vv{X}^{N,\varepsilon}_t-\vv{x}_{\mathrm{eq}}^\ast\right)\di t+\varepsilon \deriv W_t,
\end{equation*}
with $(\vv{X}^{N,\varepsilon}_0,\vv{V}^{N,\varepsilon}_0)=(\vv{x_0},\vv{v_0})$ and $ b^N\left(\left(x_j,v_j\right)_{j \in \ens}\right)=\big(b^{i,N}\left(\left(x_j,v_j\right)_{j \in \ens}\right)\big)_{i=1}^N$ for $b^{i,N}$ as in \eqref{fonction_communis}. For any $0\le s\le t\le T$, the triangle inequality gives $(\widehat{\phi}=2\lambda\sqrt{N}/\kappa^\beta)$
\begin{align*}
    \underset{0\leq s\leq t}{\max}\norme{\vv{V}^{N,\varepsilon}_s-\vv{v_0}}&\leq \inte{t}\norme{b^{N}\left(\grandeva{s}\right)}\di s +\alpha \inte{t}\norme{\vv{V}^{N,\varepsilon}_s}\di s\\
    &\quad+\gamma\inte{t}\norme{\vv{X}^{N,\varepsilon}_t-\vv{x}_{\mathrm{eq}}^\ast}\di s +\varepsilon\supre{T}\norme{W_t}\\
    &\leq \left(\widehat{\phi}+\alpha\right)\inte{t}\left(\norme{\vv{V}^{N,\varepsilon}_s-\vv{v_0}}+\norme{\vv{v_0}}\right)\di s + \varepsilon\supre{T}\norme{W_t}\\
    &\quad+\gamma\inte{t}\left(\norme{\vv{X}^{N,\varepsilon}_s-\vv{x_0}}+\norme{\vv{x}_{\mathrm{eq}}^\ast-\vv{x_0}}\right)\di s,
\end{align*}
 and where
\begin{equation}\label{CroissancePosition}
\norme{\vv{X}^{N,\varepsilon}_s-\vv{x_0}}\leq T  \underset{0\leq r\leq s}{\max}\norme{\vv{V}^{N,\varepsilon}_r-\vv{v_0}}+ T\|\vv{v}_0\|.
\end{equation}
 Recalling that $\initia\in \ensinva{R}$ $\Leftrightarrow$ $\|\vv{v_0}\|^2+\gamma\norme{\vv{x}_{\mathrm{eq}}^\ast-\vv{x_0}}^2\le R$, it follows that
\begin{align*}
    \underset{0\leq s\leq t}{\max}\norme{\vv{V}^{N,\varepsilon}_t-\vv{v_0}}&\leq \left[ \left(\widehat{\phi}+\alpha\right)+\sqrt{\gamma}\right]T\sqrt{R} +T^2\gamma\sqrt{R}+\varepsilon\supre{T}\norme{W_t}\\
    &\quad+\left(\widehat{\phi}+\alpha+\gamma T\right)\inte{t}\underset{0\leq r\leq s}{\max}\norme{\vv{V}^{N,\varepsilon}_r-\vv{v_0}}\di s.
\end{align*}
From Grönwall lemma, we obtain the upper-bound
\begin{align*}
    \underset{0\leq t\leq T}{\max}\norme{\vv{V}^{N,\varepsilon}_t-\vv{v_0}}&\leq  \left(g_R(T)+\varepsilon\supre{T}\norme{W_t}\right)h(T),
\end{align*}
where $\displaystyle g_R(T)\!\!:=\left[ \left(\widehat{\phi}+\alpha\right)+\sqrt{\gamma}\right]T\sqrt{R} +T^2\gamma\sqrt{R} $ and $h(T):=\exp{\left(\left(\widehat{\phi}+\alpha+\gamma T\right)T\right)}$. From \eqref{CroissancePosition}, we also~get
\begin{align*}
    \underset{0\leq s\leq t}{\max}\norme{\vv{X}^{N,\varepsilon}_s-\vv{x_0}}&\leq  T\left(g_R(T)+\varepsilon\supre{T}\norme{W_t}\right)h(T)+ T\sqrt{R}.
\end{align*}
Since we have $\underset{T\to 0}{\lim}\, g_R(T)=0$, one can find $T_1:=T_1(\delta,R)\in(0;1)$ such that taking $T\leq T_1$ implies $$\max\left(g_R(T)h(T),T g_R(T)h(T)+ T\sqrt{R}\right)<\dfrac{\delta}{2},$$ 
and so the inclusion: 
\begin{equation*}
 \{\supre{T}\norme{\vv{V}^{N,\varepsilon}_t-\vv{v_0}}\geq \delta\}\cup \{\supre{T}\norme{\vv{X}^{N,\varepsilon}_t-\vv{x_0}}\geq \delta\}\subset\{\supre{T}\norme{W_t}\geq \dfrac{\delta}{2\varepsilon h(T)}\}.
\end{equation*}
In addition, for all $\mu>0$ and $0<T<\infty$, we have, by Doob's and (an exponential) Markov's inequalities,
\begin{align*}
    \proba{\supre{T}\norme{W_t}\geq \mu}
    \leq 4Nd\exp{\left(-\dfrac{\mu^2}{2(Nd)^2T}\right)}.
\end{align*}
Hence, for all $T\leq T_1(\delta,R) $, we get
 \begin{align*}
     \probabi{\left(\supre{T}\norme{\vv{X}^{N,\varepsilon}_t-\vv{x_0}}\geq \delta\right)}\vee\probabi{\left(\supre{T}\norme{\vv{V}^{N,\varepsilon}_t-\vv{v_0}}\geq \delta\right)}\leq 4Nd\exp{\left(-\dfrac{\delta^2}{8(Nd)^2\varepsilon^2T h^2(T)}\right)}.
 \end{align*}
 Since  $h$ is increasing, the claim follows by reducing $T_1(\delta,R)$ into $ T_1(\delta,c)$ so that $\left(8(Nd)^2T_1 h^2(T_1)\right)^{-1}\geq \dfrac{c}{\delta^2}$. 
 \end{proof}
\bibliographystyle{plain} 
\bibliography{reference.bib}

\end{document}